\documentclass[12pt]{article}
\usepackage{amsfonts,amsmath,amsxtra}
\usepackage{latexsym,textcomp}
\usepackage{amssymb}
\usepackage{color}
\usepackage[matrix,arrow,curve]{xy}

\makeatletter
\def\@seccntformat#1{\csname
the#1\endcsname\enspace}
\makeatother

\def\hybrid{\topmargin 0pt      \oddsidemargin 0pt
        \headheight 0pt \headsep 0pt
        \textwidth 16.5cm
        \textheight 23cm
        \marginparwidth 0.0in
        \parskip 5pt plus 1pt   \jot = 1.5ex}
\catcode`\@=11
\def\marginnote#1{}
\newcount\hour
\newcount\minute
\newtoks\amorpm
\hour=\time\divide\hour by60 \minute=\time{\multiply\hour by60
\global\advance\minute by-\hour}
\edef\standardtime{{\ifnum\hour<12 \global\amorpm={am}%
        \else\global\amorpm={pm}\advance\hour by-12 \fi
        \ifnum\hour=0 \hour=12 \fi
      \number\hour:\ifnum\minute<10 0\fi\number\minute\the\amorpm}}
\edef\militarytime{\number\hour:\ifnum\minute<10 0\fi\number\minute}

\def\draftlabel#1{{\@bsphack\if@filesw {\let\thepage\relax
   \xdef\@gtempa{\write\@auxout{\string
      \newlabel{#1}{{\@currentlabel}{\thepage}}}}}\@gtempa
   \if@nobreak \ifvmode\nobreak\fi\fi\fi\@esphack}
        \gdef\@eqnlabel{#1}}
\def\@eqnlabel{}
\def\@vacuum{}
\def\draftmarginnote#1{\marginpar{\raggedright\scriptsize\tt#1}}

\def\draft{\oddsidemargin -0.1truein
        \def\@oddfoot{\sl GivGT.tex \hfil 
        \rm\thepage\hfil\sl\today\quad\militarytime}
        \let\@evenfoot\@oddfoot \overfullrule 3pt
        \let\label=\draftlabel
        \let\marginnote=\draftmarginnote
\def\@eqnnum{{\rm (\theequation)}
\rlap{\kern\marginparsep\tt\@eqnlabel}%
\global\let\@eqnlabel\@vacuum}  }

\newfont{\Bbbb}{msbm7 scaled 1\@ptsize00}
\newcommand{\zs}{\raise-1pt\hbox{$\mbox{\Bbbb Z}$}}

scaled 1\@ptsize00

\def\numberbysection{\@addtoreset{equation}{section}
        \def\theequation{\thesection.\arabic{equation}}}
\numberbysection

\renewcommand{\theequation}{\thesection.\arabic{equation}}

\def\titlepage{\@restonecolfalse\if@twocolumn\@restonecoltrue\onecolumn
     \else \newpage \fi \thispagestyle{empty}\c@page\z@
\def\thefootnote{\fnsymbol{footnote}} }
\def\endtitlepage{\if@restonecol\twocolumn \else  \fi
        \def\thefootnote{\arabic{footnote}}
        \setcounter{footnote}{0}}  
\relax \hybrid
\makeatletter
\newdimen\normalarrayskip            
\newdimen\minarrayskip               
\normalarrayskip\baselineskip \minarrayskip\jot
\newif\ifold             \oldtrue            \def\new{\oldfalse}
\def\arraymode{\ifold\relax\else\displaystyle\fi}
\def\eqnumphantom{\phantom{(\theequation)}} 
\def\@arrayskip{\ifold\baselineskip\z@\lineskip\z@
     \else
     \baselineskip\minarrayskip\lineskip1\baselineskip\fi}
\def\@arrayclassz{\ifcase \@lastchclass \@acolampacol \or
\@ampacol \or \or \or \@addamp \or
   \@acolampacol \or \@firstampfalse \@acol \fi
\edef\@preamble{\@preamble
  \ifcase \@chnum
     \hfil$\relax\arraymode\@sharp$\hfil
     \or $\relax\arraymode\@sharp$\hfil
     \or \hfil$\relax\arraymode\@sharp$\fi}}
\def\@array[#1]#2{\setbox\@arstrutbox=\hbox{\vrule
     height\arraystretch \ht\strutbox
     depth\arraystretch \dp\strutbox
width\z@}\@mkpream{#2}\edef\@preamble{\halign \noexpand\@halignto
\bgroup \tabskip\z@ \@arstrut \@preamble \tabskip\z@ \cr}%
\let\@startpbox\@@startpbox \let\@endpbox\@@endpbox
  \if #1t\vtop \else \if#1b\vbox \else \vcenter \fi\fi
  \bgroup \let\par\relax
  \let\@sharp##\let\protect\relax
  \@arrayskip\@preamble}
\def\eqnarray{\stepcounter{equation}%
              \let\@currentlabel=\theequation
              \global\@eqnswtrue
              \global\@eqcnt\z@
              \tabskip\@centering              
              \let\\=\@eqncr
              $$%
            \halign to \displaywidth  \bgroup
             \eqnumphantom \@eqnsel
      \hskip\@centering                               
    $\displaystyle  \tabskip\z@ {##}$%
    &\global\@eqcnt\@ne \hskip 2\arraycolsep
         $ \displaystyle  \arraymode{##}$\hfil
    &\global\@eqcnt\tw@ \hskip 2\arraycolsep
         $\displaystyle\tabskip\z@{##}$\hfil
         \tabskip\@centering
    &{##}\tabskip\z@\cr}
\makeatother

\def\IC{\mathbb{C}}

\def\IR{\mathbb{R}}

\def\CB {\mathcal{B}}
\def\CC {\mathcal{C}}
\def\CD {\mathcal{D}}
\def\CE {\mathcal{E}}
\def\CF {\mathcal{F}}

\def\CH {\mathcal{H}}

\def\CN {\mathcal{N}}

\def\CS {\mathcal{S}}

\def\CU {\mathcal{U}}
\def\CV {\mathcal{V}}

\def\CZ {\mathcal{Z}}
\def\fg{{\frak g}}

\def\fh{{\frak h}}

\def\fn{{\frak n}}

\def\a {{\alpha}}

\def\KE{{{}^{\vk}\!\wt{E}}}
\def\g {{\gamma}}
\def\s {{\sigma}}

\def\la{\lambda}
\def\ve{\varepsilon}
\def\vk{\varkappa}

\def\pr {\partial}

\def\wt{\widetilde}
\def\wh{\widehat}

\def\Arg{{\mathop{\rm Arg}}}

\def\frak{\mathfrak}
\def\ov {{\overline}}

\def\gl{\mathfrak{gl}}

\def\<{\langle}
\def\>{\rangle}
\def\ov{\overline}
\newtheorem{te}{Theorem}[section]

\newtheorem{prop}{Proposition}[section]           

\newtheorem{lem}{Lemma}[section]

\newcommand\rk{\operatorname{rank}}

\newcommand{\proof}{\noindent {\it Proof}. }
\newcommand\bqa{\begin{eqnarray}}
\newcommand\eqa{\end{eqnarray}}
\def\be{\begin{eqnarray}\new\begin{array}{cc}}
\def\ee{\end{array}\end{eqnarray}}
\def\beq{\begin{equation}}
\def\eeq{\end{equation}}
\def\bse{\begin{subequations}}                
\def\ese{\end{subequations}}
\def\bp{\begin{pmatrix}}
\def\ep{\end{pmatrix}}

\newcounter{pac}[section]

\newcounter{pacc}[subsection]

 0

\begin{document}

\title{\bf On equivalence of the Mellin-Barnes and\\ the  Givental integral
  representations of\\ the Whittaker functions}
\author{A.A. Gerasimov, D.R. Lebedev and S.V. Oblezin}
\date{\today}
\maketitle

\renewcommand{\abstractname}{}

\begin{abstract}
\noindent {\bf Abstract}.  We construct an integral transformation
intertwining the Gelfand-Tsetlin and the (modified) Gauss-Givental
realizations of principle series representations of
$\mathfrak{gl}_3(\IR)$. This provides a direct identification of the
corresponding integral representations for the
$\mathfrak{gl}_3(\IR)$-Whittaker function.  The construction
essentially uses integral identities due to Barnes and Gustafson
thus providing a basis for their representation theory
interpretation. The result of this paper might be useful for
constructing the explicit analytic realization of the mirror
symmetry map in the case of the flag manifold $GL_3(\IC)/B$.
\end{abstract}
\vspace{5mm}

\section{Introduction}

Analytical objects like special functions as well as  various
integral identities have been a continuing source of inspiration for
finding their hidden structures, explaining their (sometimes
mysterious) properties. Along this route, going from analysis to
algebra, we usually encounter new structures and understand better
the known ones. One of conceptual frameworks proved to be very
successful in uncovering hidden structures of special functions and
associated integral identities is provided by representation theory
of Lie groups and Lie algebras (see for example \cite{V}), as well
as by their higher dimensional analogs given by mapping spaces into
the Lie groups. Often the analytical objects are arranged in series
resembling the series of classical Lie groups, which strongly
suggests representation theory background. Following this approach,
we are looking for interpretation of special functions as matrix
elements of representations of algebraic structures, so that
associated integral identities reflect basic structures of their
representation theory. In some cases this line of reasoning leads to
constructions of various generalizations/deformations  of Lie
groups/algebras. General context for these considerations is
 non-commutative geometry modeled via non-commutative harmonic
 analysis on Lie groups.

First example of a class of (rather mysterious) integral identities
was found by Barnes, now known as the First Barnes Lemma \cite{Ba}:
 \be\label{Barnes}
  \frac{1}{2\pi}\!\int\limits_{\IR}\!d\g\,\prod\limits_{i=1}^2\Gamma(a_i-\imath\g)
  \prod\limits_{j=1}^2\Gamma(b_j+\imath\g)\,
  =\,\frac{\prod\limits_{i,j=1}^2\Gamma(a_j+b_i)}
  {\Gamma\Big(\sum\limits_{i=1}^2a_i+b_i\Big)}\,,\quad
  {\rm Re}(a_j,b_i)>0\,.
  \ee
A bit later, identities of similar type were reported by Ramanujan
in his last paper published before moving to England \cite{Ra}. This
line of development culminated in several series of identities found
by Gustafson \cite{Gu}, in particular, (formula (9.4), for $n=1$)
 \be\label{Gustafson}
  \frac{1}{2\pi}\!\int\limits_{\IR}\!d\gamma\,
  \frac{\prod\limits_{i=1}^4\Gamma(a_i+\imath \gamma)\,
  \Gamma(a_i-\imath \gamma)}{\Gamma(2\imath \gamma)\,\Gamma(-2\imath
  \gamma)}\,
  =\,\frac{2\prod\limits_{i<j}
  \Gamma(a_i+a_j)}{\Gamma(a_1+a_2+a_3+a_4)}\,,\quad{\rm Re}(a_i)>0\,.
 \ee
The Gustafson identities arise in series, which makes a relation
with representation theory of classical Lie groups quite plausible.
However  to our best knowledge up to now, no a representation theory
interpretation of these identities is found, yet.

In this paper, we report on first steps towards  such an
interpretation that has unexpectedly arisen in a problem of
identifying certain matrix elements of principle series
representations of the general linear Lie group in various
realizations. The considered realizations of representations of Lie
algebras by operators in some particular spaces are not arbitrary
but have deep meaning in terms of the Gelfand-Tsetlin construction
of bases in representations \cite{GT}, positivity structures, as
well as torifications of flag manifolds. Hence one should expect
that thus arising identities reveal deep constructions in
representation theory. Precisely, the identities of this class arise
in the setting of
  this paper as follows. We consider the $\mathfrak{gl}_{\ell+1}(\IR)$-Whittaker
function given by common eigenfunction of the
$\mathfrak{gl}_{\ell+1}(\IR)$-Toda chain Hamiltonians, well-known
finite-dimensional (quantum) integrable system (see e.g. \cite{STS}
and references therein). This entails its identification as a
special matrix element of principle series representation
$(\pi_{\mu},\,\CV_{\mu})$ of the universal enveloping algebra
$\CU(\mathfrak{gl}_{\ell+1})$ of the Lie algebra
$\mathfrak{gl}_{\ell+1}(\IR)$. The representation space $\CV_{\mu}$
allows various realizations as function spaces thus leading to
various explicit realizations of
$\CU(\mathfrak{gl}_{\ell+1})$-representations via
differential/difference operators. Two such functional realizations
having special properties were considered in \cite{GKL},\cite{GKLO}.
One realization relies on the recursive (over the rank) construction
of maximal commuting subalgebras in $\CU(\mathfrak{gl}_{\ell+1})$
and generalizes the Gelfand-Tsetlin construction of bases in
finite-dimensional representations of $\gl_{\ell+1}(\IC)$. This
realization has deep connections with the separation of variables
method for quantum integrable systems \cite{Skl} and the associated
Hecke-Baxter operators \cite{GLO08}. The other one is the
Gauss-Givental functional realization, which clarifies the works of
Givental on mirror symmetry for flag manifolds and his
 integral representations of Whittaker functions \cite{Giv}.
Along with applications to  mirror
 symmetry issues, this construction has deep relations with
 the positivity structures in flag manifolds.  Both considered
functional realizations of principle series representations lead to
corresponding integral representations for the
$\mathfrak{gl}_{\ell+1}(\IR)$-Whittaker functions. In this paper, we
are interested in constructing the
 operators intertwining  these realizations thus allowing a direct
transformation  of one integral representation into another. Our
main result  is the construction of such operators for the first
non-trivial case of  $\ell=2$  (the quite trivial  case of $\ell=1$
is also considered). It is employing the intertwining operator to
identify the two realizations, we use the Barnes identity
\eqref{Barnes} and (a degenerate form of) the Gustafson identity
\eqref{Gustafson} in essential way. We also construct a direct
integral transformation between the considered  integral
representations of the $\mathfrak{gl}_3(\IR)$-Whittaker function
without reference to its matrix element structure. Here again the
identity \eqref{Gustafson} shows up inevitably. Although we consider
only the first non-trivial case of $\mathfrak{gl}_3(\IR)$, there are
clear hints of a possibility of a recursive (over the rank $\ell$)
construction of intertwining operators for
$\mathfrak{gl}_{\ell+1}(\IR)$. It would be interesting to single out
thus arising integral identities.

There is a different perspective on the problem of explicit
identification of the two types of  integral representations for the
$\gl_{\ell+1}(\IR)$-Whittaker functions (which was our primal
  goal to start these investigations).  Actually, such an identification
provides an important ingredient  of  explicit analytic mirror
map construction  for flag manifolds. This line of research  was started by
Givental \cite{Giv} who constructed the complete set of linear
differential equations satisfied by the generating function in
certain topological quantum field theory associated with flag
manifolds for $GL_{\ell+1}(\IC)$. Namely, he identified these linear
differential operators with generators of the ring of quantum
Hamiltonians of the $\mathfrak{gl}_{\ell+1}(\IR)$-Toda chain and
provided an integral representation for their common
  eigenfunctions. In \cite{GKLO}, the Givental integral formula of
the $\gl_{\ell+1}(\IR)$-Whittaker functions was rederived using
representation theory interpretation of the Toda chain via special
choice of realization of the underlying
$\CU(\gl_{\ell+1})$-representation. Another type of integral
representation for the $\gl_{\ell+1}(\IR)$-Whittaker function is
based on the Gelfand-Tsetlin realization, and it also has a proper
meaning in the topological quantum field theory framework. The
existence of these two forms of integral representations reflects
the dual way to interpret Whittaker functions via mirror symmetry of
the underlying topological quantum field theories associated with
flag manifolds \cite{GLO11}.

In the case of toric varieties, the mirror symmetry map may be
constructed quite explicitly leading to two different types of
expressions for (deformed) generating function in the relevant
two-dimensional quantum field theories. Since in general flag
manifolds are not toric, one should first construct a toric
degeneration of the flag manifold and then apply mirror symmetry map
to the resulting toric manifold.  The approach of
Strominger-Yau-Zaslow to mirror symmetry \cite{SYZ}  provides a
conceptual explanation for such considerations. The relevant
analytic structures are as follows. We start with the Givental
integral representation which relies  on a torification of the flag
manifold $GL_{\ell+1}(\IC)/B$. Then the standard Fourier transform
applied to representation space gives rise to the mirror dual
integral representation given by the GKZ type integral (see
\cite{GG} and references therein). For the representation theoretic
account of the GKZ integrals and the relations with the Givental
integral formulas for the $\gl_{\ell+1}(\IR)$-Whittaker functions,
see \cite{GLO23}.  Now, additional transformation (constructed
in this paper for the case of $\mathfrak{gl}_3(\IR)$) identifies
thus obtained integral representation with the Gelfand-Tsetlin one.
This establishes a connection of the torified flag manifolds with
representation theory constructions after Gelfand-Tsetlin. As a
consequence, we thus provide  an  important link between mirror
symmetry and representation theory of Yangians \cite{GKL}.

Finally, let us point out that mirror symmetry in the  context of
flag manifolds might be also considered as a manifestation of the
Langlands correspondence for Archimedean fields \cite{GLO11}. It is
important to stress that here we are dealing with the analytical
form of the Langlands correspondence which not only identifies  the
sets of parameters of the underlying representations but essentially
works at the level of analytic objects, like $L$-functions and
the Whittaker functions.

The plan  of the paper is as follows. In Section 2 we start with
 recalling briefly realizations of the
$\mathfrak{gl}_{\ell+1}(\IR)$-Whittaker functions as special matrix
elements in principle series representations of
$\CU(\mathfrak{gl}_{\ell+1})$. We describe explicitly three
functional realizations of these representations: the
Gelfand-Tsetlin, the Gauss-Givental and the modified Gauss-Givental
realizations. These realizations lead to  corresponding integral
expressions for the $\mathfrak{gl}_{\ell+1}(\IR)$-Whittaker
functions. In the following two sections we construct integral
operators that intertwine different realizations of representations
of $\CU(\mathfrak{gl}_{\ell+1})$, for $\ell=1,2$.  In particular
application of the intertwining operators  to Whittaker vectors
provides an identification of integral expressions of the matrix
elements in the considered realizations. This is done in Section 3
for $\mathfrak{gl}_2(\IR)$ and in Section 4 for
$\mathfrak{gl}_3(\IR)$ (Theorem \ref{MAINTHM}). The latter case
involves applications of the identity \eqref{Barnes} and a
degenerate version of \eqref{Gustafson}. Finally, in Section 5 a
direct identification of the two integral representations for the
$\gl_2(\IR)$- and $\gl_3(\IR)$- Whittaker functions considered in
Section 2 via integral transform is given. The key step is done via
invoking the Gustafson identity \eqref{Gustafson} from \cite{Gu}.
Details on the modified Gauss-Givental realization obtained via the
Fourier transform  of the original Gauss-Givental realization are
given in Appendix A, while the main part of calculations for the
$\mathfrak{gl}_3(\IR)$-intertwining operator  is collected in
Appendix B.

\section{Matrix element realizations of the $\gl_{\ell+1}(\IR)$-Whittaker function}

In this section we recall matrix element representation for the
$\gl_{\ell+1}(\IR)$-Whittaker function. This representation is
instrumental for constructing various integral representations of
the $\gl_{\ell+1}(\IR)$-Whittaker function.

Recall from \cite{Ko} the interpretation of Whittaker functions as
particular matrix elements of infinite-dimensional representations
 of universal enveloping algebras of reductive Lie algebras.
For a real reductive Lie algebra $\fg$ with a Cartan subalgebra
$\fh\subset\fg$, let $\CU(\fg)$ be the universal enveloping algebra
of $\fg$. Let $\Phi(\fg)$ be the root system of $\fh\subset\fg$,
$\Phi_+(\fg)\subset\Phi(\fg)$ be the subset of positive roots, let
$\rho=\frac{1}{2}\sum\limits_{\a\in\Phi_+(\fg)}\a$, and let
$\Pi(\fg)\subset\Phi(\fg)$ be the set of simple roots. For
$\g\in\IC^{\rk(\fg)}$, given a principal series representation
$(\pi_{-\imath\g-\rho},\CV_{-\imath\g-\rho})$ of $\fg$, we assume
that the $\fh$-action integrates to the action of the torus
$A=\exp(\mathfrak{h})$, so that
$(\pi_{-\imath\g-\rho},\CV_{-\imath\g-\rho})$ affords a structure of
$(\fg,A)$-module. Let
$(\pi_{-\imath\g-\rho},\,\CV_{-\imath\g-\rho})$ and
$(\pi^{\vee}_{-\imath\g-\rho},\,\CV^{\vee}_{-\imath\g-\rho})$ be a
pair of dual $\CU(\fg)$-modules with a sesquilinear pairing (anti
$\IC$-linear in the first argument and $\IC$-linear in the second
argument), for $\phi_1\in\CV_{-\imath\g-\rho}^{\vee}$ and
$\phi_2\in\CV_{-\imath\g-\rho}$:
  \be\label{pair}
  \<\,,\,\>\,:\quad\CV_{-\imath\g-\rho}^{\vee}
  \times
  \CV_{-\imath\g-\rho}\,\longrightarrow\,\IC\,\,,\\
  \<\pi^{\vee}_{-\imath\g-\rho}(X)\phi_1,\,\phi_2\>\,
  =\,-\<\phi_1,\,\pi_{-\imath\g-\rho}(X)\phi_2\>\,,\qquad\forall
  X\in\fg\,.
 \ee
 The dual $\CU(\fg)$-module allows the following
identification:
 \be\label{dualrep}
  \CV^{\vee}_{-\imath\g-\rho}\,\simeq\,\CV_{-\imath\bar{\g}-\rho}\,.
 \ee
In the case $\g=\la\in\IR^{\rk(\fg)}$, \eqref{pair} defines a
non-degenerate $\fg$-invariant Hermitian pairing  on
$\CV_{-\imath\la-\rho}$, thus providing a structure of unitary
representation on it.

Now, let $\fn_{\pm}\subset\fg$ be the pair of opposite nilpotent
subalgebras corresponding to the choice of positive roots
$\Phi_+(\fg)$. Introduce the pair of Whittaker vectors
$w_L\in\CV^{\vee}_{-\imath\g-\rho}$ and $w_R\in\CV_{-\imath\g-\rho}$
to be non-degenerate characters of the nilpotent subalgebras
$\fn_-\subset\fg$ and $\fn_+\subset\fg$, respectively. Explicitly,
in terms of the Chevalley-Weyl generators
$e_{\alpha_j},\,f_{\alpha_j},\,\alpha_j\in\Pi(\fg)$ of
$\fn_{\pm}\subset\fg$, the Whittaker vectors $w_L$, $w_R$ are
defined as solutions of the following equations:
 \be\label{whit}
  \pi^{\vee}_{-\imath\g-\rho}(f_{\alpha_j})\cdot w_L\,=\,-w_L\,,\qquad
  \pi_{-\imath\g-\rho}(e_{\alpha_j})\cdot w_R\,=\,-w_R\,,\qquad
  \alpha_j\in\Pi(\fg)\,.
  \ee
Let us introduce the following matrix element in
$(\pi_{-\imath\g-\rho},\CV_{-\imath\g-\rho})$ defined as a function on
$A=\exp(\fh)$:
 \be\label{matel}
  \Phi_{\g}(e^x)\,
  =\,\bigl\<w_L,\,
  \pi_{-\imath\g-\rho}\bigl(e^{-\sum\limits_{i=1}^{\rk(\fg)}\varpi_i(x)h_i}\bigr)
  w_R\bigr\>\,,\qquad \g\in\IC^{\rk(\fg)},\quad e^x\in A\,,
 \ee
where $\varpi_i\in\fh^*$ are the fundamental weights of $\fg$. The
matrix element \eqref{matel} gives rise to the $\fg$-Whittaker
function,
 \be\label{whitfun}
  \Psi_{\g}(e^x)\,
  =\,e^{-\rho(x)}\,\Phi_{\g}(e^x)\,
  =\,e^{-\rho(x)}\,\bigl\<w_L,\,
  \pi_{-\imath\g-\rho}\bigl(e^{-\sum\limits_{i=1}^{\rk(\fg)}\varpi_i(x)h_i}\bigr)
  w_R\bigr\>\,,
 \ee
which is a solution to the quantum $\fg$-Toda chain \cite{Ko}. The
action of the quadratic quantum $\fg$-Toda Hamiltonian is defined
via
 \be
  \CH_2\cdot\Psi_{\g}(e^x)\,
  =\,e^{-\rho(x)}\,\bigl\<w_L,\,
  \pi_{-\imath\g-\rho}\bigl(e^{-\sum\limits_{i=1}^{\rk(\fg)}\varpi_i(x)h_i}\,c_2\,\bigr)
  w_R\bigr\>\,,
 \ee
where $c_2\in\CZ\CU(\gl_{\ell+1})$ is a degree two element of the
center of the enveloping algebra. For a special choice of $c_2$, we
have
 \be
  \CH_2\,
  =\,-\frac{1}{2}\sum_{i=1}^{\rk(\fg)}\frac{\pr^2}{\pr x_i^2}\,
  +\,\sum_{\alpha_j\in\Pi(\fg)}e^{\alpha_j(x)}\,.
 \ee

In the rest of this Section, we restrict ourselves to the case
$\fg=\gl_{\ell+1}(\IR)$. We consider two different  realizations of
the principal series $\CU(\gl_{\ell+1})$-representation
$(\pi_{-\imath\g-\rho},\,\CV_{-\imath\g-\rho})$ introduced in
\cite{GKL} and \cite{GKLO} that we refer to as the Gelfand-Tsetlin
and the Gauss-Givental realizations, respectively. In addition, we
provide another  version of the Gauss-Givental realization obtained
via the Fourier transform, that we refer to as the modified
Gauss-Givental realization. Note that in the Gelfand-Tsetlin and in
the modified Gauss-Givental realizations, the
$\gl_{\ell+1}(\IR)$-generators are represented by the first-order
difference operators; meanwhile in the Gauss-Givental realization,
the $\gl_{\ell+1}(\IR)$-generators are represented by the
first-order differential operators. In each case, we introduce an
explicit form of the pairing \eqref{pair} and the corresponding
integral formulas the matrix element \eqref{matel}. Thus, in the
Gelfand-Tsetlin realization, the Mellin-Barnes integral
representation arises. On the other hand, the Gauss-Givental
realization leads to the Givental integral representation
\cite{Giv},\cite{JK}. The modified Gauss-Givental realization leads
to yet another integral representation of the
$\mathfrak{gl}_{\ell+1}(\IR)$-Whittaker function.

\subsection{The Gelfand-Tsetlin realization}

According to \cite{KL}, the $\gl_{\ell+1}(\IR)$-Whittaker function
allows the following integral representation of the Mellin-Barnes
type:
 \be\label{MBint0}
  \Psi_{\g_1,\ldots,\g_{\ell+1}}(e^{x_1},\ldots,e^{x_{\ell+1}})\,
  =\,\frac{1}{(2\pi)^{\frac{\ell(\ell+1)}{2}}}
  \!\int\limits_{\CC}
  \prod_{k=1}^{\ell}\prod_{i=1}^kd\g_{k,i}\,
  e^{\imath\sum\limits_{k=1}^{\ell+1}\bigl
    (\sum\limits_i\g_{k,i}-\sum\limits_j\g_{k-1,j}\bigr)x_k}\\
  \times\prod_{n=1}^{\ell}\frac{\prod\limits_{i=1}^{n+1}\prod\limits_{j=1}^n
  \Gamma\bigl(\imath\g_{n+1,i}-\imath\g_{n,j}\bigr)}
  {\prod\limits_{a<b}\bigl|\Gamma(\imath\g_{n,a}-\imath\g_{n,b})\bigr|^2}\,,
 \ee
where $\g_{\ell+1,i}=\g_i\,,1\leq i\leq\ell+1$ and the terms
$\g_{k,j},\,j>k$ are omitted when occur. The integration domain
$\CC\subset\IC^{\frac{\ell(\ell+1)}{2}}$ in \eqref{MBint0} satisfies
the following:
 \be\label{GTcontour0}
  {\rm Im}(\g_{k,j})\,>\,\max_{1\leq i\leq k+1}\bigl\{
  {\rm Im}(\g_{k+1,\,i})\bigr\}\,,\qquad1\leq j\leq k\leq\ell\,.
 \ee
Upon choosing  a special solution to \eqref{GTcontour0}, for
$\tau_{k,j}\in\IR,\,1\leq j\leq k\leq\ell$,
 \be\label{GTcontour}
  \g_{k,j}\,=\,\tau_{k,j}\,+\,\imath\frac{\ell+1-k}{2}\,,
 \ee
we obtain  the following form of integral
 representation  (see (5.4)  in \cite{GKL})
 \be\label{MBint}
  \Psi_{\g_1,\ldots,\g_{\ell+1}}(e^{x_1},\ldots,e^{x_{\ell+1}})\,
  =\,\frac{e^{-\rho(x)}}{(2\pi)^{\frac{\ell(\ell+1)}{2}}}
  \!\int\limits_{\IR^{\frac{\ell(\ell+1)}{2}}}
  \prod_{k=1}^{\ell}\prod_{i=1}^kd\tau_{k,i}\\
  \times\,
  e^{\imath\sum\limits_{k=1}^{\ell+1}\bigl(\sum\limits_i\tau_{k,i}-\sum\limits_j\tau_{k-1,j}
    \bigr)x_k}\prod_{n=1}^{\ell}\frac{\prod\limits_{i=1}^{n+1}\prod\limits_{j=1}^n
  \Gamma\bigl(\imath(\tau_{n+1,i}-\tau_{n,j})+\frac{1}{2}\bigr)}
  {\prod\limits_{a<b}\bigl|\Gamma(\imath\tau_{n,a}-\imath\tau_{n,b})\bigr|^2}\,,
 \ee
where the terms $\tau_{k,j},\,j>k$ are omitted when occur and we
make the identification $\tau_{\ell+1,i}:=\g_i\in\IC,\,1\leq
i\leq\ell+1$. Recall that for $\fg=\fg_{\ell+1}(\IR)$, the half-sum
of positive roots $\rho=(\rho_1,\ldots,\rho_{\ell+1})$ is given by
$\rho_i=\frac{\ell}{2}+1-i$.

By \cite{GKL} the integral expression \eqref{MBint} is  obtained as
a matrix element in a particular realization  a principal series
$\gl_{\ell+1}(\IR)$-representation
$(\pi_{-\imath\g-\rho},\CV_{-\imath\g-\rho})$ labeled by
$\g=(\g_1,\ldots,\g_{\ell+1})\in\IC^{\ell+1}$.  This realization
 (called the Gelfand-Tsetlin realization in \cite{GKL}) is defined as follows. The
 images $E_{ij}=\pi_{-\imath\g-\rho}(e_{ij})$
of the elementary matrix generators $e_{ij}$, $1\leq i,j\leq\ell+1$
of $\mathfrak{gl}_{\ell+1}(\IR)$ acting in the space of functions in
variables $\tau_{k,i}\in\IR,\,1\leq i\leq k\leq\ell$ are given by
 \be\label{GTrep}
  E_{i,i}\,
  =\,-\imath\Big(\sum_{k=1}^i\tau_{i,k}\,-\,\sum_{j=1}^{i-1}\tau_{i-1,j}\Big)\,,\qquad1\leq
  i\leq\ell+1\,,\\
  E_{i,\,i+1}\,
  =\,\imath\sum_{k=1}^i\frac{\prod\limits_{j=1}^{i+1}\bigl(\tau_{i,k}-\tau_{i+1,j}-\frac{\imath}{2}\bigr)}
  {\prod\limits_{l\neq j}(\tau_{i,k}-\tau_{i,l})}\,
  e^{-\imath\pr_{\tau_{i,k}}}\,,\qquad1\leq i\leq\ell\,,\\
  E_{i+1,\,i}\,
  =\,-\imath\sum_{k=1}^i\frac{\prod\limits_{j=1}^{i-1}\bigl(\tau_{i,k}-\tau_{i-1,j}+\frac{\imath}{2}\bigr)}
  {\prod\limits_{l\neq j}(\tau_{i,k}-\tau_{i,l})}\,
  e^{\imath\pr_{\tau_{i,k}}}\,,\qquad1\leq i\leq\ell\,,
 \ee
where the terms $\tau_{k,j},\,j>k$ are omitted when occur and
 $\tau_{\ell+1,k}:=\g_k\in \IC, 1\leq k\leq \ell+1$.  The proper
   subspace in the space of functions in $\tau_{k,i}\in\IR,\,1\leq i\leq
k\leq\ell$ should be defined as the space of functions obtained by
actions of polynomials in generators \eqref{GTrep} on the Whittaker
vectors defined below. The corresponding functional models for
representations $\CV_{-\imath\g-\rho},\,\CV^{\vee}_{-\imath\g-\rho}$
of $\CU(\mathfrak{gl}_{\ell+1})$ allows the duality  pairing
\eqref{pair} (see \cite{GKL}, formulas (2.22),(2.23)):
 \be\label{GTpair}
  \<\psi_1,\,\psi_2\>_{GT}\,:
  =\!\int\limits_{\IR^{\frac{\ell(\ell+1)}{2}}}\!\prod_{k,i}d\tau_{k,i}\,\,
  \mu(\tau)\,\,
  \overline{\psi_1(\tau)}\,\psi_2(\tau)\,,
 \ee
where
 \be\label{GTmeasure}
  \mu(\tau)\,
  =\,(2\pi)^{-\sum\limits_{n=1}^{\ell}\frac{n(n+1)}{2}}
  \prod_{n=2}^{\ell}\prod_{a<b}(\tau_{n,a}-\tau_{n,b})
  \bigl(e^{2\pi\tau_{n,a}}-e^{2\pi\tau_{n,b}}\bigr)\\
  =\,\frac{1}{(2\pi)^{\frac{\ell(\ell+1)}{2}}}
  \prod_{n=2}^{\ell}\prod_{a<b}\frac{e^{\pi(\tau_{n,a}+\tau_{n,b})}}
  {\bigl|\Gamma(\imath\tau_{n,a}-\imath\tau_{n,b})\bigr|^2}\,.
 \ee
The difference operators \eqref{GTrep} satisfy, for  vectors
$\phi_1\in \CV^{\vee}_{-\imath\g-\rho}$ and $\phi_2\in
\CV_{-\imath\g-\rho}$
 \be
  \<E_{ij}^{\vee}\phi_1,\,\phi_2\>_{GT}\,
  =\,-\<\phi_1,\,E_{ij}\phi_2\>_{GT}\,,\qquad1\leq i,j\leq\ell+1\,.
 \ee

The defining relations \eqref{whit} on the Whittaker vectors
$\psi_R\in\CV_{-\imath\g-\rho}$ and
$\psi_L\in\CV^{\vee}_{-\imath\g-\rho}$ allow the following
solutions:
 \be\label{GTwhitL}
  \psi_L(\tau)\,
  =\,\prod_{n=1}^{\ell}\prod_{j=1}^ne^{-\frac{\pi}{2}\tau_{n,j}}\,,
  \ee
  \be\label{GTwhitR}
  \psi_R(\tau)\,
  =\,e^{-\pi\sum\limits_{k=1}^{\ell}(n-1)\sum\limits_{i=1}^k\tau_{k,i}}
  \prod_{n=1}^{\ell}\prod_{j=1}^ne^{\frac{\pi}{2}\tau_{n,j}}
  \prod_{i=1}^{n+1}
  \Gamma\Big(\imath(\tau_{n+1,i}-\tau_{n,j})+\frac{1}{2}\Big)\,.
 \ee

  Combining all ingredients in \eqref{whitfun}, we arrive at integral
representation  \eqref{MBint} of the
$\mathfrak{gl}_{\ell+1}(\IR)$-Whittaker function:
  \be\label{GTmatel}
   e^{-\rho(x)}\<\psi_L,\,e^{-\sum\limits_{i=1}^{\ell+1}x_iE_{ii}}\psi_R\>_{GT}\,
   =\,\frac{e^{-\rho(x)}}{(2\pi)^{\frac{\ell(\ell+1)}{2}}}
   \!\int\limits_{\IR^{\frac{\ell(\ell+1)}{2}}}
  \prod_{k=1}^{\ell}\prod_{i=1}^kd\tau_{k,i}\\
  \times\,
  e^{\imath\sum\limits_{k=1}^{\ell+1}\bigl(\sum\limits_i\tau_{k,i}-\sum\limits_j\tau_{k-1,j}\bigr)x_k}
  \prod_{n=1}^{\ell}\frac{\prod\limits_{i=1}^{n+1}\prod\limits_{j=1}^n
  \Gamma\bigl(\imath(\tau_{n+1,i}-\tau_{n,j})+\frac{1}{2}\bigr)}
  {\prod\limits_{a<b}\bigl|\Gamma(\imath\tau_{n,a}-\imath\tau_{n,b})\bigr|^2}\,,
 \ee
with $\tau_{\ell+1,i}:=\g_i\,, 1\leq i\leq \ell+1$.

Let us notice that since the $\mathfrak{gl}_{\ell+1}$-generators act
via difference operators \eqref{GTrep} in the space of functions in
variables $\{\tau_{k,i}\}$, there is a large commutant of this
action containing periodic functions in $\{\tau_{k,i}\}$. This
reflects the possibility to pick different solutions to the
equations \eqref{whit} on Whittaker vectors via multiplying the
expressions in \eqref{GTwhitL},\eqref{GTwhitR} by periodic
functions:
 \be\label{GTperiodic}
  \psi_L\longmapsto\s_L\psi_L\,,\qquad
  \psi_R\longmapsto\s_R\psi_R\,,\\
  \s_L(\tau\pm\imath)=\s_L(\tau)\,,\quad
  \s_R(\tau\pm\imath)=\s_R(\tau)\,.
  \ee
Specialization of the proper subspaces
$\CV_{-\imath\g-\rho},\,\CV^{\vee}_{-\imath\g-\rho}$ in  the
functional space corresponds to a choice of a model of the principle
series representation, so that in \eqref{GTwhitL}, \eqref{GTwhitR} a
particular choice is made. However, in the following it will be
convenient not to fix this choice beforehand and impose only the
condition on $\s_L,\,\s_R$ to satisfy the relation
$\ov{\s_L}\times\s_R=1$. This guarantees that the corresponding
matrix elements coincide with the integral expression for the matrix
element \eqref{GTmatel}.

\subsection{The Gauss-Givental realization}

In this subsection we recall from \cite{GKLO} a construction of
another integral  representation of the
$\gl_{\ell+1}(\IR)$-Whittaker function via representation theory.
This  integral representation first introduced  by Givental
\cite{Giv},\cite{JK} is given by
 \be\label{Givint}
  \Psi_{\g_1,\ldots,\g_{\ell+1}}(e^{x_1},\ldots,e^{x_{\ell+1}})\,
  =\,\int\limits_{\Gamma}\prod_{k=1}^{\ell}\prod_{i=1}^kdT_{k,i}\,
  e^{\mathcal{F}_{\g}(T)},
 \ee
where $T_{\ell+1,i}:=x_i$, and the (superpotential) function
$\mathcal{F}_{\g}(T)$ is given by
 \be\label{pot}
  \mathcal{F}_{\g}(T)\,
  =\sum_{k=1}^{\ell+1}\imath\g_k\Big(
  \sum_{i=1}^kT_{k,i}-\sum_{i=1}^{k-1}T_{k-1,i}\Big)\\
  -\,\sum_{k=1}^{\ell}\sum_{i=1}^k\Big(e^{T_{k+1,i}-T_{k,i}}
  +e^{T_{k,i}-T_{k+1,i+1}}\Big)\,.
 \ee
The integration domain  $\Gamma$ in \eqref{Givint}, is a middle
dimensional submanifold in $\IC^{\frac{\ell(\ell+1)}{2}}$
 such that the integral converges. In the following we
make a particular choice $\Gamma=\IR^{\frac{\ell(\ell+1)}{2}}$.

By \cite{GKLO}, given a principal series
$\gl_{\ell+1}(\IR)$-representation
$(\pi_{-\imath\g-\rho},\CV_{-\imath\g-\rho})$ labeled by
$\g=(\g_1,\ldots,\g_{\ell+1})\in\IC^{\ell+1}$, the Gauss-Givental
realization of $(\pi_{-\imath\g-\rho},\CV_{-\imath\g-\rho})$ is
defined by the  first-order differential operators acting in the
space of Schwartz functions in variables $T_{k,i},\,1\leq i\leq
k\leq\ell$. Precisely, the following expressions provide a
realization of the  principal series
$\gl_{\ell+1}(\IR)$-representation
 \be\label{GGrep}
  \CE_{i,i}\,
  =\,-\imath\g_i\,-\,\sum_{k=i}^{\ell}\frac{\pr}{\pr T_{k,i}}\,
  +\,\sum_{k=1}^{i-1}\frac{\pr}{\pr T_{\ell+1+k-i,k}}\,,\qquad1\leq i\leq\ell+1\,,\\
  \CE_{i,\,i+1}\,
  =\,\sum_{n=1}^i\left(\,\sum_{k=n}^ie^{T_{\ell+1+k-i,k}-T_{\ell+k-i,k}}\right)
  \!\!\Big(\frac{\pr}{\pr T_{\ell+n-i,n}}-\frac{\delta_{n,1}}{2}
  -\frac{\pr}{\pr T_{\ell+n-i,n-1}}\Big)\,,\\
  \CE_{i+1,\,i}\,
  =\,\sum_{k=i}^{\ell}e^{T_{k,i}-T_{k+1,i+1}}\Big\{\imath(\g_{i+1}-\g_i)-1+\frac{\delta_{i,\ell}}{2}\\
  +\,\sum_{s=i}^k\Big(\frac{\pr}{\pr T_{s,i+1}}+\frac{\delta_{s,1}}{2}
  -\frac{\pr}{\pr T_{s,i}}\Big)\Big\}\,,
 \ee
where the operators $\CE_{i,i+1}$ and $\CE_{i+1,i}$ are defined for
$1\leq i\leq\ell$. In \eqref{GGrep}, we impose $T_{\ell+1,k}=0$,
also for $k<j$ the terms $T_{k,j}$ and $\pr/\pr T_{k,j}$ are omitted
(when occur). Let us stress that the expressions \eqref{GGrep}
differ from the ones in \cite{GKLO} via conjugation by
$\prod\limits_{k,i}e^{T_{k,i}/2}$. The pairing \eqref{pair} in the
Gauss-Givental representation is given by the standard
$L^2$-pairing:
 \be\label{GGpair}
  \<\phi_1,\,\phi_2\>\,
  =\int\limits_{\IR^{\frac{\ell(\ell+1)}{2}}}\!\prod_{k,i}dT_{k,i}\,
  \overline{\phi_1(T)}\,
  \phi_2(T)\,.
 \ee
The defining relations \eqref{whit} on the Whittaker vectors
$\phi_R\in\CV_{-\imath\g-\rho}$ and
$\phi_L\in\CV^{\vee}_{-\imath\g-\rho}$ allow the following
solutions:
 \be\label{GGwhitL}
  \phi_L(T)\,
  =\,\prod_{k=1}^{\ell}
  e^{[\imath(\bar{\g}_{k+1}-\bar{\g}_k)-\frac{1}{2}]\sum\limits_{i=1}^kT_{k,i}\,
  -\sum\limits_{i=1}^ke^{T_{k+1,i}-T_{k,i}}}\,,
 \ee
 \be\label{GGwhitR}
  \phi_R(T)\,
  =\,\prod_{k=1}^{\ell}\prod_{i=1}^k
  e^{\frac{T_{k,i}}{2}\,-\,e^{T_{k,i}-T_{k+1,i+1}}}\,.
 \ee
Now the integral \eqref{Givint} coincides with the matrix element
representation in the Gauss-Givental realization:
 \be\label{GGmatel}
  \Psi_{\g_1,\ldots,\g_{\ell+1}}(e^{x_1},\ldots,e^{x_{\ell+1}})\,
  =\,e^{-\rho(x)}\<\phi_L,\,
  e^{-\sum\limits_{i=1}^{\ell+1}x_i\CE_{ii}}\phi_R\>\,.
 \ee

In the following we consider
 yet one more    realization of the principal series
 $\CU(\gl_{\ell+1})$-representation.
This realization is  straightforwardly obtained from the
Gauss-Givental realization via the Fourier transform and will be
called the modified Gauss-Givental realization. The details can be
found in Appendix A.

\begin{prop}\label{GGFourier} The following realization of
the principal series $\CU(\gl_{\ell+1})$-representation
$(\pi_{-\imath\g-\rho},\CV_{-\imath\g-\rho})$ in the space of
Schwartz functions $\CS(\IR^{\frac{\ell(\ell+1)}{2}})$ holds:
 \be\label{GGrepF}
  \hat{\CE}_{i,i}\,
  =\,-\imath\g_i\,-\,\imath\sum_{k=i}^{\ell}s_{k,i}\,
  +\,\imath\sum_{k=1}^{i-1}s_{\ell+1+k-i,k}\,,\qquad1\leq i\leq\ell+1\,,\\
  \hat{\CE}_{i,\,i+1}\,
  =\,\sum_{n=1}^i\Big(\imath s_{\ell-i+n,n}\\
  +\sum_{k=1}^{n-1}\imath(s_{\ell-i+n-k,k}-s_{\ell+1-i+n-k,k})+\frac{1}{2}\Big)\,
  e^{\imath\pr_{s_{\ell+1-i+n,n}}-\imath\pr_{s_{\ell-i+n,n}}}\,,\\
  \hat{\CE}_{i+1,\,i}\,
  =\sum_{n=i}^{\ell}\Big(\imath(\g_{i+1}-\g_i)-\imath
  s_{ii}\\
  +\sum_{k=i+1}^n\imath(s_{k,i+1}-s_{k,i})+\frac{1}{2}\Big)\,
  e^{\imath\pr_{s_{n,i}}-\imath\pr_{s_{n+1,i+1}}}\,,
 \ee
where the operators $\hat{\CE}_{i,i+1}$ and $\hat{\CE}_{i+1,i}$ for
$1\leq i\leq\ell$. In \eqref{GGrepF}, we impose $s_{k,j}=0,\,k<j$
and the terms $\pr_{s_{\ell+1,k}}$ are omitted when occur.
\end{prop}
\proof See Appendix A. $\Box$

The Whittaker vectors in this realization  are given by the
following solutions to \eqref{whit}:
 \be\label{GGwhitLF}
  \hat{\phi}_L(s)\,
  =\,\frac{1}{(2\pi)^{\frac{\ell(\ell+1)}{4}}}\prod_{k=1}^{\ell}\prod_{i=1}^{\ell+1-k}
  \Gamma\Big(\imath(\bar{\g}_i-\bar{\g}_{i+k})\,
  +\,\imath\sum_{n=0}^{k-1}s_{i+k-n-1,\,i}\,+\,\frac{k}{2}\Big)\,,
 \ee
and
 \be\label{GGwhitRF}
  \hat{\phi}_R(s)\,
  =\,\frac{1}{(2\pi)^{\frac{\ell(\ell+1)}{4}}}\prod_{k=1}^{\ell}\prod_{i=k}^{\ell}
  \Gamma\Big(\frac{k}{2}\,-\,\imath\sum_{n=0}^{k-1}s_{i-n,\,k-n}\Big)\,.
 \ee
Here we assume that the imaginary parts of the representation
parameters $\g\in\IC^{\ell+1}$ are small enough. In this
realization, the  matrix element in \eqref{GGmatel} affords the
following integral expression (via the Parseval identity):
 \be\label{GGintF}
  \Psi_{\g_1,\ldots,\g_{\ell+1}}(e^{x_1},\ldots,e^{x_{\ell+1}})\,
  =\,e^{-\rho(x)}\<\hat{\phi}_L,\,
  e^{-\sum\limits_{i=1}^{\ell+1}x_i\hat{\CE}_{ii}}\hat{\phi}_R\>\\
  =\,\frac{e^{\sum\limits_{n=1}^{\ell+1}(\imath\g_n-\rho_n)x_n}}{(2\pi)^{\frac{\ell(\ell+1)}{2}}}\!
  \int\limits_{\IR^{\frac{\ell(\ell+1)}{2}}}\prod_{k=1}^{\ell}\prod_{i=1}^kds_{k,i}\,
  e^{\imath\sum\limits_{n=1}^{\ell+1}\bigl(\sum\limits_{k=n}^{\ell}s_{k,n}\,
  -\sum\limits_{k=1}^{n-1}s_{\ell+1+k-n,k}\bigr)\,x_n}\\
  \times\prod_{k=1}^{\ell}\prod_{i=k}^{\ell}
  \Gamma\Big(-\imath\sum_{n=0}^{k-1}s_{i-n,\,k-n}\,+\,\frac{k}{2}\Big)\\
  \times\prod_{k=1}^{\ell}
  \prod_{i=1}^{\ell+1-k}
  \Gamma\Big(-\imath(\g_i-\g_{i+k})\,
  -\,\imath\sum_{n=0}^{k-1}s_{i+k-n-1,\,i}\,+\,\frac{k}{2}\Big)\,.
  \ee
Let us stress that the  transformation of the integral
representation \eqref{Givint} into \eqref{GGmatel},\eqref{GGintF} is
straightforward and is given by the  Fourier transform (see Appendix
A for details). On the other hand, the relation of the Mellin-Barnes
integral representation \eqref{MBint}  and the Givental
representation \eqref{Givint} as well as its Fourier transformed
version \eqref{GGintF} is a more delicate matter.
 In the rest of the note, we provide a direct identification of
the Mellin-Barnes integral representation and the modified Givental
representation  for  small ranks $\ell=1,2$ thus proving an
equivalence of all the three integral expressions of the
$\mathfrak{gl}_{\ell+1}(\IR)$-Whittaker function considered above,
for $\ell=1,2$ .


\section{Identification of matrix element realizations for $\gl_2(\IR)$}


As a preparation for a more complicated case of the
$\mathfrak{gl}_3(\IR)$-Whittaker function, in this  section we
consider an identification of two realizations of the
$\mathfrak{gl}_2(\IR)$-Whittaker function. Precisely, we transform
the corresponding matrix element written explicitly using  the
modified Gauss-Givental realization into the same matrix element
written in the Gelfand-Tsetlin realization. Although relevant
calculations are trivial in the $\mathfrak{gl}_2(\IR)$ case,  we use
this exercise as a warm-up for a more involved case of the
$\mathfrak{gl}_3(\IR)$-Whittaker function. In particular, in view of
 the $\mathfrak{gl}_3(\IR)$-case studied further, we consider a slightly more general
constructions of representations than it is needed for the
$\mathfrak{gl}_2(\IR)$-case.

For $\g=(\g_1,\g_2)\in\IC^2$, consider the principal series
$\CU(\gl_2)$-representation
$(\pi_{-\imath\g-\rho},\,\CV_{-\imath\g-\rho})$. Then the
Mellin-Barnes integral representation \eqref{MBint0} of the
$\gl_2(\IR)$-Whittaker function reads
 \be\label{MBgl2int0}
  \Psi_{\g_1,\g_2}(e^{x_1},e^{x_2})\,
  =\,\frac{1}{2\pi}\!\int\limits_{\CC}\!d\xi\,
  e^{\imath\xi x_1+\imath(\g_1+\g_2-\xi)x_2}\,
  \Gamma(\imath\g_1-\imath\xi)\,
  \Gamma(\imath\g_2-\imath\xi)\,,
 \ee
where the integration domain  $\CC\subset\IC$ satisfies the
conditions \eqref{GTcontour0}:
 \be\label{gl2contour0}
  {\rm Im}(\xi)\,>\,\max\bigl\{{\rm Im}(\g_1),\,{\rm Im}(\g_2)\bigr\}\,.
 \ee
Let us specify the representation parameters
$\g=(\g_1,\g_2)\in\IC^2$ in terms of $\la=(\la_1,\la_2)\in\IR^2$ as
follows:
 \be\label{gl2lambda}
  \g\,:=\,\la+2\imath\ve\rho\,,\quad
  \rho=\Big(\frac{1}{2}\,,-\frac{1}{2}\Big)\,,
  \ee
for some fixed $\ve$ such that $0<\ve<\frac{1}{2}$. Then
\eqref{MBgl2int0} takes the following form:
 \be\label{MBgl2int1}
  \Psi_{\la_1+\imath\ve,\,\la_2-\imath\ve}(e^{x_1},e^{x_2})\,
  =\,\frac{1}{2\pi}\!\int\limits_{\CC}\!d\xi\,
  e^{\imath\xi x_1+\imath(\la_1+\la_2-\xi)x_2}\\
  \times\,\Gamma\bigl(\imath(\la_1-\xi)-\ve\bigr)\,
  \Gamma\bigl(\imath(\la_2-\xi)+\ve\bigr)\,,
 \ee
where for the integration domain $\CC$ we have ${\rm Im}(\xi)>\ve$.
Now by \eqref{GTcontour}, we specify the integration domain to be
parameterized by variable $\tau\in \IR$ as follows:
 \be
  \xi\,=\,\left(\tau\,+\,\frac{\imath}{2}\right)\,\in  \CC=\left(\IR+\frac{\imath}{2}\right)\,,
 \ee
so that after substitution, \eqref{MBgl2int1} takes the form:
 \be\label{MBgl2int11}
  \Psi_{\la_1+\imath\ve,\,\la_2-\imath\ve}(e^{x_1},e^{x_2})\,
  =\,\frac{e^{\frac{x_2-x_1}{2}}}{2\pi}\!\int\limits_{\IR}\!d\tau\,
  e^{\imath\tau x_1+\imath(\la_1+\la_2-\tau)x_2}\,\\
  \times\,\Gamma\Big(\imath(\la_1-\tau)+\ve+\frac{1}{2}\Big)\,
  \Gamma\Big(\imath(\la_2-\tau)-\ve+\frac{1}{2}\Big)\,.
  \ee

\subsection{The Gelfand-Tsetlin realization}

By \eqref{GTrep} for $\ell=1$, the Gelfand-Tsetlin realization of
the principle series $\mathfrak{gl}_2(\IR)$-representa-tion
$(\pi_{-\imath\g-\rho},\CV_{-\imath\g-\rho})=(\pi_{-\imath\la-(1-2\ve)\rho},\CV_{-\imath\la-(1-2\ve)\rho})$
 reads
 \be\label{GTgl2rep}
  E_{11}=-\imath\tau,\qquad E_{22}=-\imath(\la_1+\la_2)+\imath\tau,\\
  E_{12}=-\imath\Big(\imath(\tau-\la_1)+\ve+\frac{1}{2}\Big)
  \Big(\imath(\tau-\la_2)-\ve+\frac{1}{2}\Big)\,e^{-\imath\pr_{\tau}},\quad
  E_{21}=-\imath e^{\imath\pr_{\tau}}\,.
 \ee
The Whittaker vectors satisfy the defining equations \eqref{whit}
and are given by \eqref{GTwhitL},\eqref{GTwhitR} for $\ell=1$:
 \be\label{GTgl2whit}
  \psi_L(\tau)=e^{-\frac{\pi\tau}{2}}\,,\quad
  \psi_R(\tau)=e^{\frac{\pi\tau}{2}}
  \Gamma\Big(\imath(\la_1-\tau)-\ve+\frac{1}{2}\Big)\,
  \Gamma\Big(\imath(\la_2-\tau)+\ve+\frac{1}{2}\Big)\,.
 \ee
For $\ell=1$, \eqref{GTmeasure} reads $\mu(\tau)=\frac{1}{2\pi}$.
Then \eqref{MBgl2int11} allows for the matrix element
representation:
 \be\label{GTgl2matel0}
  \Psi_{\la_1+\imath\ve,\,\la_2-\imath\ve}(e^{x_1},e^{x_2})\,
  =\,e^{-\rho(x)}\<\psi_L,\,e^{-x_1E_{11}-x_2E_{22}}\,\psi_R\>_{GT}\\
  =\,\frac{e^{-\rho(x)}}{2\pi}\!\int\limits_{\IR}\!d\tau\,
  e^{\imath\tau x_1+\imath(\la_1+\la_2-\tau)x_2}\,
  \Gamma\Big(\imath(\la_1-\tau)-\ve+\frac{1}{2}\Big)\,
  \Gamma\Big(\imath(\la_2-\tau)+\ve+\frac{1}{2}\Big)\,.
 \ee

To account for a more general integral representation we consider a
slightly modified construction of the matrix element. For
$\vk\in\IR$, we introduce the following deformation of the
$L^2$-pairing, for $\psi_1\in\CV_{-\imath\g-\rho}^{\vee}$ and
$\psi_2\in\CV_{-\imath\g-\rho}$:
 \be\label{DefPAIR}
  \<\,\,,\,\>_{\vk}\,:\,\,
  \CV^{\vee}_{-\imath\g-\rho}\times
  \CV_{-\imath\g-\rho}\,\longrightarrow\,\IC\,,\qquad
  \<\psi_1,\psi_2\>_{\vk}\,:
  =\!\int\limits_{\IR-\imath\vk}\!\!\frac{d\tau}{2\pi}\,
  \bar{\psi_1}(\tau)\,\psi_2(\tau)\,,
 \ee
so that $\<\,\,,\,\>_0=\<\,\,,\,\>_{GT}$. Then \eqref{DefPAIR} might
be rewritten as follows
 \be
  \<\psi_1,\psi_2\>_{\vk}\,
  =\!\int\limits_{\IR}\!\frac{d\tau}{2\pi}\,
  \bar{\psi_1}(\tau-\imath\vk)\,\psi_2(\tau-\imath\vk)\,
  =\!\int\limits_{\IR}\!\frac{d\tau}{2\pi}\,
  \overline{\psi_1(\tau+\imath\vk)}\,\psi_2(\tau-\imath\vk)\,.
 \ee
Given the  Gelfand-Tsetlin realization \eqref{GTgl2rep}, consider
the operators $E_{ij}$ with shifted argument $\tau_{11}\mapsto
(\tau_{11}-\imath\vk)$:
 \be\label{GTgl2repk}
  \KE_{ij}(\tau)\,:
  =\,E_{ij}(\tau-\imath\vk)\,,\qquad
  1\leq i,j\leq2\,.
 \ee
\begin{lem}\label{gl2GTcontour}
Let $\psi_L,\,\psi_R$ be the Whittaker vectors given by
\eqref{GTgl2whit}. Then for $\vk\in\IR$ subjected
$\vk<\frac{1}{2}-\ve$, the following holds:
 \be\label{IDINT}
  \<\psi_L,e^{-\sum x_i\KE_{ii}}\psi_R\>_{\vk}\,
  =\,\<\psi_L,e^{-\sum x_iE_{ii}}\psi_R\>_{GT}\,.
 \ee
\end{lem}
\proof We need to show that for $\vk<\frac{1}{2}-\ve$ the following
holds:
 \be
  \int\limits_{\IR-\imath\vk}\!\!d\tau\,
  \bar{\psi_1}(\tau)\,\psi_2(\tau)\,
  -\,\int\limits_{\IR}\!d\tau\,
  \bar{\psi_1}(\tau)\,\psi_2(\tau)\,
  =\!\int\limits_{\CC_{\vk}}\!d\tau\,
  \bar{\psi_1}(\tau)\,\psi_2(\tau)
  =\,0\,,
 \ee
 where $\CC_{\vk}=\IR\sqcup(\IR-\imath\vk)$ and
 \be
  \psi_1(\tau)\,=\,\psi_L(\tau)=e^{-\frac{\pi\tau}{2}}\,,\\
  \psi_2(\tau)\,=\,e^{-\sum x_iE_{ii}(\tau)}\psi_R(\tau)\,
  =e^{\frac{\pi\tau}{2}}e^{\imath\tau(x_1-x_2)+\imath(\g_1+\g_2)x_2}\,
  \prod_{i=1}^2\Gamma\Big(\imath(\g_i-\tau)+\frac{1}{2}\Big)\,,
 \ee
with $\g=\la+2\imath\ve\rho$, so that the integrand function reads
 \be
  f(\tau)\,
  =\,\ov{\psi}_1(\tau)\,\psi_2(\tau)\,
  =\,e^{\imath\tau x_1+\imath(\g_1+\g_2-\tau)x_2}
  \prod_{i=1}^2\Gamma\Big(\imath(\g_i-\tau)+\frac{1}{2}\Big)\,.
 \ee
Consider the contour
 \be
  \CC_{\vk,R}\,=\,[-R;R]\cup[R;R-\imath\vk]
  \cup[-R-\imath\vk;R-\imath\vk]\cup[-R;-R-\imath\vk]\,,
 \ee
then we claim that
 \be
  \lim_{R\to+\infty}\int\limits_{-R-\imath\vk}^{-R}\!
  d\tau\,f(\tau)\,
  =\lim_{R\to+\infty}\int\limits_{R-\imath\vk}^R\!
  d\tau\,f(\tau)\,
  =\,0\,.
 \ee
Indeed, we apply the following version of the Stirling asymptotic
formula (see \cite{BE}, 1.18, formula (6)),
 \be
  \bigl|\Gamma(s+\imath t)\bigr|\,
  \sim\,\sqrt{2\pi}\,e^{-\frac{\pi}{2}\,|t|}\,|t|^{s-\frac{1}{2}}\,,\quad|t|\to\infty\,,
 \ee
which implies
 \be
  |f(\tau)|\,
  =\,\prod_{i=1}^2\Big|\Gamma\Big(\imath(\g_i-\tau)+\frac{1}{2}\Big)\Big|\,
  \sim\,2\pi\,e^{-\frac{\pi}{2}|\g_1-\tau|-\frac{\pi}{2}|\g_2-\tau|}\,,\quad
  |\tau|\to\infty\,.
 \ee
Hence the integrals over vertical parts give zero contributions,
therefore
 \be
  \lim_{R\to+\infty}\int\limits_{-R}^R\!
  d\tau\,f(\tau)\,
  =\lim_{R\to+\infty}\int\limits_{-R-\imath\vk}^{R-\imath\vk}\!
  d\tau\,f(\tau)\,.
 \ee
Next, $f(\tau)$ has the poles, all below the real line
$\IR\subset\IC$, since $\g=\la+2\imath\ve\rho$ with
$0<\ve<\frac{1}{2}$:
 \be
  \tau\,=\,\g_i-\imath\Big(n+\frac{1}{2}\Big)\,
  =\,\la_i\,+\,\imath\Big(2\ve\rho_i-n-\frac{1}{2}\Big)\,,\quad
  n\geq0\,,\quad i=1,2\,.
 \ee
None of the poles are sitting inside the strip bounded by
$\CC_{\vk}$ if and only if $\vk<\frac{1}{2}-\ve$, hence we may
safely deform the integration contour to establish the equivalence
of the integral expressions in \eqref{IDINT}. $\Box$

Thus, for $\vk<\frac{1}{2}-\ve$ by Lemma \ref{gl2GTcontour}, we
arrive at the  following matrix element representation for the
$\gl_2(\IR)$-Whittaker function:
 \be\label{GTgl2matel}
  \Psi_{\la_1+\imath\ve,\la_2-\imath\ve}(e^{x_1},e^{x_2})\,
  =\,e^{-\rho(x)}\<\psi_L,\,e^{-x_1\KE_{11}-x_2\KE_{22}}\psi_R\>_{\vk}\\
  =\,\frac{e^{(\vk-\frac{1}{2})(x_1-x_2)}}{2\pi}\!\int\limits_{\IR}\!d\tau\,
  e^{\imath\tau x_1+\imath(\la_1+\la_2-\tau)x_2}\\
  \times\,\Gamma\Big(\imath(\la_1-\tau)-\ve-\vk+\frac{1}{2}\Big)\,
  \Gamma\Big(\imath(\la_2-\tau)+\ve-\vk+\frac{1}{2}\Big)\,.
 \ee
This  reproduces a more general integral expression  for the
 matrix element \eqref{MBgl2int11} in representation theory framework.

\subsection{The Gauss-Givental realization}

Now we compare the  integral representation \eqref{GTgl2matel} for
$\mathfrak{gl}_2(\IR)$-Whittaker function considered above with the
corresponding  modified Givental integral representation
\eqref{GGintF}. By \eqref{GGrepF} for $\ell=1$, the modified
Gauss-Givental realization of the principle series
$\mathfrak{gl}_2(\IR)$-representation
$(\pi_{-\imath\g-\rho},\CV_{-\imath\g-\rho})$ is given by
 \be\label{GGgl2repF}
  \hat{\CE}_{11}=-\imath(\g_1+s)\,,\qquad
  \hat{\CE}_{22}=-\imath(\g_2-s)\,,\\
  \hat{\CE}_{12}=\Big(\imath s+\frac{1}{2}\Big)\,e^{-\imath\pr_s}\,,\qquad
  \hat{\CE}_{21}=\Big(\imath(\g_2-\g_1-s)+\frac{1}{2}\Big)\,e^{\imath\pr_s}\,.
 \ee
By \eqref{dualrep}, the dual representation satisfies
$\CV^{\vee}_{-\imath\g-\rho}\simeq\CV_{-\imath\bar{\g}-\rho}$, so
that
 \be
  \hat{\CE}^{\vee}_{11}=-\imath(\bar{\g}_1+s)\,,\qquad
  \hat{\CE}^{\vee}_{22}=-\imath(\bar{\g}_2-s)\,,\\
  \hat{\CE}^{\vee}_{12}=\Big(\imath s+\frac{1}{2}\Big)\,e^{-\imath\pr_s}\,,\qquad
  \hat{\CE}^{\vee}_{21}=\Big(\imath(\bar{\g}_2-\bar{\g}_1-s)+\frac{1}{2}\Big)\,e^{\imath\pr_s}\,.
 \ee
Then for $\g=\la+2\imath\ve\rho$, the corresponding Whittaker
vectors are given by \eqref{GGwhitLF},\eqref{GGwhitRF} for $\ell=1$:
 \be\label{GGgl2Fwhit}
  \hat{\phi}_L(s)
  =\,\frac{1}{\sqrt{2\pi}}\,\Gamma\Big(\imath(\la_1-\la_2+s)+2\ve+\frac{1}{2}\Big)\,,
 \quad
  \hat{\phi}_R(s)\,=\,\frac{1}{\sqrt{2\pi}}\,\Gamma\Big(\frac{1}{2}-\imath s\Big)\,,
 \ee
and the resulting integral representation of the
$\mathfrak{gl}_2(\IR)$-Whittaker function is as follows
\be\label{GGgl2Fmatel}
  \Psi_{\la_1+\imath\ve,\,\la_2-\imath\ve}(e^{x_1},e^{x_2})\,
   =\,e^{-\rho(x)}\<\hat{\phi}_L,\,
  e^{-x_1\hat{\CE}_{11}-x_2\hat{\CE}_{22}}\hat{\phi}_R\>\\
  =\,\frac{e^{-(\ve+\frac{1}{2})(x_1-x_2)} }{2\pi}\!\!
  \int\limits_{\IR}\!\!ds\,e^{\imath(\la_1x_1+\la_2x_2)+\imath
  s(x_1-x_2)}\\
  \times\,\Gamma\Big(-\imath(\la_1-\la_2+s)+2\ve+\frac{1}{2}\Big)\,
  \Gamma\Big(-\imath s+\frac{1}{2}\Big)\,.
  \ee
The integrals \eqref{GGgl2Fmatel} and \eqref{GTgl2matel} can be
easily identified via changing the integration variable:
 \be\label{gl2ID}
  s\,:=\,\tau-\la_1\,\in\,\IR\,,\quad\text{for}\quad \vk=-\ve\,.
 \ee
Although the transformation \eqref{gl2ID} is quite trivial, we would
like to represent it by an ope-rator intertwining the two
realizations of the principal series representation
$(\pi_{-\imath\g-\rho},\CV_{-\imath\g-\rho})$.

\subsection{Intertwining operators}

The form of the intertwining operator may be easily inferred by
comparing Whittaker vectors in the two realizations of the
$\gl_2(\IR)$-representation
$\CV_{-\imath\g-\rho},\,\g=\la+2\imath\ve\rho,\,\la\in\IR^2$. The
standard delta-function is defined by
 \be
  \int\limits_{\IR}\!ds\,\,\delta(s-\tau)\,\phi(s)\,
  =\,\phi(\tau)\,,\qquad\phi\in\CS(\IR)\,.
 \ee
Introduce the delta-function of complex argument defined by (see
e.g. \cite{GS})
 \be\label{Delta}
  \int\limits_{\IR}\!ds\,\,\delta(s-\tau+\imath\ve)\,\phi(s)\,
  =\,e^{-\imath\ve\pr_{\tau}}\!\int\limits_{\IR}\!ds\,\,
  \delta(s-\tau)\,\phi(s)\,
  =\,\phi(\tau-\imath\ve)\,.
 \ee
The corresponding space of test functions consists of Schwartz
functions allowing analytic continuation into the strip bounded by
$\IR$ and $\IR-\imath\ve$. Now, introduce the operators $\CB_R$ and
$\CB_L$, acting in $\CV_{-\imath\g-\rho}$ and
$\CV_{-\imath\g-\rho}^{\vee}\simeq\CV_{-\imath\bar{\g}-\rho}$,
respectively, defined by their integral kernels:
 \be\label{gl2BR}
  B_R(\tau;s)\,
  =\,e^{\frac{\pi\tau}{2}}\,
  \Gamma\Big(\imath(\g_2-\g_1-s)+\frac{1}{2}\Big)\,
  \delta(\g_1+s-\tau)\\
  =\,e^{\frac{\pi\tau}{2}}\,
  \Gamma\Big(\imath(\la_2-\la_1-s)+2\ve+\frac{1}{2}\Big)\,
  \delta(\la_1+s-\tau+\imath\ve)\,,
 \ee
and
 \be\label{gl2BL}
  B_L(\tau;s)\, =\,\frac{e^{-\frac{\pi\tau}{2}}}
  {\Gamma\bigl(\imath(\bar{\g}_1-\bar{\g}_2+s)+\frac{1}{2}\bigr)}\,
  \delta(\bar{\g}_1+s-\tau)\\
  =\,\frac{e^{-\frac{\pi\tau}{2}}}
  {\Gamma\bigl(\imath(\la_1-\la_2+s)+2\ve+\frac{1}{2}\bigr)}\,
  \delta(\la_1+s-\tau-\imath\ve)\,.
 \ee
Let $\CB_L^{\dag}$ be the adjoint operator w.r.t. the standard
$L^2$-pairing, then its kernel reads
 \be\label{gl2BLdag}
  B_L^{\dag}(s;\tau)\,
  =\,\ov{B_L(\tau;s)}\,
  =\,\frac{e^{-\frac{\pi\tau}{2}}\,}
 {\Gamma\bigl(\imath(\g_2-\g_1-s)+\frac{1}{2}\bigr)}\,
  \delta(\g_1+s-\tau)\\
  =\,\frac{e^{-\pi\tau}}{\Gamma\bigl(\imath(\g_2-\tau)+\frac{1}{2}\bigr)^2}\,
  \times\,B_R(\tau;s)\,,
 \ee
and we have  $\CB_L^{\dag}=\CB_R^{-1}$\,.

\begin{lem}\label{gl2B}
The operators $\CB_R,\,\CB_L$ and $\CB^{\dag}_L$ possess the
following intertwining properties:
 \be\label{gl2Bintertw}
  E_{ij}\circ\CB_R\,=\,\CB_R\circ\hat{\CE}_{ij}\,,\qquad
  E^{\vee}_{ij}\circ\CB_L\,=\,\CB_L\circ\hat{\CE}^{\vee}_{ij}\,,\\
  \hat{\CE}_{ij}\circ\CB^{\dag}_L\,=\,\CB_L^{\dag}\circ E_{ij}\,,\qquad
  1\leq i,j\leq2\,.
 \ee
\end{lem}
\proof The intertwining relations for the Cartan generators $E_{ii}$
and $\hat{\CE}_{ii}$ follow from \eqref{gl2ID} due to the presence
of the delta-factors $\delta(\g_1+s-\tau)$ and
$\delta(\bar{\g}_1+s-\tau)$ in \eqref{gl2BR} and \eqref{gl2BL},
respectively. Relations for $E_{ij},\hat{\CE}_{ij},\,i\neq j$ and
$E^{\vee}_{ij},\hat{\CE}^{\vee}_{ij},\,i\neq j$ can be verified by
direct computations. $\Box$

\begin{lem}\label{gl2Bid}
The operators $\CB_R,\,\CB_L$ and $\CB^{\dag}_L$ satisfy the
following:
 \be\label{gl2Bwhit}
  \CB_R\cdot\hat{\phi}_R\,=\,\psi_R\,,\qquad
  \CB_L\cdot\hat{\phi}_L\,=\,\psi_L\,,\\
  (\CB^{\dag}_L\circ\CB_R)\cdot\hat{\phi}_R\,
  =\,\CB^{\dag}_L\cdot\psi_R\,=\,\hat{\phi}_R\,.
 \ee
\end{lem}
\proof All the three identities can be easily verified in
straightforward way. $\Box$

Now we identify the realizations \eqref{GTgl2matel0} and
\eqref{GGgl2Fmatel} of the $\gl_2(\IR)$-matrix element as follows:
 \be
  \<\psi_L,\,e^{-\sum x_iE_{ii}}\psi_R\>_{GT}\\
  =\,\<\CB_L\hat{\phi}_L,\,
  e^{-\sum x_iE_{ii}}\CB_R\hat{\phi}_R\>_{GT}\,
  =\,\<\hat{\phi}_L,\,
  \CB^{\dag}_L\,e^{-\sum x_iE_{ii}}\CB_R\hat{\phi}_R\>\\
  =\,\<\hat{\phi}_L,\,
  e^{-\sum x_i\hat{\CE}_{ii}}\CB^{\dag}_L\circ\CB_R\hat{\phi}_R\>\,
  =\,\<\hat{\phi}_L,\,
  e^{-\sum x_i\hat{\CE}_{ii}}\hat{\phi}_R\>\,.
 \ee
where we apply the third identity in Lemma \ref{gl2B} providing
 \be
  \CB^{\dag}_L\,e^{-\sum x_iE_{ii}}(\CB^{\dag}_L)^{-1}\,
  =\,e^{-\sum x_i\hat{\CE}_{ii}}\,.
 \ee

It is instructive to  consider the special case of $\g=\la\in\IR^2$
in \eqref{gl2lambda}, that is setting $\ve=0$. Upon this choice, the
$\gl_2(\IR)$-representation $\CV_{-\imath\la-\rho}$ has unitary
structure and both sets of $\gl_2(\IR)$-generators \eqref{GTgl2rep}
and \eqref{GGgl2repF} are Hermitian w.r.t. the appropriate pairings.
The transformation between the two realizations of
  matrix elements  constructed above
does not however take into account the unitary structure of the
underlying representation $\CV_{-\imath\la-\rho}$. To rectify this,
let us introduce the following operator:
 \be\label{gl2N}
  \CN\,:
  =\,e^{\imath\Arg\bigl\{\Gamma\bigl(\imath(\la_2-\tau)+\frac{1}{2}\bigr)\bigr\}}\,
  e^{-\la_1\pr_{\tau}}\,
  =\,\s_{\la_2}(\tau)\times\CB_R\,
  =\,\frac{1}{\s_{\la_2}(\tau)}\times\CB_L\,,
 \ee
where $\s_{\la_2}(\tau)$ is the $\imath$-periodic function in $\tau$
given by
 \be\label{gl2sigma}
  \s_{\la_2}(\tau)\,=\,\frac{e^{-\frac{\pi\tau}{2}}}
  {\Big|\Gamma\bigl(\imath(\la_2-\tau)+\frac{1}{2}\bigr)\Big|}\,
  =\,\frac{e^{-\frac{\pi\la_2}{2}}}{\sqrt{2\pi}}\,\sqrt{1+e^{2\pi(\la_2-\tau)}}\,.
 \ee

\begin{prop}
The operator $\CN$ is the unitary operator intertwining the two
realizations of the $\mathfrak{gl}_2(\IR)$-representation
$\CV_{-\imath\la-\rho}$ given by \eqref{GTgl2rep} and
\eqref{GGgl2repF}:
 \be\label{gl2Nintert}
  E_{ij}\circ\CN\,=\,\CN\circ\hat{\CE}_{ij}\,,\qquad1\leq i,j\leq2\,.
 \ee
The $\CN$-action on the Whittaker vectors \eqref{GGgl2Fwhit} reads
 \be\label{gl2Nwhit}
  \CN\cdot\hat{\phi}_R\,=\,\s_{\la_2}\,\psi_R\,,\qquad
  \CN\cdot\hat{\phi}_L\,=\,\frac{1}{\s_{\la_2}}\,\psi_L\,\,.
 \ee
\end{prop}
\proof By \eqref{gl2N}, the operator $\CN$ has the following
integral kernel:
 \be\label{gl2Nker}
  N(\tau;s)\,
  =\,\frac{\Gamma\bigl(\imath(\la_2-\tau)+\frac{1}{2}\bigr)}
  {\Big|\Gamma\bigl(\imath(\la_2-\tau)+\frac{1}{2}\bigr)\Big|}\,
  \delta(\la_1+s-\tau)\\
  =\,e^{\imath\Arg\bigl(\Gamma\bigl(\imath(\la_2-\tau)+\frac{1}{2}\bigr)\bigr)}\,
  \delta(\la_1+s-\tau)\,.
 \ee
The intertwining relations \eqref{gl2Nintert} follow from
\eqref{gl2Bintertw} in straightforward way. The relations
\eqref{gl2Nwhit} follow from \eqref{gl2Bwhit} and the identity
\eqref{gl2N}. $\Box$

The unitary  operator $\CN$ provides a unitary transformation
identifying the two realizations of $\gl_2(\IR)$-matrix elements as
follows:
 \be
  \<\psi_L,\,e^{-\sum x_iE_{ii}}\psi_R\>_{GT}\,
  =\,\<\s_{\la_2}^{-1}\psi_L,\,e^{-\sum x_iE_{ii}}\,\s_{\la_2}\,\psi_R\>_{GT}\\
  =\,\<\CN\hat{\phi}_L,\,
  e^{-\sum x_iE_{ii}}\CN\hat{\phi}_R\>\,
  =\,\<\hat{\phi}_L,\,
  \CN^{\dag}\,e^{-\sum x_iE_{ii}}\CN\hat{\phi}_R\>\\
  =\,\<\hat{\phi}_L,\,
  e^{-\sum x_i\hat{\CE}_{ii}}\CN^{\dag}\circ\CN\hat{\phi}_R\>\,
  =\,\<\hat{\phi}_L,\,
  e^{-\sum x_i\hat{\CE}_{ii}}\hat{\phi}_R\>\\
  =\,\<\phi_L,\,
  e^{-\sum x_i\CE_{ii}}\phi_R\>\,.
 \ee


\section{Identification of matrix element realizations for $\gl_3(\IR)$}


In this section we give a representation-theoretic identification of
the two integral formulas for the $\gl_3(\IR)$-Whittaker function
provided by the two realizations of principal series
$\CU(\gl_3)$-representation $\CV_{-\imath\g-\rho}$, the
Gelfand-Tsetlin and the modified Gauss-Givental ones.

For $\g=(\g_1,\g_2,\g_3)\in\IC^3$, consider the principal series
 $\CU(\gl_3)$-representation \\
$(\pi_{-\imath\g-\rho},\,\CV_{-\imath\g-\rho})$. The Mellin-Barnes
integral representation \eqref{MBint0} of the $\gl_3(\IR)$-Whittaker
function reads
 \be\label{MBgl3int0}
  \Psi_{\g_1,\g_2,\g_3}(e^{x_1},e^{x_2},e^{x_3})\,
  =\,\frac{e^{\imath(\g_1+\g_2+\g_3)x_3}}{(2\pi)^3}\!
  \int\limits_{\CC}\!
  \frac{d\g_{11}d\g_{21}d\g_{22}}
  {\Gamma(\imath\g_{21}-\imath\g_{22})\,
  \Gamma(\imath\g_{22}-\imath\g_{21})}\\
  \times\,
  e^{\imath(\g_{21}+\g_{22})(x_2-x_3)+\imath\g_{11}(x_1-x_2)}
  \prod_{j=1}^2\Gamma(\imath\g_{2,j}-\imath\g_{11})\,
  \prod\limits_{i=1}^3\Gamma(\imath\g_i-\imath\g_{2,j})\,,
 \ee
where the integration domain $\CC$ is subjected to the following
conditions
 \be\label{gl3contour}
  {\rm Im}(\g_{11})\,
  >\,\max\bigl\{{\rm Im}(\g_{21}),\,{\rm Im}(\g_{22})\bigr\}\,,\\
  {\rm Im}(\g_{2j})\,>\,\max\bigl\{{\rm Im}(\g_1),\,{\rm Im}(\g_2),\,{\rm
  Im}(\g_3)\bigr\}\,,\quad j=1,2\,.
 \ee
We specify the representation parameters
$\g=(\g_1,\g_2,\g_3)\in\IC^3$ in terms of
$\la=(\la_1,\la_2,\la_3)\in\IR^3$ as follows:
 \be\label{gl3lambda}
  \g\,:=\,\la+\imath\ve\a_1\,,\quad\a_1=(1,-1,0),
 \ee
for some fixed $\ve$ such that $0<\ve<\frac{1}{2}$.
 Then \eqref{MBgl3int0} takes the following form:
 \be\label{MBgl3int1}
  \Psi_{\la_1+\imath\ve,\,\la_2-\imath\ve,\la_3}(e^{x_1},e^{x_2},e^{x_3})\,
  =\,\frac{e^{\imath(\la_1+\la_2+\la_3)x_3}}{(2\pi)^3}
  \!\int\limits_{\CC}\!
  \frac{d\g_{11}d\g_{21}d\g_{22}}
  {\Gamma(\imath\g_{21}-\imath\g_{22})\,
  \Gamma(\imath\g_{22}-\imath\g_{21})}\\
  \times\,e^{\imath(\g_{21}(x_2-x_3)+\g_{22})+\imath\g_{11}(x_1-x_2)}\,
  \prod_{j=1}^2\Gamma(\imath\g_{21}-\imath\g_{11})\\
  \times\prod\limits_{j=1}^2\Gamma\bigl(\imath\la_1-\imath\g_{2,j}-\ve\bigr)\,
  \Gamma\bigl(\imath\la_2-\imath\g_{2,j}+\ve\bigr)\,
  \Gamma\bigl(\imath\la_3-\imath\g_{2,j}\bigr)\,,
 \ee
where the integration domain $\CC$ is defined by
 \be\label{gl3contour11}
  {\rm Im}(\g_{11})>\max\bigl\{{\rm Im}(\g_{21}),\,{\rm Im}(\g_{22})\bigr\},
  \qquad
  {\rm Im}(\g_{2j})>\ve,\quad j=1,2\,.
 \ee
By \eqref{GTcontour}, we specify the integration domain $\CC$ to be
parameterized by  $\tau_{k,i}\in\IR$ as follows
\be\label{GTgl3contour}
  \g_{11}=\tau_{11}\,+\,\imath,\quad
  \g_{2,j}=\tau_{2,j}+\frac{\imath}{2}\,,\quad j=1,2\,,
 \ee
so that after substitution, \eqref{MBgl3int1} becomes
 \be\label{MBgl3int}
  \Psi_{\la_1+\imath\ve,\,\la_2-\imath\ve,\la_3}
  (e^{x_1},e^{x_2},e^{x_3})\,
  =\,\frac{e^{\imath(\la_1+\la_2+\la_3)x_3-x_1+x_3}}{(2\pi)^3}
  \!\int\limits_{\IR^3}\!
  d\tau_{11}d\tau_{21}d\tau_{22}\\
  \times\,\frac{e^{\imath(\tau_{21}+\tau_{22})(x_2-x_3)+\imath\tau_{11}(x_1-x_2)}}
  {\Gamma(-\imath\tau_{21}+\imath\tau_{22})\,
  \Gamma(\imath\tau_{21}-\imath\tau_{22})}\,
  \prod_{j=1}^2\Gamma\Big(\imath(\tau_{21}-\tau_{11})+\frac{1}{2}\Big)\\
  \times\prod\limits_{j=1}^2\Gamma\Big(\imath(\la_1-\tau_{2,j})-\ve+\frac{1}{2}\Big)\,
  \Gamma\Big(\imath\la_2-\tau_{2,j})+\ve+\frac{1}{2}\Big)\,
  \Gamma\Bigl(\imath\la_3-\tau_{2,j})+\frac{1}{2}\Big)\,.
 \ee

\subsection{The Gelfand-Tsetlin realization}

By \eqref{GTrep} for $\ell=2$, the Gelfand-Tsetlin realization of
the principal series $\CU(\gl_3)$-representa-tion
$(\pi_{-\imath\g-\rho},\,\CV_{-\imath\g-\rho})$ with
$\g=\la+\imath\ve\a_1,\,\la\in\IR^3$ reads
 \be\label{GTgl3rep}
  E_{11}=-\imath\tau_{11},\qquad
  E_{22}=-\imath(\tau_{21}+\tau_{22})+\imath\tau_{11},\\
  E_{33}=-\imath(\g_1+\g_2+\g_3)+\imath(\tau_{21}+\tau_{22}),\\
  E_{12}=-\imath\Big(\imath(\tau_{11}-\tau_{21})+\frac{1}{2}\Big)
  \Big(\imath(\tau_{11}-\tau_{22})+\frac{1}{2}\Big)\,e^{-\imath\pr_{\tau_{11}}},\\
  E_{23}
  =-\imath\Big\{\frac{\prod\limits_{i=1}^3\bigl(\imath\tau_{21}-\imath\g_i+\frac{1}{2}\bigr)}
  {\imath\tau_{21}-\imath\tau_{22}}\,e^{-\imath\pr_{\tau_{21}}}\,
  +\,\frac{\prod\limits_{j=1}^3\bigl(\imath\tau_{22}-\imath\g_i+\frac{1}{2}\bigr)}
  {\imath\tau_{22}-\imath\tau_{21}}\,e^{-\imath\pr_{\tau_{22}}}\Big\}\,,\\
  E_{21}=-\imath e^{\imath\pr_{\tau_{11}}}\,,\qquad
  E_{32}=\imath\Big\{\frac{\imath\tau_{11}-\imath\tau_{21}+\frac{1}{2}}
  {\imath\tau_{21}-\imath\tau_{22}}\,e^{\imath\pr_{\tau_{21}}}\,
  +\,\frac{\imath\tau_{11}-\imath\tau_{22}+\frac{1}{2}}
  {\imath\tau_{22}-\imath\tau_{21}}\,e^{\imath\pr_{\tau_{22}}}\Big\}\,.
 \ee
The defining equations \eqref{whit} on Whittaker vectors allow the
solutions \eqref{GTwhitL},\eqref{GTwhitR}:
 \be\label{GTgl3vec}
  \psi_L(\tau)\,=\,e^{-\frac{\pi}{2}(\tau_{11}+\tau_{21}+\tau_{22})}\,,\\
  \psi_R(\tau)=e^{\frac{\pi}{2}(\tau_{11}-\tau_{21}-\tau_{22})}
  \prod_{j=1}^2\Gamma\Big(\imath(\tau_{2,j}-\tau_{11})+\frac{1}{2}\Big)
  \prod_{i=1}^3\Gamma\Big(\imath\g_i-\tau_{21})+\frac{1}{2}\Big)\,.
 \ee
Taking into account the  measure \eqref{GTmeasure} entering the
pairing \eqref{GTpair},
 \be\label{GTgl3Measure}
  \mu(\tau)\,
  =\,\frac{1}{(2\pi)^3}\,\frac{e^{\pi(\tau_{21}+\tau_{22})}}
  {\Gamma(\imath\tau_{21}-\imath\tau_{22})\,
  \Gamma(\imath\tau_{22}-\imath\tau_{21})}\,,
 \ee
and substituting \eqref{GTgl3vec}, we identify \eqref{MBgl3int} for
$\g=\la+\imath\ve\alpha_1$ with the matrix element as follows:
 \be\label{GTgl3matel0}
  \Psi_{\g_1,\g_2,\g_3}
  (e^{x_1},e^{x_2},e^{x_3})\,
  =\,e^{-\rho(x)}\<\psi_L,\,e^{-\sum\,x_iE_{ii}}\psi_R\>_{GT}\\
  =\,e^{-\rho(x)}\!\int\limits_{\IR^3}\!d\tau\,
  \mu(\tau)\,\ov{\psi_L(\tau)}\,
  e^{-\sum\,x_iE_{ii}}\psi_R(\tau)\,,\qquad\rho(x)=x_1-x_3\,.
 \ee
To make an identification of the two realizations of principle
series
  representations of $\CU(\mathfrak{gl}_3)$ and of the
  corresponding integral expressions for matrix elements, it is useful to rewrite
  the Gelfand-Tsetlin realization in a bit different form.

First, we transform the pairing \eqref{GTpair} with a non-trivial
measure \eqref{GTgl3Measure} into the standard $L^2$ one. Precisely,
let us write the integration measure function \eqref{GTgl3Measure}
as follows:
 \be
  \mu(\tau)\,=\,\mu_1(\tau)\,\ov{\mu_1(\tau)}\,,\qquad
  \mu_1(\tau)\,
  =\,\frac{1}{(2\pi)^{3/2}}\,\frac{e^{\frac{\pi}{2}(\tau_{21}+\tau_{22})}}
  {\Gamma(\imath\tau_{21}-\imath\tau_{22})}\,,
 \ee
and  consider  the following modification of the Gelfand-Tsetlin
realization \eqref{GTgl3rep}:
 \be\label{GTgl3S}
  \wt{E}_{ij}\,:=\,\mu_1E_{ij}\,\mu^{-1}_1\,,\qquad
  1\leq i,j\leq2\,\,.
 \ee
The modified difference operators \eqref{GTgl3S} have the following  form:
 \be\label{GTgl3Srep}
  \wt{E}_{11}=-\imath\tau_{11},\qquad
  \wt{E}_{22}=-\imath(\tau_{21}+\tau_{22})+\imath\tau_{11},\\
  \wt{E}_{33}=-\imath(\la_1+\la_2+\la_3)+\imath(\tau_{21}+\tau_{22}),\\
  \wt{E}_{12}=-\imath\Big(\imath(\tau_{11}-\tau_{21})+\frac{1}{2}\Big)
  \Big(\imath(\tau_{11}-\tau_{22})+\frac{1}{2}\Big)\,e^{-\imath\pr_{\tau_{11}}}\,,\qquad
  \wt{E}_{21}=-\imath\,e^{\imath\pr_{\tau_{11}}}\,,\\
  \wt{E}_{23}=\prod_{i=1}^3\Big(\imath(\tau_{21}-\g_i)+\frac{1}{2}\Big)\,
  e^{-\imath\pr_{\tau_{21}}}
  +\frac{\prod\limits_{i=1}^3\bigl(\imath(\tau_{22}-\g_i)+\frac{1}{2}\bigr)}
  {\imath(\tau_{22}-\tau_{21})[\imath(\tau_{21}-\tau_{22})-1]}\,
  e^{-\imath\pr_{\tau_{22}}}\,,\\
  \wt{E}_{32}=\frac{\imath(\tau_{11}-\tau_{21})+\frac{1}{2}}
  {\imath(\tau_{21}-\tau_{22})[\imath(\tau_{21}-\tau_{22})-1]}\,
  e^{\imath\pr_{\tau_{21}}}
  -\Big(\imath(\tau_{11}-\tau_{22})+\frac{1}{2}\Big)\,
  e^{\imath\pr_{\tau_{22}}}\,.
 \ee
The corresponding  Whittaker vectors \eqref{GTgl3vec} are given by
 \be\label{GTgl3whit}
  \wt{\psi}_L(\tau)
  =\mu_1(\tau)\,\psi_L(\tau)
  =\,\frac{1}{(2\pi)^{3/2}}\,
  \frac{e^{-\frac{\pi\,\tau_{11}}{2}}}
  {\Gamma\bigl(\imath\tau_{21}-\imath\tau_{22}\bigr)}\,,\\
  \wt{\psi}_R(\tau)
  =\mu_1(\tau)\,\psi_R(\tau)\,
  =\,\frac{1}{(2\pi)^{3/2}}\,
  \frac{e^{\frac{\pi\,\tau_{11}}{2}}}
  {\Gamma\bigl(\imath\tau_{21}-\imath\tau_{22}\bigr)}
  \prod_{j=1}^2\Gamma\Big(\imath(\tau_{2,j}-\tau_{11})+\frac{1}{2}\Big)\\
  \times\prod_{i=1}^3 \prod_{j=1}^2\Gamma\Big(\imath(\g_i-\tau_{2,j})+\frac{1}{2}\Big)\,.
  \ee
The matrix element \eqref{GTgl3matel0} in the considered realization
now takes the following form
 \be\label{GTgl3Smatel1}
  \Psi_{\g_1,\g_2,\g_3}(e^{x_1},e^{x_2},e^{x_3})\,
  =\,e^{-\rho(x)}\!\!\int_{\IR}\!d\tau\,\ov{\wt{\psi}_L(\tau)}\,
  e^{-\sum\,x_i\wt{E}_{ii}}\wt{\psi}_R(\tau)\,.
 \ee

Second, we extend the matrix element representation to encompass the
general integral representation \eqref{MBgl3int0} for the case of
$\mathfrak{gl}_3(\IR)$. Similarly to \eqref{DefPAIR}, for
$\vk\in\IR$, we introduce the following deformation of the standard
$L^2$-pairing, for $\psi_1\in\CV_{-\imath\g-\rho}^{\vee}$ and
$\psi_2\in\CV_{-\imath\g-\rho}$:
 \be\label{gl3DefPAIR}
  \<\,\,,\,\>_{\vk}\,:\quad\CV^{\vee}_{-\imath\g-\rho}
  \times\CV_{-\imath\g-\rho}\,\longrightarrow\,\IC\,,\\
  \<\psi_1,\psi_2\>_{\vk}\,:
  =\!\int\limits_{\IR-\imath\vk}\!\!d\tau_{11}
  \!\int\limits_{\IR^2}\!d\tau_{21}d\tau_{22}\,
  \bar{\psi_1}(\tau_{11},\tau_{21},\tau_{22})\,
  \psi_2(\tau_{11},\tau_{21},\tau_{22})\,,
 \ee
which might be rewritten as
 \be
  \<\psi_1,\psi_2\>_{\vk}\,
  =\int\limits_{\IR^3}\!d\tau\,
  \bar{\psi_1}(\tau-\imath\vk)\,\psi_2(\tau-\imath\vk)\,
  =\!\int\limits_{\IR^3}\!d\tau\,
  \overline{\psi_1(\tau+\imath\vk)}\,\psi_2(\tau-\imath\vk)\,.
 \ee
Given the  Gelfand-Tsetlin realization \eqref{GTgl3Srep}, consider
the operators $\wt{E}_{ij}$ with shifted argument $\tau_{11}\mapsto
(\tau_{11}-\imath\vk)$:
 \be\label{GTgl3repk}
  \KE_{ij}(\tau)\,:
  =\,\wt{E}_{ij}(\tau_{11}-\imath\vk,\,\tau_{21},\tau_{22})\,,\qquad
  1\leq i,j\leq3\,.
 \ee
\begin{lem}\label{gl3GTcontour}
Let $\wt{\psi}_L,\,\wt{\psi}_R$ be the Whittaker vectors given by
\eqref{GTgl3whit}. Then for $\vk\in\IR$ subjected $\vk<\frac{1}{2}$,
the following holds:
 \be\label{gl3IDINT}
  \<\wt{\psi}_L,\,e^{-\sum\,x_i\KE_{ii}}\wt{\psi}_R\>_{\vk}\,
  =\,\<\wt{\psi}_L,\,e^{-\sum\,x_i\wt{E}_{ii}}\wt{\psi}_R\>_0\,.
 \ee
\end{lem}
\proof Analogous to the proof of Lemma \ref{gl2GTcontour} in the
$\mathfrak{gl}_2(\IR)$-case. $\Box$

Thus, for $\g_1,\g_2,\g_3)=(\la_1+\imath\ve,\la_2-\imath\ve,\la_3)$
with $\la\in\IR^3$ and for $\vk<\frac{1}{2}$ using Lemma
\ref{gl3GTcontour}, we arrive at the  following matrix element
representation for the $\gl_3(\IR)$-Whittaker function:
 \be\label{GTgl3matel}
  \Psi_{\la_1+\imath\ve,\,\la_2-\imath\ve,\la_3}
  (e^{x_1},e^{x_2},e^{x_3})\,
  =\,e^{-\rho(x)}\<\wt{\psi}_L,\,e^{-\sum\,x_i\KE_{ii}}\wt{\psi}_R\>_{\vk}\\
  =\,\frac{e^{\imath(\la_1+\la_2+\la_3)x_3+\vk(x_1-x_2)-x_1+x_3}}{(2\pi)^3}\,
  \!\int\limits_{\IR^3}\!\frac{d\tau_{11}d\tau_{21}d\tau_{22}}
  {\Gamma(-\imath\tau_{21}+\imath\tau_{22})\,
  \Gamma(\imath\tau_{21}-\imath\tau_{22})}\\
  \times\,e^{\imath(\tau_{21}+\tau_{22})(x_2-x_3)+\imath\tau_{11}(x_1-x_2)}\,
  \prod_{j=1}^2\Gamma\Big(\imath(\tau_{21}-\tau_{11})-\vk+\frac{1}{2}\Big)\\
  \times\prod\limits_{j=1}^2\Gamma\Big(\imath(\la_1-\tau_{2,j})+\ve+\frac{1}{2}\Big)\,
  \Gamma\Big(\imath(\la_2-\tau_{2,j})-\ve+\frac{1}{2}\Big)\,
  \Gamma\Bigl(\imath(\la_3-\tau_{2,j})+\frac{1}{2}\Big)\,.
 \ee
This  reproduces a more general integral expression for the matrix
element \eqref{MBgl3int} in representation theory framework.

\subsection{The Gauss-Givental realization}

The analogous construction for the Gauss-Givental realization of the
  $\mathfrak{gl}_3(\IR)$-Whittaker matrix elements goes as follows.
For $\g=(\g_1,\g_2,\g_3)\in\IC^3$, consider the principal series
$\CU(\gl_3)$-representation
$(\pi_{-\imath\g-\rho},\,\CV_{-\imath\g-\rho})$. For some fixed
$\ve$ such that $0<\ve<\frac{1}{2}$ and for
 \be
  \g\,=\,\la+\imath\ve\a_1\,,\quad
  (\g_1,\g_2,\g_3)\,=\,(\la_1+\imath\ve,\la_2-\imath\ve,\la_3)\,,\quad
  \la\in\IR^3\,,
 \ee
the modified  Gauss-Givental realization of the
$\CU(\gl_3)$-representation is given by specialization of the
  expressions \eqref{GGrepF} to $\ell=2$
 \be\label{GGgl3repF}
  \hat{\CE}_{11}=-\imath(\la_1+s_{11}+s_{21})+\ve,\qquad
  \hat{\CE}_{22}=\imath(-\la_2+s_{21}-s_{22})-\ve,\\
  \hat{\CE}_{33}=\imath(-\la_3+s_{11}+s_{22}),\qquad
  \hat{\CE}_{12}=\Big(\imath s_{21}+\frac{1}{2}\Big)e^{-\imath\pr_{s_{21}}},\\
  \hat{\CE}_{23}
  =\Big(\imath s_{11}+\frac{1}{2}\Big)e^{\imath\pr_{s_{21}}-\imath\pr_{s_{11}}}\,
  +\,\Big(\imath(s_{11}-s_{21}+s_{22})+\frac{1}{2}\Big)e^{-\imath\pr_{s_{22}}},\\
  \hat{\CE}_{21}
  =\Big(\imath(\la_2-\la_1-s_{11})+2\ve+\frac{1}{2}\Big)\,e^{\imath\pr_{s_{11}}-\imath\pr_{s_{22}}}\\
  +\,\Big(\imath(\la_2-\la_1-s_{11}-s_{21}+s_{22})+2\ve+\frac{1}{2}\Big)\,e^{\imath\pr_{s_{21}}},\\
  \hat{\CE}_{32}
  =\Big(\imath(\la_3-\la_2-s_{22})-\ve+\frac{1}{2}\Big)\,e^{\imath\pr_{s_{22}}}\,,\\
 \ee
The corresponding  Whittaker vectors read from
\eqref{GGwhitLF},\eqref{GGwhitRF} for $\ell=2$:
 \be\label{GGgl3Fwhit}
 \begin{array}{lc}
  \hat{\phi}_L(s)\,
  =\,&\frac{1}{(2\pi)^{3/2}}\,\Gamma\Big(\imath(\bar{\g}_1-\bar{\g}_2+s_{11})+\frac{1}{2}\Big)\,
  \Gamma\Big(\imath(\bar{\g}_2-\bar{\g}_3+s_{22})+\frac{1}{2}\Big)\\
  &\times\,\Gamma\Big(\imath(\bar{\g}_1-\bar{\g}_3+s_{11}+s_{21})+1\Big)\,,\\
  \hat{\phi}_R(s)\,
  =\,&\frac{1}{(2\pi)^{3/2}}\,\Gamma\Big(\frac{1}{2}-\imath s_{21}\Big)\,
  \Gamma\Big(\frac{1}{2}-\imath s_{11}\Big)\,
  \Gamma\bigl(1-\imath s_{11}-\imath s_{22}\bigr)\,.
 \end{array}
 \ee
Then the  modified Gauss-Givental realization \eqref{GGintF} of the
$\mathfrak{gl}_3(\IR)$-Whittaker function reads:
 \be\label{GGgl3intF}
  \Psi_{\g_1,\g_2,\g_3}
  (e^{x_1},e^{x_2},e^{x_3})\,
  =\,e^{-\rho(x)}\<\hat{\phi}_L,\,e^{-x_1\hat{\CE}_{11}-x_1\hat{\CE}_{22}-x_1\hat{\CE}_{33}}
  \hat{\phi}_R\>\\
  =\,\frac{e^{\imath(\g_1x_1+\g_2x_2+\g_3x_3)-x_1+x_3}}{(2\pi)^3}\!
  \int\limits_{\IR^3}\!\prod_{k,i}ds_{k,i}\,
  e^{\imath s_{21}(x_1-x_2)+\imath s_{11}(x_1-x_3)
  +\imath s_{22}(x_2-x_3)}\\
  \times\Gamma\Big(\frac{1}{2}-\imath s_{21}\Big)\,
  \Gamma\Big(\frac{1}{2}-\imath s_{11}\Big)\,
  \Gamma\bigl(1-\imath s_{11}-\imath s_{22}\bigr)\\
  \times\Gamma\Big(-\imath(\g_1-\g_3+s_{11}+s_{21})+1\Big)\\
  \times\Gamma\Big(-\imath(\g_2-\g_3+s_{22})+\frac{1}{2}\Big)\,
  \Gamma\Big(-\imath(\g_1-\g_2+s_{11})+\frac{1}{2}\Big)\,.
 \ee
We identify the integral representation \eqref{GGgl3intF} with
\eqref{MBgl3int} via the first Barnes integral identity
\eqref{Barnes} in Proposition \ref{GGgl3MB} below.

\subsection{Intertwining operators}

For $0<\vk<\frac{1}{2}$, introduce operators $\CB_R$ and $\CB_L$
acting in $\CV_{-\imath\g-\rho}$ and
$\CV_{-\imath\g-\rho}^{\vee}\simeq\CV_{-\imath\bar{\g}-\rho}$,
respectively. Namely, the operators are defined by their integral
kernels (see \eqref{Delta} for the definition of delta-function in
complex argument):
 \be\label{gl3BR}
  B_R(\tau;s)\,
  =\,\frac{1}{2\pi}\frac{e^{\frac{\pi(\tau_{11}-\imath\vk)}{2}}}
  {\Gamma\bigl(\imath\tau_{21}-\imath\tau_{22}\bigr)}\,
  \frac{\Gamma\bigl(\imath(\g_2-\g_1-s_{11})+\frac{1}{2}\bigr)}
  {\Gamma(-\imath s_{21}+\frac{1}{2})}\\
  \times\prod_{j=1}^2\Gamma\Big(\imath(\tau_{11}-\tau_{2,j}-s_{21})+\vk\Big)\,
  \Gamma\Big(\imath(\tau_{2,j}-\tau_{11})-\vk+\frac{1}{2}\Big)\\
  \times\prod_{j=1}^2\Gamma\Big(\imath(\g_3-\tau_{2,j})+\frac{1}{2}\Big)\\
  \times\,\delta(\g_1+s_{11}+s_{21}-\tau_{11}+\imath\vk)\,
  \delta\bigl(\g_1+\g_2+s_{11}+s_{22}-\tau_{21}-\tau_{22}\bigr)\,.
 \ee
and
 \be\label{gl3BL}
  B_L(\tau;s)\,
  =\,\frac{1}{2\pi}\frac{e^{-\frac{\pi(\tau_{11}-\imath\vk)}{2}}}
  {\Gamma\bigl(\imath\tau_{21}-\imath\tau_{22}\bigr)}\,
  \frac{\Gamma\bigl(\imath s_{21}+\frac{1}{2}\bigr)}
  {\Gamma\bigl(\imath(\bar{\g}_1-\bar{\g}_2+s_{11})+\frac{1}{2}\bigr)}\\
  \times\prod_{j=1}^2\frac{\Gamma\bigl(\imath(\tau_{11}-\tau_{2,j}-s_{21})+\vk\bigr)}
  {\Gamma\big(\imath(\tau_{11}-\tau_{2,j})+\vk+\frac{1}{2}\bigr)\,
  \Gamma\bigl(\imath(\tau_{2,j}-\g_3)+\frac{1}{2}\bigr)}\\
  \times\,\delta(\bar{\g}_1+s_{11}+s_{21}-\tau_{11}+\imath\vk)\,
  \delta\bigl(\bar{\g}_1+\bar{\g}_2+s_{11}+s_{22}-\tau_{21}-\tau_{22}\bigr)\,.
 \ee
\begin{lem} Given the  Gelfand-Tsetlin realization \eqref{GTgl3repk}
and the modified Gauss-Givental realization \eqref{GGgl3repF}, the
operators $\CB_L$ and $\CB_R$ satisfy the following intertwining
relations:
 \be\label{gl3intert}
  \KE_{ij}\circ\CB_R\,=\,\CB_R\circ\hat{\CE}_{ij}\,,\quad
  \KE_{ij}^{\vee}\circ\CB_L\,=\,\CB_L\circ\hat{\CE}^{\vee}_{ij}\,,\qquad
  1\leq i,j\leq3\,.
 \ee
\end{lem}
\proof See Appendix B. $\Box$

\begin{prop}\label{BRwhit}
The action of operator $\CB_R$ on $\hat{\phi}_R$ reads
 \be\label{gl3BwhitR}\phantom{\int}
  (\CB_R\cdot\hat{\phi}_R)(\tau)\,
  =\,\wt{\psi}_R(\tau_{11}-\imath\vk,\tau_{21},\tau_{22})\,
  =\,\frac{1}{(2\pi)^{3/2}}\,\frac{e^{\frac{\pi(\tau_{11}-\imath\vk)}{2}}}
  {\Gamma\bigl(\imath\tau_{21}-\imath\tau_{22}\bigr)}\\
  \times\prod_{j=1}^2\Gamma\Big(\imath(\tau_{2,j}-\tau_{11})-\vk+\frac{1}{2}\Big)
  \prod_{i=1}^3\Gamma\Big(\imath(\g_i-\tau_{2,j})+\frac{1}{2}\Big)\,.
 \ee
\end{prop}
\proof By \eqref{GGgl3Fwhit}, we have
 \be
  \hat{\phi}_R(s)\,
  =\,\frac{1}{(2\pi)^{3/2}}\,\Gamma\Big(\frac{1}{2}-\imath s_{11}\Big)\,
  \Gamma\Big(\frac{1}{2}-\imath s_{21}\Big)\,
  \Gamma\bigl(1-\imath(s_{11}+s_{22})\bigr)\,.
 \ee
Then substituting \eqref{gl3BR} and canceling out the
$\Gamma\bigl(\frac{1}{2}-\imath s_{21}\bigr)$-factor, one finds out:
 \be
  (\CB_R\cdot\hat{\phi}_R)(\tau)
  =\,\frac{1}{(2\pi)^{3/2}}\!\int\limits_{\IR^3}\!ds\,B_R(\tau;s)\,
  \Gamma\Big(\frac{1}{2}-\imath s_{21}\Big)\\
  \times\,\Gamma\Big(\frac{1}{2}-\imath s_{11}\Big)\,
  \Gamma\bigl(1-\imath(s_{11}+s_{22})\bigr)
 \ee
 \be
  =\,
  \frac{e^{\frac{\pi(\tau_{11}-\imath\vk)}{2}}}
  {\Gamma\bigl(\imath\tau_{21}-\imath\tau_{22}\bigr)}
  \prod_{j=1}^2\Gamma\Big(\imath(\tau_{2,j}-\tau_{11})-\vk+\frac{1}{2}\Big)\,
  \Gamma\Big(\imath(\g_3-\tau_{2,j})+\frac{1}{2}\Big)\\
  \times\,\frac{1}{(2\pi)^{5/2}}\int\limits_{\IR^3}\!ds\,\,
  \Gamma\Big(\frac{1}{2}-\imath s_{11}\Big)\,
  \Gamma\Big(\imath(s_{21}+\g_2-\tau_{11})-\vk+\frac{1}{2}\Big)\\
  \times\,
  \Gamma\bigl(1-\imath(s_{11}+s_{22}\bigr)
  \prod_{j=1}^2\Gamma\Big(\imath(\tau_{11}-\tau_{2,j}-s_{21})+\vk\Big)\\
  \times\,\delta(\g_1+s_{11}+s_{21}-\tau_{11}+\imath\vk)\,
  \delta\bigl(\g_1+\g_2+s_{11}+s_{22}-\tau_{21}-\tau_{22}\bigr)
 \ee
 \be
  =\,\frac{1}{(2\pi)^{5/2}}\,
  \frac{e^{\frac{\pi(\tau_{11}-\imath\vk)}{2}}}
  {\Gamma\bigl(\imath\tau_{21}-\imath\tau_{22}\bigr)}
  \prod_{j=1}^2\Gamma\Big(\imath(\tau_{2,j}-\tau_{11})-\vk+\frac{1}{2}\Big)\,
  \Gamma\Big(\imath(\g_3-\tau_{2,j})+\frac{1}{2}\Big)\\
  \times\,
  \Gamma\bigl(\imath(\g_1+\g_2-\tau_{21}-\tau_{22})+1\bigr)\!
  \int\limits_{\IR}\!ds_{21}
  \prod_{i=1}^2\Gamma\Big(\imath(\g_i-\tau_{11}+s_{21})-\vk+\frac{1}{2}\Big)\\
  \times
  \prod_{j=1}^2\Gamma\Big(\imath(\tau_{11}-\tau_{2,j}-s_{21})+\vk\Big)\,.
 \ee
The integral above can be calculated via the First Barnes Identity
\eqref{Barnes} with the following specified parameters $a_j,\,b_i$:
 \be\label{gl3Rab}
  a_j=\imath(\tau_{11}-\tau_{2,j})+\vk\,,\quad
  b_i=\imath(\g_i-\tau_{11}+s_{21})+\frac{1}{2}-\vk\,,\quad
  i,j=1,2\,.
 \ee
For $\frac{1}{2}-\vk-\ve>0$, applying \eqref{Barnes} one derives
\eqref{gl3BwhitR}. $\Box$

\begin{prop}\label{BLwhit}
The action of $\CB_L$ on $\hat{\phi}_L$ reads
 \be\label{gl3BwhitL}
  (\CB_L\cdot\hat{\phi}_L)(\tau)\,
  =\,\wt{\psi}_L(\tau_{11}-\imath\vk,\tau_{21},\tau_{22})\,
  =\,\frac{1}{(2\pi)^{3/2}}\,
  \frac{e^{-\frac{\pi(\tau_{11}-\imath\vk)}{2}}}
  {\Gamma\bigl(\imath\tau_{21}-\imath\tau_{22}\bigr)}\,.
 \ee
\end{prop}
\proof Recall from \eqref{GGgl3Fwhit} the expression for the left
Whittaker vector,
 \be
  \hat{\phi}_L(s)\,
  =\,\frac{1}{(2\pi)^{3/2}}\,
  \Gamma\Big(\imath(\bar{\g}_2-\bar{\g}_3+s_{22})+\frac{1}{2}\Big)\,
  \Gamma\Big(\imath(\bar{\g}_1-\bar{\g}_2+s_{11})+\frac{1}{2}\Big)\\
  \times\,\Gamma\Big(\imath(\bar{\g}_1-\bar{\g}_3+s_{11}+s_{21})+1\Big)\,.
 \ee
Substituting \eqref{gl3BL} and canceling the
$\Gamma\bigl(\imath(\bar{\g}_1-\bar{\g}_2+s_{11})+\frac{1}{2}\bigr)$-factor,
we derive the following:
 \be
  (\CB_L\cdot\hat{\phi}_L)(\tau)\,
  =\,\frac{1}{(2\pi)^{3/2}}\!\int\limits_{\IR^3}\!ds\,B_L(\tau;s)\,
  \Gamma\Big(\imath(\bar{\g}_1-\bar{\g}_2+s_{11})+\frac{1}{2}\Big)\\
  \times\,
  \Gamma\Big(\imath(\bar{\g}_2-\bar{\g}_3+s_{22})+\frac{1}{2}\Big)\,
  \Gamma\Big(\imath(\bar{\g}_1-\bar{\g}_3+s_{11}+s_{21})+1\Big)
 \ee
 \be
  =\,\frac{1}{(2\pi)^{5/2}}\,
  \frac{e^{-\frac{\pi(\tau_{11}-\imath\vk)}{2}}}
  {\Gamma\bigl(\imath\tau_{21}-\imath\tau_{22}\bigr)}
  \prod_{j=1}^2\frac{1}
  {\Gamma\bigl(\imath(\tau_{11}-\tau_{2,j})+\vk+\frac{1}{2}\bigr)\,
  \Gamma\bigl(\imath(\tau_{2,j}-\g_3)+\frac{1}{2}\bigr)}\\
  \times\!\int\limits_{\IR^3}\!ds\,
  \Gamma\Big(\imath s_{21}+\frac{1}{2}\Big)\,
  \Gamma\Big(\imath(\bar{\g}_2-\bar{\g}_3+s_{22})+\frac{1}{2}\Big)\\
  \times\,
  \Gamma\Big(\imath(-\bar{\g}_3+\bar{\g}_1+s_{11}+s_{21})+1\Big)
  \prod_{j=1}^2
  \Gamma\Big(\imath(\tau_{11}-\imath\vk-s_{21}-\tau_{2,j})\Big)\\
  \times\,\delta(\bar{\g}_1+s_{11}+s_{21}-\tau_{11}+\imath\vk)\,
  \delta\bigl(\bar{\g}_1+\bar{\g}_2+s_{11}+s_{22}-\tau_{21}-\tau_{22}\bigr)
 \ee
 \be
  =\,\frac{1}{(2\pi)^{5/2}}\,
  \frac{e^{-\frac{\pi(\tau_{11}-\imath\vk)}{2}}}
  {\Gamma\bigl(\imath\tau_{21}-\imath\tau_{22}\bigr)}
  \frac{\Gamma\bigl(\imath(\tau_{11}-\bar{\g}_3)+\vk+1\bigr)}
  {\prod\limits_{j=1}^2\Gamma\bigl(\imath(\tau_{11}-\tau_{2,j})+\vk+\frac{1}{2}\bigr)\,
  \Gamma\bigl(\imath(\tau_{2,j}-\g_3)+\frac{1}{2}\bigr)}\\
  \times\!\int\limits_{\IR}\!ds_{11}\,
  \Gamma\Big(\imath(\tau_{11}-\bar{\g}_1-s_{11})+\vk+\frac{1}{2}\Big)\,
  \Gamma\Big(\imath(\tau_{21}+\tau_{22}-\bar{\g}_1-\bar{\g}_3-s_{11})+\frac{1}{2}\Big)\\
  \times\prod_{j=1}^2\Gamma\Big(\imath(\bar{\g}_1-\tau_{2,j}+s_{11})\Big)\,.
 \ee
The integral above can be calculated via the first Barnes identity
\eqref{Barnes} with the following specified parameters $a_i,\,b_j$:
 \be\label{gl3Rab11}
  a_1=\imath(\tau_{11}-\bar{\g}_1)+\vk+\frac{1}{2}
  =\imath(\tau_{11}-\la_1)+\vk+\frac{1}{2}-\ve\,,\\
  a_2=\imath(\tau_{21}+\tau_{22}-\bar{\g}_1-\bar{\g}_3)+\frac{1}{2}
  =\imath(\tau_{21}+\tau_{22}-\la_1-\la_3)+\frac{1}{2}-\ve\,,\\
  b_j=\imath(\bar{\g}_1-\tau_{2,j})=\imath(\la_1-\tau_{2,j})+\ve\,,\qquad j=1,2\,.
 \ee
Namely, for $\frac{1}{2}>\ve>0$ we apply \eqref{Barnes} and taking
into account $\bar{\g}_3=\g_3=\la_3\in\IR$ we deduce the assertion
\eqref{gl3BwhitL}. $\Box$

Let $\CB_L^{\dag}$ be the adjoint operator w.r.t.  the pairing
\eqref{gl3DefPAIR}, so its integral kernel reads
 \be\label{gl3BLdag}
  B_L^{\dag}(s;\tau)\,
  =\,\frac{1}{2\pi}\frac{e^{-\frac{\pi(\tau_{11}-\imath\vk)}{2}}}
  {\Gamma\bigl(-\imath\tau_{21}+\imath\tau_{22}\bigr)}\,
  \frac{\Gamma\bigl(-\imath s_{21}+\frac{1}{2}\bigr)}
  {\Gamma\bigl(-\imath(\g_1-\g_2+s_{11})+\frac{1}{2}\bigr)}\\
  \times\prod_{j=1}^2\frac{\Gamma\bigl(-\imath(\tau_{11}-\tau_{2,j}-s_{21})+\vk\bigr)}
  {\Gamma\big(-\imath(\tau_{11}-\tau_{2,j})-\vk+\frac{1}{2}\bigr)\,
  \Gamma\bigl(-\imath(\tau_{2,j}-\g_3)+\frac{1}{2}\bigr)}\\
  \times\,\delta(\g_1+s_{11}+s_{21}-\tau_{11}+\imath\vk)\,
  \delta\bigl(\g_1+\g_2+s_{11}+s_{22}-\tau_{21}-\tau_{22}\bigr)\,.
 \ee
\begin{lem}\label{BRLrels}
Given the  Gelfand-Tsetlin realization \eqref{GTgl3repk} and the
modified Gauss-Givental realization \eqref{GGgl3repF}, the operator
$\CB^{\dag}_L$ satisfies the following intertwining relations:
 \be\label{gl3intert1}
  \hat{\CE}_{ij}\circ\CB^{\dag}_L\,
  =\,\CB^{\dag}_L\circ\KE_{ij}\,,\quad1\leq i,j\leq3\,.
 \ee
\end{lem}
\proof See Appendix B, Section \ref{BLintert}. $\Box$

Let us calculate the $\CB_L^{\dag}\circ\CB_R$-action on the right
Whittaker vector $\hat{\phi}_R$.
\begin{lem}\label{BLBRwhit}
For the operator $\CB^{\dag}_L$ defined in \eqref{gl3BLdag}, the
following holds:
  \be\label{gl3BBwhitR}
  \phantom{\int}
  (\CB_L^{\dag}\circ\CB_R)\cdot\hat{\phi}_R\,
  =\,\hat{\phi}_R\,.
 \ee
\end{lem}
\proof Combining \eqref{gl3BLdag} and \eqref{gl3BwhitR}, we find out
 \be
 \CB_L^{\dag}\cdot
  (\CB_R\cdot\hat{\phi}_R)(s)\,
  =\,\bigl(\CB_L^{\dag}\cdot\wt{\psi}_R\bigr)(s)\,
  =\!\int\limits_{\IR^3}\!d\tau\,\CB^{\dag}_L(s;\tau)\,
  \wt{\psi}_R(\tau_{11}-\imath\vk,\,\tau_{21},\tau_{22})\\
  =\,\frac{1}{2\pi}\!\int\limits_{\IR^3}\!d\tau\,
  e^{-\frac{\pi(\tau_{11}-\imath\vk)}{2}}\,
  \frac{\wt{\psi}_R(\tau_{11}-\imath\vk,\,\tau_{21},\tau_{22})}
  {\Gamma\bigl(-\imath\tau_{21}+\imath\tau_{22}\bigr)}\,
  \frac{\Gamma\bigl(-\imath s_{21}+\frac{1}{2}\bigr)}
  {\Gamma\bigl(-\imath(\g_1-\g_2+s_{11})+\frac{1}{2}\bigr)}\\
  \times\prod_{j=1}^2\frac{\Gamma\bigl(-\imath(\tau_{11}-\tau_{2,j}-s_{21})+\vk\bigr)}
  {\Gamma\big(-\imath(\tau_{11}-\tau_{2,j})-\vk+\frac{1}{2}\bigr)\,
  \Gamma\bigl(-\imath(\tau_{2,j}-\g_3)+\frac{1}{2}\bigr)}\\
  \times\,\delta(\g_1+s_{11}+s_{21}-\tau_{11}+\imath\vk)\,
  \delta\bigl(\g_1+\g_2+s_{11}+s_{22}-\tau_{21}-\tau_{22}\bigr)\,.
 \ee
Substituting from \eqref{GTgl3whit}
 \be
  \wt{\psi}_R(\tau_{11}-\imath\vk)
  =\,\frac{1}{(2\pi)^{3/2}}\,\frac{e^{\frac{\pi(\tau_{11}-\imath\vk)}{2}}}
  {\Gamma\bigl(\imath\tau_{21}-\imath\tau_{22}\bigr)}\\
  \times\prod_{j=1}^2\Gamma\Big(\imath(\tau_{2,j}-\tau_{11})-\vk+\frac{1}{2}\Big)\,
  \prod_{i=1}^3\Gamma\Big(\imath(\g_i-\tau_{2,j})+\frac{1}{2}\Big)
 \ee
and canceling out the factor
$$
 e^{\frac{\pi(\tau_{11}-\imath\vk)}{2}}
 \prod_{j=1}^2\Gamma\Big(\imath(\tau_{2,j}-\tau_{11})-\vk+\frac{1}{2}\Big)\\
 \Gamma\Big(\imath(\g_3-\tau_{2,j})+\frac{1}{2}\Big)\,,
$$
we make the integration over $\tau_{11}$:
 \be
  \bigl(\CB_L^{\dag}\cdot\wt{\psi}_R\bigr)(s)\,
  =\,\frac{1}{(2\pi)^{5/2}}\,
  \frac{\Gamma\bigl(-\imath s_{21}+\frac{1}{2}\bigr)}
  {\Gamma\bigl(\imath(\g_2-\g_1-s_{11})+\frac{1}{2}\bigr)}
  \!\int\limits_{\IR^3}\!d\tau_{11}d\tau_{21}d\tau_{22}\\
  \times
  \frac{\prod\limits_{j=1}^2\Gamma\Big(\imath\tau_{2,j}+
  \imath(s_{21}-\tau_{11}+\imath\vk)\Big)
  \prod\limits_{i=1}^2\Gamma\Big(\imath(\g_i-\tau_{2,j})+\frac{1}{2}\Big)}
  {\Gamma\bigl(-\imath\tau_{21}+\imath\tau_{22}\bigr)\,
  \Gamma\bigl(\imath\tau_{21}-\imath\tau_{22}\bigr)}\\
  \times\,\delta(\g_1+s_{11}+s_{21}-\tau_{11}+\imath\vk)\,
  \delta\bigl(\g_1+\g_2+s_{11}+s_{22}-\tau_{21}-\tau_{22}\bigr)
 \ee
 \be
  =\,\frac{1}{(2\pi)^{5/2}}\,
  \frac{\Gamma\bigl(-\imath s_{21}+\frac{1}{2}\bigr)}
  {\Gamma\bigl(\imath(\g_2-\g_1-s_{11})+\frac{1}{2}\bigr)}
  \int\limits_{\IR^2}\!\frac{d\tau_{21}d\tau_{22}}
  {\Gamma\bigl(-\imath\tau_{21}+\imath\tau_{22}\bigr)\,
  \Gamma\bigl(\imath\tau_{21}-\imath\tau_{22}\bigr)}\\
  \times
  \prod_{j=1}^2\Gamma\Big(\imath\tau_{2,j}-\imath(s_{11}+\g_1)\Big)
  \prod_{i=1}^2\Gamma\Big(\imath(\g_i-\tau_{2,j})+\frac{1}{2}\Big)\\
  \times\,
  \delta\bigl(\g_1+\g_2+s_{11}+s_{22}-\tau_{21}-\tau_{22}\bigr)\,.
 \ee
Now we introduce variables
 \be
  \tau_{\pm}:=\frac{\tau_{21}\pm\tau_{22}}{2}\,,\qquad
  d\tau_{21}d\tau_{22}=d\tau_+d\tau_-\,,
 \ee
so that after substitution we integrate over $\tau_+$:
 \be
  \bigl((\CB_L)^{\dag}\cdot\wt{\psi}_R\bigr)(s)\,
  =\,\frac{1}{(2\pi)^{5/2}}\,
  \frac{\Gamma\bigl(-\imath s_{21}+\frac{1}{2}\bigr)}
  {\Gamma\bigl(\imath(\g_2-\g_1-s_{11})+\frac{1}{2}\bigr)}
  \int\limits_{\IR^2}\!\frac{d\tau_+d\tau_-}
  {\Gamma(-2\imath\tau_-)\,\Gamma(2\imath\tau_-)}\\
  \times
  \Gamma\Big(\imath(\tau_++\tau_-)-\imath(s_{11}+\g_1)\Big)\,
  \Gamma\Big(\imath(\tau_+-\tau_-)-\imath(s_{11}+\g_1)\Big)\\
  \times\prod_{i=1}^2\Gamma\Big(\imath(\g_i-\tau_+-\tau_-)+\frac{1}{2}\Big)\,
  \Gamma\Big(\imath(\g_i-\tau_++\tau_-)+\frac{1}{2}\Big)\\
  \times\,
  \delta\Big(\frac{\g_1+\g_2+s_{11}+s_{22}}{2}-\tau_+\Big)
 \ee
 \be
  =\,\frac{1}{(2\pi)^{5/2}}\,
  \frac{\Gamma\bigl(-\imath s_{21}+\frac{1}{2}\bigr)}
  {\Gamma\bigl(\imath(\g_2-\g_1-s_{11})+\frac{1}{2}\bigr)}
  \int\limits_{\IR}\!\frac{d\tau_-}
  {\Gamma(-2\imath\tau_-)\,\Gamma(2\imath\tau_-)}\\
  \times
  \Gamma\Big(\imath(\frac{\g_1+\g_2+s_{11}+s_{22}}{2}-\g_1-s_{11})+\imath\tau_-\Big)\\
  \times
  \Gamma\Big(\imath(\frac{\g_1+\g_2+s_{11}+s_{22}}{2}-\g_1-s_{11})-\imath\tau_-\Big)\\
  \times\prod_{i=1}^2\Gamma\Big(\imath(\g_i-\frac{\g_1+\g_2+s_{11}+s_{22}}{2}-\tau_-)+\frac{1}{2}\Big)\\
\hspace{-2cm}
  \times
  \Gamma\Big(\imath(\g_i-\frac{\g_1+\g_2+s_{11}+s_{22}}{2}+\tau_-)+\frac{1}{2}\Big)\,.
 \ee
We calculate the above integral via the identity
 (see \cite{GLO11}, Appendix D, page 187):
 \be\label{GLO11}
  \frac{1}{2\pi}\!\int\limits_{\IR}\!\!\frac{d\tau}
  {\Gamma(-2\imath\tau)\,\Gamma(2\imath\tau)}
  \prod_{i=1}^3\Gamma(a_i-\imath\tau)\,\Gamma(a_i+\imath\tau)\\
  =\prod_{i<j}\Gamma(a_i+a_j)\,,\qquad {\rm Re}(a_i)>0\,,
 \ee
with the parameters $a_i$ given by
 \be
  a_i\,=\,\imath\Big(\g_i-\frac{\g_1+\g_2+s_{11}+s_{22}}{2}\Big)+\frac{1}{2}\,,\quad i=1,2\,,\\
  a_3\,=\,\imath\Big(\frac{\g_1+\g_2+s_{11}+s_{22}}{2}-\g_1-s_{11}\Big)\,
  =\,\imath\frac{-\la_1+\la_2-s_{11}+s_{22}}{2}+\ve\,.
 \ee
Thus applying \eqref{GLO11} we deduce the assertion
\eqref{gl3BBwhitR}. $\Box$

\subsection{The equivalence of integral realizations of the matrix
  element}

Given the  operators $\CB_L,\,\CB_R$ defined by
\eqref{gl3BL},\eqref{gl3BR} and the operator $\CB^{\dag}_L$ defined
by \eqref{gl3BLdag}, we combine the constructions above in this
section to state and prove the main result.

\begin{te}\label{MAINTHM} Transformation of the left and right $\mathfrak{gl}_3(\IR)$-Whittaker vectors by the
  operators $\CB_L,\,\CB_R$ provides a direct
  identification of the integral forms   \eqref{GTgl3matel0} and
  \eqref{GGgl3intF}   of the matrix element arising
  via the Gelfand-Tsetlin and the modified Gauss-Givental
  realizations.
\end{te}
\proof Starting with \eqref{GTgl3Smatel1} and followed by
\eqref{GTgl3matel}, we apply Propositions \ref{BRwhit} and
\ref{BLwhit}:
 \be
  \Psi_{\g_1,\g_2,\g_3}(e^{x_1},e^{x_2},e^{x_3})\,
  =\,\<\psi_L,\,e^{-\sum x_iE_{ii}}\psi_R\>_{GT}\,
  =\,\<\wt{\psi}_L,\,e^{-\sum x_i\wt{E}_{ii}}\,\wt{\psi}_R\>\\
  =\,\<\wt{\psi}_L,\,e^{-\sum x_i\wt{E}_{ii}}\,\wt{\psi}_R\>_{\vk}\,
  =\,\<\CB_L\hat{\phi}_L,\,
  e^{-\sum x_i\KE_{ii}}\CB_R\hat{\phi}_R\>\,,
 \ee
where we use the notation \eqref{GTgl3repk}. Next, we note that
\eqref{gl3intert1} implies
 \be
  \CB_L^{\dag}\,e^{-\sum x_i\KE_{ii}}\bigl(\CB_L^{\dag}\bigr)^{-1}\,
  =\,e^{-\sum x_i\hat{\CE}_{ii}}\,.
 \ee
Therefore, combining this with  Lemma \ref{BLBRwhit} we arrive at
\eqref{GGgl3intF}:
 \be
  \Psi_{\g_1,\g_2,\g_3}(e^{x_1},e^{x_2},e^{x_3})\,
  =\,\<\hat{\phi}_L,\,
  \CB_L^{\dag}\,e^{-\sum x_i\KE_{ii}}\CB_R\hat{\phi}_R\>\\
  =\,\<\hat{\phi}_L,\,
  e^{-\sum x_i\hat{\CE}_{ii}}\CB_L^{\dag}\circ\CB_R\hat{\phi}_R\>\,
  =\,\<\hat{\phi}_L,\,
  e^{-\sum x_i\hat{\CE}_{ii}}\hat{\phi}_R\>\,.
 \ee
This completes our proof. $\Box$

Although we have completed our task to establish the desired
equivalences, in the following section we give  an alternative
identification of the Mellin-Barnes and the Givental integral
representations of the $\mathfrak{gl}_{\ell+1}(\IR)$-Whittaker
functions for $\ell=1,2$. This time we consider transformations of
one integral expressions into another  disregarding its matrix
element structure. In the process we encounter yet another
application of integral identity \eqref{Gustafson} mentioned in
Introduction.

\section{ Direct identification of two integral representations}

\subsection{The $\gl_2(\IR)$ case}

In this section we directly identify the Gelfand-Tsetlin and the
Gauss-Givental representations of the $\gl_2(\IR)$-Whittaker
function.

For $(\g_1,\g_2)\in\IC^2$ and $(x_1,x_2)\in\IR^2$, let
$\Psi_{\g_1,\g_2}(e^{x_1},e^{x_2})$ be the $\gl_2(\IR)$-Whittaker
function. Then it affords two integral representations: the Givental
integral representation,
 \be\label{gl2Giv}
  \Psi_{\g_1,\g_2}(e^{x_1},e^{x_2})\,
  =\!\int\limits_{\IR}\!dT\,
  e^{\imath\g_2(x_1+x_2-T)+\imath\g_1T\,-\,e^{x_1-T}-e^{T-x_2}}\,,
 \ee
and the Mellin-Barnes integral representation:
 \be\label{gl2MB}
  \Psi_{\g_1,\g_2}(e^{x_1},e^{x_2})\,
  =\,\frac{e^{\imath(\g_1+\g_2)x_2-\frac{x_1-x_2}{2}}}{2\pi}
  \!\int\limits_{\IR}\!\!d\tau\,
  e^{\imath\tau (x_1-x_2)}\,
  \prod_{j=1}^2\Gamma\Big(\imath\g_j-\imath \tau+\frac{1}{2}\Big)\,.
 \ee
The two integral representations, \eqref{gl2Giv} and \eqref{gl2MB},
may be explicitly identified via series of integral transformations.

\begin{prop} The integrals \eqref{gl2Giv} and
   \eqref{gl2MB} coincide.
\end{prop}

\proof To verify an equivalence of the two integral
  representations,
  we consider  a bit different  expression
  \be
  \Phi_{\g_1,\g_2}(e^{u})=e^{-\imath(\g_1+\g_2)x_2+\frac{x_1-x_2}{2}}
  \Psi_{\g_1,\g_2}(e^{x_1},e^{x_2})\,,\qquad u=x_1-x_2\,,
  \ee
and compare the two integral expressions depending on $u=x_1-x_2$:
 \be\label{gl2GivA}
  \Phi_{\g_1,\g_2}(e^{u})\,
  =e^{-\imath(\g_1+\g_2)+\frac{x_1-x_2}{2}}\,
    \!\int\limits_{\IR}\!dT\,
  e^{\imath(\g_1-\g_2)T\,-\,e^{x_1-T}-e^{T-x_2}}\,,
  \ee
and
 \be\label{gl2MBA}
  \Phi_{\g_1,\g_2}(e^{u})\,
  =\,\frac{1}{2\pi}\!\int\limits_{\IR}\!d\tau\,
  e^{\imath \tau  u}\,
  \prod_{j=1}^2\Gamma\Big(\imath\g_j-\imath \tau+\frac{1}{2}\Big)\,.
  \ee
After simple transformations, \eqref{gl2GivA} may be turned into the
following expression
 \be\label{gl2Givred}
  \Phi_{\g_1,\g_2}(e^u)\,
  =\,e^{(\imath\g_2+\frac{1}{2})u}\!\int\limits_{\IR}\!dT\,
    e^{\imath(\g_1-\g_2)T\,-\,e^{u-T}-e^{T}}\,.
 \ee
We proceed by applying the Fourier transform
$\CF:\,\CS(\IR_u)\to\CS(\IR_p)$ from \eqref{Fourier} for $d=1$ to
both expressions \eqref{gl2MBA} and \eqref{gl2Givred}. On the one
hand, the Fourier transform of  \eqref{gl2Givred} reads:
 \be
  \wh{\Phi}_{\g_1,\g_2}(p)\,
  =\,\frac{1}{\sqrt{2\pi}}\!\int\limits_{\IR}\!dT\,
  e^{\imath(\g_1-\g_2)T\,-\,e^{T}}\!\int\limits_{\IR}\!du\,
  e^{[\imath(\g_2-p)+\frac{1}{2}]u\,-\,e^{u-T}}\,.
  \ee
Changing the integration variable by $u':=u-T\,\in\IR$, after the
integration we obtain
 \be
  \wh{\Phi}_{\g_1,\g_2}(p)\,
  =\,\frac{1}{\sqrt{2\pi}}\,\Gamma\Big(\imath(\g_1-p)+\frac{1}{2}\Big)\,
  \Gamma\Big(\imath(\g_2-p)+\frac{1}{2}\Big)\,,
 \ee
On the other hand, the Fourier transform of  \eqref{gl2MBA} reads:
 \be\label{gl2Fourier}
  \wh{\Phi}_{\g_1,\g_2}(p)\,
  =\,\frac{1}{(2\pi)^{3/2}}\!\int\limits_{\IR}\!du\,e^{-\imath pu}\!
  \int\limits_{\IR}\!d\tau\,e^{\imath\tau u}\,
  \prod_{j=1}^2\Gamma\Big(\imath\g_j-\imath\tau+\frac{1}{2}\Big)\\
  =\,\frac{1}{\sqrt{2\pi}}\!\int\limits_{\IR}\!d\tau
  \prod_{j=1}^2\Gamma\Big(\imath\g_j-\imath\tau+\frac{1}{2}\Big)\,
  \delta(\tau-p)\\
 =\,\frac{1}{\sqrt{2\pi}}\,\Gamma\bigl(\imath(\g_1-p)+\frac{1}{2}\bigr)\,
  \Gamma\Big(\imath(\g_2-p)+\frac{1}{2}\Big)\,,
 \ee
where we use the integral representation of the standard
delta-function:
 \be\label{Deltaint}
  \delta(s-\tau)\,=\,\frac{1}{2\pi}\!\int\limits_{\IR}\!du\,
  e^{\imath(s-\tau)u}\,,\qquad s,\tau\in\IR\,.
 \ee
This establishes the equivalence claimed in Proposition. $\Box$

\subsection{The $\gl_3(\IR)$ case}

Now we compare the two  integral representations for the
$\mathfrak{gl}_3(\IR)$-Whittaker functions. Precisely we establish
an equivalence of their Fourier transforms by using the Gustafson
identity \eqref{Gustafson}.

The $\mathfrak{gl}_3(\IR)$-Whittaker function allows the Givental
integral representation \cite{Giv},\cite{JK}:
 \be\label{gl3Giv}
  \Psi_{\g_1,\g_2,\g_3}(e^{x_1},e^{x_2},e^{x_3})\,
  =\,e^{\imath\g_3(x_1+x_2+x_3)}\\
  \times\!\int\limits_{\IR^3}\!
  \prod_{k=1}^2\prod_{i=1}^kdT_{k,i}\,e^{\imath(\g_2-\g_3)(T_{21}+T_{22})
  +\imath(\g_1-\g_2)T_{11}}\\
  \times\exp\Big\{-e^{x_1-T_{21}}-e^{T_{21}-x_2}-e^{x_2-T_{22}}-e^{T_{22}-x_3}
  -e^{T_{21}-T_{11}}-e^{T_{11}-T_{22}}\Big\}\,.
 \ee
On the other hand, we have the Mellin-Barnes integral representation
\cite{KL}:
 \be\label{gl3MB}
  \Psi_{\g_1,\g_2,\g_3}(e^{x_1},e^{x_2},e^{x_3})\,
  =\,\frac{e^{\imath(\g_1+\g_2+\g_3)x_3-\rho(x)}}{(2\pi)^3}
  \!\int\limits_{\IR^3}\!
  \frac{d\tau_{11}d\tau_{21}d\tau_{22}}{\Gamma(\imath\tau_{21}-\imath\tau_{22})\,
  \Gamma(\imath\tau_{22}-\imath\tau_{21})}\\
  \times\,e^{\imath(x_2-x_3)(\tau_{21}+\tau_{22})+\imath(x_1-x_2)\tau_{11}}\\
  \times\Gamma\Big(\imath(\tau_{21}-\tau_{11})+\frac{1}{2}\Big)
  \Gamma\Big(\imath(\tau_{22}-\tau_{11})+\frac{1}{2}\Big)\,
  \prod\limits_{i=1}^3\prod_{j=1}^2\
  \Gamma\Big(\imath(\g_i-\tau_{2,j})+\frac{1}{2}\Big)\,.
 \ee
To identify the two integrals, let us introduce the following
normalized matrix element:
 \be\label{gl3NORMatel}
  \Phi_{\g_1,\g_2,\g_3}(e^{u_1},e^{u_2})\,:
  =\,e^{-\imath(\g_1+\g_2+\g_3)x_3+x_1-x_3}\times
  \Psi_{\g_1,\g_2,\g_3}(e^{x_1},e^{x_2},e^{x_3})\Big|_{u_1=x_1-x_2\atop
  u_2=x_2-x_3}\,,
 \ee
so that its Mellin-Barnes representation reads from \eqref{gl3MB}:
 \be\label{gl3MBR}
  \Phi_{\g_1,\g_2,\g_3}(e^{u_1},e^{u_2})\,
  =\,\frac{1}{(2\pi)^3}\!\int\limits_{\IR^3}\!
  \frac{d\tau_{11}d\tau_{21}d\tau_{22}}{\Gamma(\imath\tau_{21}-\imath\tau_{22})\,
  \Gamma(\imath\tau_{22}-\imath\tau_{21})}\,
  e^{\imath(\tau_{21}+\tau_{22})u_2+\imath\tau_{11}u_1}\\
  \times\Gamma\Big(\imath(\tau_{21}-\tau_{11})+\frac{1}{2}\Big)
  \Gamma\Big(\imath(\tau_{22}-\tau_{11})+\frac{1}{2}\Big)\,
  \prod\limits_{i=1}^3\prod_{j=1}^2\
  \Gamma\Big(\imath(\g_i-\tau_{2,j})+\frac{1}{2}\Big)\,.
 \ee
\begin{prop} The integrals \eqref{gl3Giv} and \eqref{gl3MB}
   coincide.
\end{prop}
\proof The Givental representation of \eqref{gl3NORMatel} reads from
\eqref{gl3Giv}:
 \be\label{gl3Giv1}
  \Phi_{\g_1,\g_2,\g_3}(e^{u_1},e^{u_2})\,
  =\,e^{-\imath(\g_1+\g_2+\g_3)x_3+x_1-x_3}
  e^{\imath\g_3(x_1+x_2+x_3)}\\
  \times\!\int\limits_{\IR^3}\!
  \prod_{k,i}dT_{k,i}\,e^{\imath(\g_2-\g_3)(T_{21}+T_{22})
  +\imath(\g_1-\g_2)T_{11}}\\
  \times\,e^{-e^{x_1-T_{21}}-e^{T_{21}-x_2}-e^{x_2-T_{22}}-e^{T_{22}-x_3}
  -e^{T_{21}-T_{11}}-e^{T_{11}-T_{22}}}\,.
 \ee
By changing the integration variables
 \be
  T'_{11}:=T_{11}-T_{22}\,,\qquad
  T'_{2,j}:=T_{2,j}-x_{j+1}\,,\quad j=1,2\,,
 \ee
the expression \eqref{gl3Giv1} may be transformed into the following
form
 \be\label{gl3Givred}
  \Phi_{\g_1,\g_2,\g_3}(e^{u_1},e^{u_2})\,
  =\,e^{(\imath\g_3+1)u_1+[\imath(\g_2+\g_3)+1]u_2}\\
  \times\!\int\limits_{\IR^3}\!
  \prod_{k,i}dT'_{k,i}\,e^{\imath(\g_1-\g_2)T'_{11}+
  \imath(\g_2-\g_3)T'_{21}+\imath(\g_1-\g_3)T'_{22}}\\
  \times\,e^{-e^{u_1-T'_{21}}-e^{T'_{21}}-e^{u_2-T'_{22}}-e^{T'_{22}}
  -e^{u_2+T'_{21}-T'_{22}-T'_{11}}-e^{T'_{11}}}\,.
 \ee
We identify the integrals \eqref{gl3Giv} and \eqref{gl3MB} by
identifying the Fourier transforms of \eqref{gl3MBR} and of
\eqref{gl3Givred}. Applying the Fourier transform
$\CF:\,\CS(\IR^2_u)\to\CS(\IR^2_p)$ from \eqref{Fourier} for $d=2$
to \eqref{gl3Givred} we obtain
 \be
  \widehat{\Phi}_{\g_1,\g_2,\g_3}(p_1,p_2)\,
  =\,\frac{1}{2\pi}\!\int\limits_{\IR^2}\!du_1du_2\,e^{-\imath p_1u_1-\imath  p_2u_2}\,
  \Phi^{\gl_3}_{\g}(u_1,u_2)\\
  =\,\frac{1}{2\pi}\!\int\limits_{\IR^2}\!du_1du_2\,e^{[\imath(\g_2+\g_3-p_2)+1]u_2+[\imath(\g_3-p_1)+1]u_1}\\
  \times\!\int\limits_{\IR^3}\!
  \prod_{k,i}dT'_{k,i}\,e^{\imath(\g_2-\g_3)T'_{21}+\imath(\g_1-\g_3)T'_{22}
  +\imath(\g_1-\g_2)T'_{11}}\\
  \times\,e^{-e^{u_1-T'_{21}}-e^{T'_{21}}-e^{u_2-T'_{22}}-e^{T'_{22}}
  -e^{u_2+T'_{21}-T'_{22}-T'_{11}}-e^{T'_{11}}}
 \ee
 \be\label{gl3GivF1}
  =\,\frac{1}{2\pi}\!\int\limits_{\IR^3}\!
  \prod_{k,i}dT'_{k,i}\,e^{\imath(\g_2-\g_3)T'_{21}+\imath(\g_1-\g_3)T'_{22}
  +\imath(\g_1-\g_2)T'_{11}\,-\,e^{T'_{21}}-e^{T'_{22}}-e^{T'_{11}}}\\
  \times\!\int\limits_{\IR}\!du_1\,e^{[\imath(\g_3-p_1)+1]u_1\,-\,e^{u_1-T'_{21}}}
  \!\int\limits_{\IR}\!du_2\,e^{[\imath(\g_2+\g_3-p_2)+1]u_2\,
  -\,e^{u_2-T'_{22}}\bigl(1+e^{T'_{21}-T'_{11}}\bigr)}\,.
 \ee
Introduce the following integration variables:
 \be
  u'_1\,:=u_1-T'_{21}\,\in\IR\,,\quad
  u'_2\,:=u_2-T'_{22}+\ln\bigl(1+e^{T'_{21}-T'_{11}}\bigr)\,
  \in\IR\,,
 \ee
so that after substitution \eqref{gl3GivF1} reads
 \be
  \widehat{\Phi}_{\g_1,\g_2,\g_3}(p_1,p_2)\,
  =\,\frac{1}{2\pi}\!\int\limits_{\IR^3}\!
  \prod_{k,i}dT'_{k,i}\,e^{\imath(\g_2-\g_3)T'_{21}+\imath(\g_1-\g_3)T'_{22}
  +\imath(\g_1-\g_2)T'_{11}\,-\,e^{T'_{21}}-e^{T'_{22}}-e^{T'_{11}}}\\
  \times\!\int\limits_{\IR}\!du'_1\,e^{[\imath(\g_3-p_1)+1](u'_1+T'_{21})\,-\,e^{u_1'}}
  \!\int\limits_{\IR}\!du'_2\,
  e^{[\imath(\g_2+\g_3-p_2)+1](u'_2+T'_{22}-\ln\bigl(1+e^{T'_{21}-T'_{11}}\bigr))\,-\,e^{u_2'}}\\
  =\,\frac{1}{2\pi}\!\int\limits_{\IR^3}\!
  \prod_{k,i}dT'_{k,i}\,\Big(1+e^{T'_{21}-T'_{11}}\Big)^{-[\imath(\g_2+\g_3-p_2)+1]}\\
  \times\,e^{[\imath(\g_2-p_1)+1]T'_{21}+[\imath(\g_1+\g_2-p_2)+1]T'_{22}
  +\imath(\g_1-\g_2)T'_{11}\,-\,e^{T'_{21}}-e^{T'_{22}}-e^{T'_{11}}}\\
  \times\!\int\limits_{\IR}\!du'_1\,e^{[\imath(\g_3-p_1)+1]u'_1\,-\,e^{u_1'}}
  \!\int\limits_{\IR}\!du'_2\,
  e^{[\imath(\g_2+\g_3-p_2)+1]u'_2\,-\,e^{u_2'}}\,.
 \ee
We integrate over $T'_{22},u'_1,u'_2$ applying  the Euler formula,
 \be\label{Euler}
  \Gamma(z)\,=\!\int\limits_{\IR}\!du\,e^{zu\,-\,e^u}\,,\quad{\rm
  Re}(z)>0\,,
 \ee
which results in
 \be\label{gl3FGiv1}
  \widehat{\Phi}^{\gl_3}_{\g}(p_1,p_2)\,
  =\,\frac{1}{2\pi}\,
  \Gamma\bigl(\imath(\g_1+\g_2-p_2)+1\bigr)\,
  \Gamma\bigl(\imath(\g_3-p_1)+1\bigr)\\
  \times\,\Gamma\bigl(\imath(\g_2+\g_3-p_2)+1\bigr)
  \!\int\limits_{\IR^2}\!dT'_{11}dT'_{21}\,
  \frac{e^{[\imath(\g_2-p_1)+1]T'_{21}+\imath(\g_1-\g_2)T'_{11}\,-\,e^{T'_{21}}-e^{T'_{11}}}}
  {\Big(1+e^{T'_{21}-T'_{11}}\Big)^{\imath(\g_2+\g_3-p_2)+1}}\,.
 \ee
Let us substitute into \eqref{gl3FGiv1} the following integration
variables:
 \be
 \tau_1\,:=T'_{21}-T'_{11}\,\in\IR\,,\qquad\tau_2\,:=T'_{11}\,\in\IR\,,
 \ee
which yields
 \be
  \widehat{\Phi}_{\g_1,\g_2,\g_3}(p_1,p_2)\,
  =\,\frac{1}{2\pi}\,\Gamma\bigl(\imath(\g_1+\g_2-p_2)+1\bigr)\,
  \Gamma\bigl(\imath(\g_3-p_1)+1\bigr)\\
  \times\,\Gamma\bigl(\imath(\g_2+\g_3-p_2)+1\bigr)
  \!\int\limits_{\IR^2}\!d\tau_1d\tau_2\,
  \frac{e^{[\imath(\g_2-p_1)+1](\tau_1+\tau_2)
  +\imath(\g_1-\g_2)\tau_2\,-\,e^{\tau_2}(1+e^{\tau_1})}}
  {\bigl(1+e^{\tau_1}\bigr)^{\imath(\g_2+\g_3-p_2)+1}}\\
  =\,\frac{1}{2\pi}\,\Gamma\bigl(\imath(\g_1+\g_2-p_2)+1\bigr)\,
  \Gamma\bigl(\imath(\g_3-p_1)+1\bigr)\\
  \times\,\Gamma\bigl(\imath(\g_2+\g_3-p_2)+1\bigr)
  \!\int\limits_{\IR^2}\!d\tau_1d\tau_2\,
  \frac{e^{[\imath(\g_2-p_1)+1]\tau_1
  +[\imath(\g_1-p_1)+1]\tau_2\,-\,e^{\tau_2}(1+e^{\tau_1})}}
  {\bigl(1+e^{\tau_1}\bigr)^{\imath(\g_2+\g_3-p_2)+1}}\,.
 \ee
Now we introduce $\tau'_2\,:=\tau_2+\ln(1+e^{\tau_1})\in\IR$ and
after substitution we obtain
 \be
  \widehat{\Phi}_{\g_1,\g_2,\g_3}(p_1,p_2)\,
  =\,\frac{1}{2\pi}\,\Gamma\bigl(\imath(\g_1+\g_2-p_2)+1\bigr)\,
  \Gamma\bigl(\imath(\g_3-p_1)+1\bigr)\\
  \times\,\Gamma\bigl(\imath(\g_2+\g_3-p_2)+1\bigr)
  \!\int\limits_{\IR^2}\!d\tau_1d\tau'_2\,
  \frac{e^{[\imath(\g_2-p_1)+1]\tau_1
  +[\imath(\g_1-p_1)+1](\tau'_2-\ln(1+e^{\tau_1}))\,-\,e^{\tau'_2}}}
  {\bigl(1+e^{\tau_1}\bigr)^{\imath(\g_2+\g_3-p_2)+1}}\\
  =\,\frac{1}{2\pi}\,\Gamma\bigl(\imath(\g_1+\g_2-p_2)+1\bigr)\,\Gamma\bigl(\imath(\g_2+\g_3-p_2)+1\bigr)\,
  \Gamma\bigl(\imath(\g_3-p_1)+1\bigr)\\
  \times\int\limits_{\IR}\!d\tau'_2\,e^{[\imath(\g_1-p_1)+1]\tau'_2\,-\,e^{\tau'_2}}
  \int\limits_{\IR}\!d\tau_1\,\frac{e^{[\imath(\g_2-p_1)+1]\tau_1}}
  {\bigl(1+e^{\tau_1}\bigr)^{\imath(\g_1+\g_2+\g_3-p_1-p_2)+2}}\\
  =\,\frac{1}{2\pi}\,\Gamma\bigl(\imath(\g_1+\g_2-p_2)+1\bigr)\,
  \Gamma\bigl(\imath(\g_2+\g_3-p_2)+1\bigr)\,
  \Gamma\bigl(\imath(\g_3-p_1)+1\bigr)\\
  \times\,\Gamma\bigl(\imath(\g_1-p_1)+1\bigr)\times\,B\bigl(\imath(\g_2-p_1)+1,\,
  \imath(\g_1+\g_3-p_2)+1\bigr)\,,
  \ee
where in the latter equality we apply the integral representation of
the Beta-function:
 \be\label{Beta}
  B(a,b)=\frac{\Gamma(a)\,\Gamma(b)}{\Gamma(a+b)}\,
  =\!\int\limits_{\IR}\!d\tau\,\frac{e^{a\tau}}{(1+e^{\tau})^{a+b}}\,,\qquad
  {\rm Re}(a),\,{\rm Re}(b)>0\,.
  \ee
Using expression of Beta-function in terms of Gamma-functions we
arrive at the explicit expression for the Fourier transform  of
\eqref{gl3Givred} \be\label{gl3GivFTA}
 \widehat{\Phi}_{\g_1,\g_2,\g_3}(p_1,p_2)\,
  =\,\frac{1}{2\pi}\,
  \frac{\prod\limits_{i=1}^3\Gamma\bigl(\imath(\g_i-p_1)+1\bigr)
  \prod\limits_{i<j}\Gamma\bigl(\imath(\g_i+\g_j-p_2)+1\bigr)}
  {\Gamma\bigl(\imath(\g_1+\g_2+\g_3-p_1-p_2)+2\bigr)}\,.
 \ee

The Fourier transform of the integral \eqref{gl3MBR} reads as
follows:
 \be
  \widehat{\Phi}_{\g_1,\g_2,\g_3}(p_1,p_2)\,
  =\,\frac{1}{2\pi}\!\int\limits_{\IR^2}\!du_1du_2\,
  e^{-\imath p_1u_1-\imath p_2u_2}\,\Phi^{\gl_3}_{\g}(e^{u_1},e^{u_2})\\
  =\,\frac{1}{(2\pi)^4}\!\int\limits_{\IR^3}\!\frac{d\tau_{11}d\tau_{21}d\tau_{22}}{\Gamma(\imath\tau_{21}-\imath\tau_{22})\,
  \Gamma(\imath\tau_{22}-\imath\tau_{21})}
  \prod\limits_{i=1}^3\prod_{j=1}^2\Gamma\Big(\imath(\g_i-\tau_{2,j})+\frac{1}{2}\Big)\\
  \times\Gamma\Big(\imath(\tau_{21}-\tau_{11})+\frac{1}{2}\Big)
  \Gamma\Big(\imath(\tau_{22}-\tau_{11})+\frac{1}{2}\Big)\!
  \int\limits_{\IR^2}\!du_1du_2\,e^{\imath(\tau_{21}+\tau_{22}-p_2)u_2+\imath(\tau_{11}-p_1)u_1}
 \ee
 \be
  =\,\frac{1}{(2\pi)^2}\!\int\limits_{\IR^3}\!\frac{d\tau_{11}d\tau_{21}d\tau_{22}}{\Gamma(\imath\tau_{21}-\imath\tau_{22})\,
  \Gamma(\imath\tau_{22}-\imath\tau_{21})}
  \prod\limits_{i=1}^3\prod_{j=1}^2\Gamma\Big(\imath(\g_i-\imath\tau_{2,j})+\frac{1}{2}\Big)\\
  \times\Gamma\Big(\imath(\tau_{21}-\imath\tau_{11})+\frac{1}{2}\Big)
  \Gamma\Big(\imath(\tau_{22}-\tau_{11})+\frac{1}{2}\Big)\,\delta(\tau_{11}-p_1)\,
  \delta(\tau_{21}+\tau_{22}-p_2)
 \ee
 \be
  =\,\frac{1}{(2\pi)^4}\!\int\limits_{\IR^2}\!\frac{d\tau_{21}d\tau_{22}}{\Gamma(\imath\tau_{21}-\imath\tau_{22})\,
  \Gamma(\imath\tau_{22}-\imath\tau_{21})}
  \prod\limits_{i=1}^3\prod_{j=1}^2\Gamma\Big(\imath(\g_i-\imath\tau_{2,j})+\frac{1}{2}\Big)\\
  \times\Gamma\Big(\imath(\tau_{21}-p_1)+\frac{1}{2}\Big)
  \Gamma\Big(\imath(\tau_{22}-p_1)+\frac{1}{2}\Big)\,\delta(\tau_{21}+\tau_{22}-p_2)
 \ee
where in the last equality we integrate out $\tau_{11}\in\IR$.
Introduce the integration variables,
 \be
  2\tau_+:=\tau_{21}+\tau_{22}\,\in\IR\,,\quad
  2\tau_-:=\tau_{21}-\tau_{22}\,\in\IR\,,\quad
  d\tau_+\,d\tau_-\,=\,d\tau_{21}d\tau_{22}\,.
 \ee
Then after substitution, and applying
$\delta(2p)=\frac{1}{2}\delta(p)$ the integral takes the following
form:
 \be
  \widehat{\Phi}_{\g_1,\g_2,\g_3}(p_1,p_2)\,
  =\,\frac{1}{2(2\pi)^2}\!\int\limits_{\IR}\!d\tau_+\!\int\limits_{\IR}\!
  \frac{d\tau_-}
  {\Gamma(2\imath\tau_-)\,\Gamma(-2\imath\tau_-)}\,
  \delta\Big(\tau_+-\frac{p_2}{2}\Big)\\
  \times
  \Gamma\Big(\imath(\tau_++\tau_-)-\imath p_1+\frac{1}{2}\Big)
  \Gamma\Big(\imath(\tau_+-\tau_-)-\imath p_1+\frac{1}{2}\Big)\\
  \times\prod\limits_{i=1}^3\Gamma\Big(\imath\g_i-\imath(\tau_++\tau_-)+\frac{1}{2}\Big)\,
  \Gamma\Big(\imath\g_i-\imath(\tau_+-\tau_-)+\frac{1}{2}\Big)
 \ee
Integrating out $\tau_+$ we arrive at the following expression for
Fourier transform of \eqref{gl3MBR}:
 \be\label{GL3MBFourier}
  \widehat{\Phi}_{\g_1,\g_2,\g_3}(p_1,p_2)\,
  =\,\frac{1}{2(2\pi)^2}\!\int\limits_{\IR}\!\frac{d\tau_-}
  {\Gamma(2\imath\tau_-)\,\Gamma(-2\imath\tau_-)}\\
  \times\,\Gamma\Big(\imath(\frac{p_2}{2}-p_1)+\imath\tau_- +\frac{1}{2}\Big)
  \Gamma\Big(\imath(\frac{p_2}{2}-p_1)-\imath\tau_- +\frac{1}{2}\Big)\\
  \times\prod\limits_{i=1}^3\Gamma\Big(\imath(\g_i-\frac{p_2}{2})+\imath\tau_-+\frac{1}{2}\Big)\,
  \Gamma\Big(\imath(\g_i-\frac{p_2}{2})-\imath\tau_-+\frac{1}{2}\Big)\,.
 \ee
To compute the integral, we use the Gustafson identity
\eqref{Gustafson} with the following parameters:
 \be
  a_i=\imath\Big(\g_i-\frac{p_2}{2}\Big)+\frac{1}{2},\quad
  i=1,2,3\,,\qquad a_4=\imath\Big(\frac{p_2}{2}-p_1\Big)+\frac{1}{2}\,.
 \ee
Namely, we arrive at
 \be\label{ID}
\hspace{-1.5cm}
 \begin{array}{c}
  \frac{1}{2(2\pi)^2}\!\int\limits_{\IR}\!\frac{d\tau}
  {\Gamma(2\imath\tau)\,\Gamma(-2\imath\tau)}
  \Gamma\Big(\imath(\frac{ p_2}{2}-p_1)+\imath\tau+\frac{1}{2}\Big)
  \Gamma\Big(\imath(\frac{p_2}{2}-p_1)-\imath\tau+\frac{1}{2}\Big)\\
  \times\prod\limits_{i=1}^3
  \Gamma\Big(\imath(\g_i-\frac{p_2}{2})-\imath\tau+\frac{1}{2}\Big)\,
  \Gamma\Big(\imath(\g_i-\frac{p_2}{2})+\imath\tau+\frac{1}{2}\Big)\\
  =\,\frac{1}{2\pi}\,\frac{\prod\limits_{i=1}^3\Gamma\bigl(\imath(\g_i-p_1)+1\bigr)
  \prod\limits_{i<j}\Gamma\bigl(\imath(\g_i+\g_j-p_2)+1\bigr)}
  {\Gamma\bigl(\imath(\g_1+\g_2+\g_3-p_1-p_2)+2\bigr)}\,,
 \end{array}
 \ee
which coincides with the Fourier transform \eqref{gl3GivFTA} of
\eqref{gl3Givred} and thus  completes the proof of  equivalence of
integral representations \eqref{gl3Giv} and \eqref{gl3MBR}. $\Box$

\begin{prop}\label{GGgl3MB} The integral \eqref{GGgl3intF} coincides
with \eqref{gl3MB}.
\end{prop}
\proof From \eqref{GGgl3intF}, the normalized matrix element
\eqref{gl3NORMatel} in the modified Gauss-Givental representation
reads
 \be\label{GGgl3MBred}
  \Phi_{\g_1,\g_2,\g_3}
  (e^{u_1},e^{u_2})\,
  =\,e^{-\imath(\g_1+\g_2+\g_3)x_3}\<\hat{\phi}_L,\,
  e^{-\sum x_i\hat{\CE}_{ii}}\hat{\phi}_R\>\Big|_{u_1=x_1-x_2\atop u_2=x_2-x_3}\\
  =\,\frac{1}{(2\pi)^3}\!\int\limits_{\IR^3}\!\prod_{k,i}ds_{k,i}\,
  e^{\imath s_{21}u_1+\imath(\g_1+s_{11})(u_1+u_2)
  +\imath(\g_2+s_{22})u_2}\\
  \times\Gamma\Big(\frac{1}{2}-\imath s_{21}\Big)\,
  \Gamma\Big(\frac{1}{2}-\imath s_{11}\Big)\,
  \Gamma\bigl(1-\imath s_{11}-\imath s_{22}\bigr)\\
  \times\Gamma\Big(-\imath(\g_1-\g_3+s_{11}+s_{21})+1\Big)\\
  \times\Gamma\Big(-\imath(\g_2-\g_3+s_{22})+\frac{1}{2}\Big)\,
  \Gamma\Big(-\imath(\g_1-\g_2+s_{11})+\frac{1}{2}\Big)\,.
 \ee
Then in straightforward way, the Fourier transform of
\eqref{GGgl3MBred} reads
 \be
  \wh{\Phi}_{\g_1,\g_2,\g_3}(p_1,p_2)\,
  =\,\frac{1}{2\pi}\!\int_{\IR^2}\!du_1du_2\,
  e^{-\imath p_1u_1-\imath p_2u_2}\,
  \Phi_{\g_1,\g_2,\g_3}(e^{u_1},e^{u_2})\\
  =\,\frac{1}{(2\pi)^4}\!\int\limits_{\IR^3}\!\prod_{k,i}ds_{k,i}
   \!\int_{\IR^3}\!du_1du_2\,
  e^{\imath(\g_1+s_{11}+s_{21}-p_1)u_1+\imath(\g_1+\g_2+s_{11}+s_{22}-p_2)u_2}\\
  \times\,
  \Gamma\Big(\frac{1}{2}-\imath s_{21}\Big)\,
  \Gamma\Big(\frac{1}{2}-\imath s_{11}\Big)\,
  \Gamma\bigl(1-\imath s_{11}-\imath s_{22}\bigr)\\
  \times\Gamma\Big(-\imath(\g_1-\g_3+s_{11}+s_{21})+1\Big)\\
  \times\Gamma\Big(-\imath(\g_2-\g_3+s_{22})+\frac{1}{2}\Big)\,
  \Gamma\Big(-\imath(\g_1-\g_2+s_{11})+\frac{1}{2}\Big)
 \ee
 \be
  =\,\frac{1}{(2\pi)^2}\!\int\limits_{\IR^3}\!\prod_{k,i}ds_{k,i}\,
  \Gamma\Big(\frac{1}{2}-\imath s_{21}\Big)\,
  \Gamma\Big(\frac{1}{2}-\imath s_{11}\Big)\,
  \Gamma\bigl(1-\imath s_{11}-\imath s_{22}\bigr)\\
  \times\Gamma\Big(-\imath(\g_1-\g_3+s_{11}+s_{21})+1\Big)\\
  \times\Gamma\Big(-\imath(\g_2-\g_3+s_{22})+\frac{1}{2}\Big)\,
  \Gamma\Big(-\imath(\g_1-\g_2+s_{11})+\frac{1}{2}\Big)\\
  \times\,\delta(\g_1+s_{11}+s_{21}-p_1)\,\delta(\g_1+\g_2+s_{11}+s_{22}-p_2)
 \ee
 \be
  =\,\frac{1}{(2\pi)^2}\,\Gamma\Big(\imath(\g_1+\g_2-p_2)+1\Big)\,
  \Gamma\Big(\imath(\g_3-p_1)+1\Big)\\
  \times\!\int\limits_{\IR^2}\!ds_{11}ds_{21}\,
  \delta(\g_1+s_{11}+s_{21}-p_1)\,
  \Gamma\Big(\frac{1}{2}-\imath s_{21}\Big)\,
  \Gamma\Big(\frac{1}{2}-\imath s_{11}\Big)\\
  \times\Gamma\Big(\imath(\g_1+\g_3+s_{11}-p_2)+\frac{1}{2}\Big)\,
  \Gamma\Big(-\imath(\g_1-\g_2+s_{11})+\frac{1}{2}\Big)
 \ee
 \be
  =\,\frac{1}{(2\pi)^2}\,\Gamma\Big(\imath(\g_1+\g_2-p_2)+1\Big)\,
  \Gamma\Big(\imath(\g_3-p_1)+1\Big)\\
  \times\!\int\limits_{\IR}\!ds_{11}\,
  \Gamma\Big(\imath(\g_1+s_{11}-p_1)+\frac{1}{2}\Big)\,
  \Gamma\Big(\frac{1}{2}-\imath s_{11}\Big)\\
  \times\Gamma\Big(\imath(\g_1+\g_3+s_{11}-p_2)+\frac{1}{2}\Big)\,
  \Gamma\Big(-\imath(\g_1-\g_2+s_{11})+\frac{1}{2}\Big)
 \ee
Thus the above integral can be identified with the r.h.s. in
\eqref{ID} via the First Barnes Identity \eqref{Barnes} with the
parameters $a_i,b_j$ given by
 \be
  a_1=\frac{1}{2}\,,\qquad a_2=\imath(\g_2-\g_1)+\frac{1}{2}\,,\\
  b_1=\imath(\g_1-p_1)+\frac{1}{2}\,,\quad
  b_2=\imath(\g_1+\g_3-p_2)+\frac{1}{2}\,.
 \ee
This completes our proof. $\Box$


\section{Appendix A: Modified Gauss-Givental realization }

In this Appendix we provide details on the construction of the
modified Gauss-Givental realization of
$\CU(\mathfrak{gl}_{\ell+1})$-representations and the corresponding
integral representation of the
$\mathfrak{gl}_{\ell+1}(\IR)$-Whittaker function. This realization
may be obtained by application of the Fourier transform to the
Gauss-Givental realization defined in \eqref{GGrep}. Recall that the
Fourier transform in the Schwartz space of functions in variables
$\xi\in\IR^d$ is given by
 \be\label{Fourier}
  \CF\,:\quad\CS(\IR^d_{\xi})\,\longrightarrow\,
  \CS(\IR^d_{\eta})\,,\\
  \hat{\phi}(\eta)\,:=\,(\CF\cdot\phi)(\eta)\,
  =\,\frac{1}{(2\pi)^{\frac{d}{2}}}\!
  \int\limits_{\IR^d}\!\prod_{a=1}^dd\xi_a\,
  e^{-\imath\,\eta_a\xi_a}\,\phi(\xi)\,.
 \ee
For the standard $L^2$-pairing $\<\,,\,\>$ in $\CS(\IR^d)$\,, the
Parseval identity holds:
 \be\label{Parseval}
  \<\phi_1,\,\phi_2\>\,
  =\,\<\hat{\phi}_1,\,\hat{\phi}_2\>\,.
 \ee
Now we apply the Fourier transform \eqref{Fourier} for
$d=\frac{\ell(\ell+1)}{2}$, to the Schwartz space of functions in
variables $T_{k,i},\,1\leq i\leq k\leq\ell$. The Fourier transforms
$\hat{\CE}_{ij}$ of the operators $\CE_{ij}$ given by \eqref{GGrep}
are defined as follows
 \be\label{FT}
  \hat{\CE}_{ij}\cdot(\CF\cdot\phi)(s)\,
  =\,\CF\bigl(\CE_{ij}\cdot\phi\bigr)(s)\,,\qquad1\leq
  i,j\leq\ell+1\,,
  \ee
for any $\phi\in\CS(\IR^{\frac{\ell(\ell+1)}{2}})$. This gives the
following  realization of $\CU(\gl_{\ell+1})$:
 \be\label{GGrepF0}
  \hat{\CE}_{i,i}\,
  =\,-\imath\Big(\g_i\,-\,\sum_{k=i}^{\ell}s_{k,i}\,
  +\,\sum_{k=1}^{i-1}s_{\ell+1+k-i,\,k}\Big)\,,\qquad1\leq
  i\leq\ell+1\,,\\
  \hat{\CE}_{i,\,i+1}\,
  =\,\sum_{n=1}^i\Big(\imath s_{\ell-i+n,n}\\
  +\imath\sum_{k=1}^{n-1}(s_{\ell-i+n-k,k}-s_{\ell+1-i+n-k,k})+\frac{1}{2}\Big)\,
  e^{\imath\pr_{s_{\ell+1-i+n,n}}-\imath\pr_{s_{\ell-i+n,n}}}\,,\\
  \hat{\CE}_{i+1,\,i}\,
  =\sum_{n=i}^{\ell}\Big[\imath(\g_{i+1}-\g_i)-\imath s_{ii}
  +\imath\sum_{k=i+1}^n(s_{k,i+1}-s_{k,i})+\frac{1}{2}\Big]\,
  e^{\imath\pr_{s_{n,i}}-\imath\pr_{s_{n+1,i+1}}}\,,
 \ee
where the operators $\CE_{i,i+1}$ and $\CE_{i+1,i}$ are defined for
$1\leq i\leq\ell$. In \eqref{GGrepF0}, we impose $s_{k,j}=0,\,k<j$
and the terms $\pr_{s_{\ell+1,k}}$ are omitted when occur. The
explicit expressions for the Whittaker vectors defined by
\eqref{whit} in this realizations might be found applying Fourier
transform to the corresponding expressions
\eqref{GGwhitL},\eqref{GGwhitR} in the Gauss-Givental realization.
Recall that in the Gauss-Givental realization the Whittaker vectors
are given by
 \be\label{GGwhitRL}
  \phi_L(T)\,
  =\,\prod_{k=1}^{\ell}e^{[\imath(\bar{\g}_{k+1}-\bar{\g}_k)-\frac{1}{2}]
    \sum\limits_{i=1}^kT_{k,i}\,
  -\sum\limits_{i=1}^ke^{T_{k+1,i}-T_{k,i}}}\,,\\
  \phi_R(T)\,
  =\,\prod_{k=1}^{\ell}\prod_{i=1}^k
  e^{\frac{T_{k,i}}{2}\,-\,e^{T_{k,i}-T_{k+1,\,i+1}}}\,,
 \ee
where we impose $T_{\ell+1,\,i}=0,\,1\leq i\leq\ell+1$ in the r.h.s.
The Fourier transform of \eqref{GGwhitRL} is easily computed by
applying the Euler integral \eqref{Euler}, assuming that the
imaginary parts of the representation parameters $\g\in\IC^{\ell+1}$
are small enough:
 \be
  \hat{\phi}_L(s)\,
  =\,\frac{1}{(2\pi)^{\frac{\ell(\ell+1)}{4}}}
  \prod_{k=1}^{\ell}\prod_{i=1}^{\ell+1-k}
  \Gamma\Big(\imath(\bar{\g}_i-\bar{\g}_{i+k})\,+\,
  \imath\sum_{n=0}^{k-1}s_{i+k-n-1,\,i}\,+\,\frac{k}{2}\Big)\,,\\
  \hat{\phi}_R(s)\,
  =\,\frac{1}{(2\pi)^{\frac{\ell(\ell+1)}{4}}}
  \prod_{k=1}^{\ell}\prod_{i=k}^{\ell}
  \Gamma\Big(-\imath\sum_{n=0}^{k-1}s_{i-n,\,k-n}\,+\,\frac{k}{2}\Big)\,.
 \ee
Thus the modified Gauss-Givental matrix element representation
\eqref{GGmatel} for the $\gl_{\ell+1}(\IR)$-Whittaker function reads
as follows:
 \be
  \Psi_{\g_1,\ldots,\g_{\ell+1}}
  (e^{x_1},\ldots,e^{x_{\ell+1}})\,
  =\,e^{-\rho(x)}\<\hat{\phi}_L,\,
  e^{-\sum\limits_{i=1}^{\ell+1}x_i\hat{\CE}_{ii}}\hat{\phi}_R\>\\
  =\,\frac{e^{\sum\limits_{n=1}^{\ell+1}\imath\g_nx_n\,-\,\rho(x)}}
  {(2\pi)^{\frac{\ell(\ell+1)}{2}}}\!
  \int\limits_{\IR^{\frac{\ell(\ell+1)}{2}}}\prod_{k=1}^{\ell}\prod_{i=1}^kds_{k,i}\,
  e^{\imath\sum\limits_{n=1}^{\ell+1}\bigl(\sum\limits_{k=n}^{\ell}s_{k,n}\,
  -\sum\limits_{k=1}^{n-1}s_{\ell+1+k-n,k}\bigr)\,x_n}\\
  \times\prod_{k=1}^{\ell}\prod_{i=k}^{\ell}
  \Gamma\Big(-\imath\sum_{n=0}^{k-1}s_{i-n,\,k-n}\,+\,\frac{k}{2}\Big)\\
  \times\prod_{k=1}^{\ell}
  \prod_{i=1}^{\ell+1-k}
  \Gamma\Big(-\imath(\g_i-\g_{i+k})\,
  -\,\imath\sum_{n=0}^{k-1}s_{i+k-n-1,\,i}\,+\,\frac{k}{2}\Big)\,.
 \ee

\section{Appendix B: Intertwining relations for $\gl_3(\IR)$}

In this Appendix, we verify the intertwining relations
\eqref{gl3intert} and \eqref{gl3intert1} satisfied by the operators
$\CB_{R}$ and $\CB^{\dag}_L$. The intertwining relations
\eqref{gl3intert} for $\CB_L$ follow from \eqref{gl3intert1}.

\subsection{Verifying the relations
$\KE_{ij}\circ\CB_R=\CB_R\circ\hat{\CE}_{ij}$}\label{BRintert}

Given the operator $\CB_R$ defined in \eqref{gl3BR}, the
intertwining relations take the form
 \be
  \KE_{ij}\cdot\!\int\limits_{\IR^3}\!\!ds\,B_R(\tau;\,s)\,\phi(s)\,
  =\!\int\limits_{\IR^3}\!\!ds\,B_R(\tau;\,s)\,(\hat{\CE}_{ij}\cdot\phi)(s)\,,\quad
  1\leq i,j\leq3\,,
  \ee
for any $\phi\in\CV_{-\imath\g-\rho}$. Using the shift of the
argument $\tau_{11}\to\tau_{11}+\imath\vk$, we obtain
 \be\label{gl3intertw0}
  \wt{E}_{ij}\cdot\!\int\limits_{\IR^3}\!\!ds\,
  B_R(\tau_{11}+\imath\vk,\tau_{21},\tau_{22};\,s)\,\phi(s)\\
  =\!\int\limits_{\IR^3}\!\!ds\,B_R(\tau_{11}+\imath\vk,\tau_{21},\tau_{22};\,s)\,
  (\hat{\CE}_{ij}\cdot\phi)(s)\,.
  \ee
Next, we rewrite \eqref{gl3intertw0} in the form of difference
equations on the kernel function $B_R(\tau;s)$.

\begin{lem}\label{GGgl3right}
The intertwining relations \eqref{gl3intertw0} are equivalent to the
following identities:
 \be\label{gl3intertw}
  \wt{E}_{ij}\cdot B_R\bigl(\tau_{11}+\imath\vk,\tau_{21},\tau_{22};\,s\bigr)\,
  =\,\hat{\CE}'_{ij}\cdot B_R\bigl(\tau_{11}+\imath\vk,\tau_{21},\tau_{22};\,s\bigr)\,,
 \ee
for $1\leq i,j\leq3$, where the operators $\hat{\CE}'_{ij}$ in
\eqref{gl3intertw} are given by
 \be\label{GGgl3Frep}
  \hat{\CE}'_{11}=-\imath(\g_1+s_{11}+s_{21}),\qquad
  \hat{\CE}'_{22}=\imath(-\g_2+s_{21}-s_{22}),\\
  \hat{\CE}'_{33}=\imath(-\g_3+s_{11}+s_{22}),\qquad
  \hat{\CE}'_{12}=\Big(\imath s_{21}-\frac{1}{2}\Big)e^{\imath\pr_{s_{21}}},\\
  \hat{\CE}'_{23}
  =\Big(\imath s_{11}-\frac{1}{2}\Big)e^{-\imath\pr_{s_{21}}+\imath\pr_{s_{11}}}\,
  +\,\Big(\imath(s_{11}+s_{22}-s_{21})-\frac{1}{2}\Big)e^{\imath\pr_{s_{22}}},\\
  \hat{\CE}'_{21}
  =\Big(\imath(\g_2-\g_1-s_{11})-\frac{1}{2}\Big)\,
  e^{-\imath\pr_{s_{11}}+\imath\pr_{s_{22}}}\\
  +\,\Big(\imath(\g_2-\g_1-s_{11}+s_{22}-s_{21})-\frac{1}{2}\Big)\,e^{-\imath\pr_{s_{21}}},\\
  \hat{\CE}'_{32}
  =\Big(\imath(\g_3-\g_2-s_{22})-\frac{1}{2}\Big)\,e^{-\imath\pr_{s_{22}}}\,.
 \ee
\end{lem}
\proof Every $\hat{\CE}_{ij},\,i\neq j$ in \eqref{GGgl3repF} is a
combination of operators of the form
$\CD^{\pm}_{\a}=\Delta^{\pm}_{\a}(s)\,e^{\pm\imath\pr_{s_{\a}}}$ and
$\CD_{\a\beta}=\Delta_{\a\beta}(s)\,e^{\imath(\pr_{s_{\a}}-\pr_{s_{\beta}})}$.
For $\CD^{\pm}_{\a}$, substitution into the r.h.s. of
\eqref{gl3intertw0}, for any $\phi\in\CV_{\imath\g}$,
 \be
   \int\limits_{\IR^3}\!\!ds\,B_R(\tau;\,s)\,(\CD^{\pm}_{\a}\phi)(s)\,
   =\!\int\limits_{\IR^3}\!\!ds_{\a}ds_{\beta}ds_{\g}\,B_R(\tau;\,s_{\a},s_{\beta},s_{\g})\,
   \Delta^{\pm}_{\a}(s)\,\phi(s_{\a}\pm\imath,s_{\beta},s_{\g})\\
   =\!\int\limits_{\IR^2}\!\!ds_{\beta}ds_{\g}\!\int\limits_{\IR\pm\imath}\!\!ds'_{\a}\,
   B_R(\tau;\,s'_{\a}\mp\imath,s_{\beta},s_{\g})\,
   \Delta^{\pm}_{\a}(s'_{\a}\mp\imath,s_{\beta},s_{\g})\,\phi(s'_{\a},s_{\beta},s_{\g})\\
   =\!\int\limits_{\IR^3}\!\!ds_{\a}ds_{\beta}ds_{\g}\,
   \bigl((\CD^{\pm}_{\a})'B_R\bigr)(\tau;\,s_{\a},s_{\beta},s_{\g})\,\phi(s_{\a},s_{\beta},s_{\g})\,,
 \ee
where
$(\CD^{\pm}_{\a})'=\Delta^{\pm}_{\a}(s'_{\a}\mp\imath,s_{\beta},s_{\g})\,e^{\mp\imath\pr_{s_{\a}}}$
providing \eqref{gl3intertw}. Similarly, for $\CD_{\a\beta}$, we
have
 \be
\hspace{-1cm}
   \int\limits_{\IR^3}\!\!ds\,B_R(\tau;\,s)\,(\CD_{\a\beta}\phi)(s)\,
   =\!\int\limits_{\IR^3}\!\!ds_{\a}ds_{\beta}ds_{\g}\,B_R(\tau;\,s_{\a},s_{\beta},s_{\g})\,
   \Delta_{\a\beta}(s)\,\phi(s_{\a}+\imath\,,
   s_{\beta}-\imath,s_{\g})\\
\hspace{-1cm}
  =\!\int\limits_{\IR}\!ds_{\g}\!\int\limits_{\IR-\imath}ds'_{\beta}\
  \!\int\limits_{\IR+\imath}\!ds'_{\a}\,
  B_R(\tau;\,s'_{\beta}+\imath,\,s'_{\a}-\imath,s_{\g})\,
  \Delta_{\a\beta}(s'_{\beta}+\imath,\,s'_{\a}-\imath,s_{\g})\,
  \phi(s'_{\a},s'_{\beta},s_{\g})\\
\hspace{-1cm}
   =\!\int\limits_{\IR^3}\!\!ds_{\a}ds_{\beta}ds_{\g}\,
   \bigl((\CD_{\a\beta})'B_R\bigr)(\tau;\,s_{\a},s_{\beta},s_{\g})\,
   \phi(s_{\a},s_{\beta},s_{\g})\,,
 \ee
where
$(\CD_{\a\beta})'=\Delta_{\a\beta}(s_{\g},s_{\beta}+\imath,\,s_{\a}-\imath)\,
  e^{\imath(-\pr_{s_{\a}}+\pr_{s_{\beta}})}$ providing \eqref{gl3intertw}. $\Box$

Note that the operators \eqref{GGgl3Frep} satisfy the ``opposite''
commutation relations:
  \be
   [\hat{\CE}'_{ij}\,,\hat{\CE}'_{k,l}]\,
   =\,-[\hat{\CE}_{ij},\,\hat{\CE}_{kl}]'\,,\qquad1\leq i,j,k,l\leq3\,.
  \ee
By \eqref{GTgl3S}, we have
 \be\label{GTrepS}
  E_{ij}=\mu_1^{-1}\,\wt{E}_{ij}\,\mu_1\,,\qquad
  \mu_1(\tau)=\frac{1}{(2\pi)^{3/2}}\frac{e^{\frac{\pi}{2}(\tau_{21}+\tau_{22})}}
  {\Gamma\bigl(\imath\tau_{21}-\imath\tau_{22}\bigr)}\,.
 \ee
After substituting \eqref{GTrepS} into the intertwining relation
\eqref{gl3intertw}, the latter takes the form
 \be\label{gl3intertw1}
  E_{ij}\cdot C_R(\tau;\,s)\,
  =\,\hat{\CE}'_{ij}\cdot C_R(\tau;\,s)\,,\qquad
  1\leq i,j\leq3\,,
 \ee
where taking into account the expression \eqref{gl3BR} we obtain
 \be\label{gl3CR}
  C_R(\tau;\,s)\,:
  =\,\mu_1^{-1}(\tau)\,
  B_R\bigl(\tau_{11}+\imath\vk,\tau_{21},\tau_{22};\,s\bigr)\\
  =\,\sqrt{2\pi}\,e^{\frac{\pi}{2}(\tau_{11}-\tau_{21}-\tau_{22})}\,
  \frac{\Gamma\bigl(\imath(\g_2-\g_1-s_{11})+\frac{1}{2}\bigr)}
  {\Gamma(-\imath s_{21}+\frac{1}{2})}\\
  \times\prod_{j=1}^2\Gamma\Big(\imath(\tau_{11}-\tau_{2,j}-s_{21})\Big)\,
  \Gamma\Big(\imath(\tau_{2,j}-\tau_{11})+\frac{1}{2}\Big)\\
  \times\prod_{j=1}^2\Gamma\Big(\imath(\g_3-\tau_{2,j})+\frac{1}{2}\Big)\\
  \times\,\delta(\g_1+s_{11}+s_{21}-\tau_{11})\,
  \delta\bigl(\g_1+\g_2+s_{11}+s_{22}-\tau_{21}-\tau_{22}\bigr)\,.
 \ee
Below we  verify the intertwining relations \eqref{gl3intertw1} by
direct calculations.

At first, for the Cartan generators, the intertwining property
\eqref{gl3intertw1} reads:
 \be
  -\imath\tau_{11}\,C_R(\tau;\,s)\,
  =\,-\imath(\g_1+s_{11}+s_{21})\,C_R(\tau,\,s)\,,\\
  -\imath(\tau_{21}+\tau_{22}-\tau_{11})\,C_R(\tau,\,s)\,
  =\,\imath(-\g_2+s_{21}-s_{22})\,C_R(\tau,\,s)\,,\\
  -\imath(\g_1+\g_2+\g_3-\tau_{21}-\tau_{22})\,C_R(\tau,\,s)\,
  =\,\imath(-\g_3+s_{11}+s_{22})\,C_R(\tau;\,s)\,,
 \ee
so that all the three relations above hold due to the presence of
delta-factors in $C_R(\tau;s)$:
$$
 \delta(\g_1+s_{11}+s_{21}-\tau_{11}+\imath\ve)\,
  \delta\bigl(\g_1+\g_2+s_{11}+s_{22}-\tau_{21}-\tau_{22}\bigr)\,.
$$

The intertwining relation
$E_{12}\,C_R(\tau;s)=\hat{\CE}'_{12}\,C_R(\tau;s)$ holds due to
 \be
  E_{12}\cdot C_R(\tau;s)\,
  =\,-\imath\Big(\imath\tau_{11}-\imath\tau_{21}+\frac{1}{2}\Big)
  \Big(\imath\tau_{11}-\imath\tau_{22}+\frac{1}{2}\Big)\,
  C_R(\tau_{11}-\imath,\tau_{21},\tau_{22}\,;s)\\
  =\,-\imath\,
  e^{\frac{\pi}{2}(\tau_{11}-\imath-\tau_{21}-\tau_{22})}\,
  \Big(\imath\tau_{11}-\imath\tau_{21}+\frac{1}{2}\Big)
  \Big(\imath\tau_{11}-\imath\tau_{22}+\frac{1}{2}\Big)
  \\
  \times\prod_{j=1}^2\frac{\imath(\tau_{11}-\tau_{2,j}-s_{21})\,
  \Gamma\Big(\imath(\tau_{11}-\tau_{2,j}-s_{21})\Big)\,
  \Gamma\Big(\imath(\tau_{2,j}-\tau_{11})+\frac{1}{2}\Big)}
  {\imath(\tau_{2,j}-\tau_{11})-\frac{1}{2}}\\
  \times\,\sqrt{2\pi}\,\frac{\Gamma\bigl(\imath(\g_2-\g_1-s_{11})+\frac{1}{2}\bigr)}
  {\Gamma(-\imath s_{21}+\frac{1}{2})}
  \prod_{j=1}^2\Gamma\Big(\imath(\g_3-\tau_{2,j})+\frac{1}{2}\Big)\\
  \times\,\delta(\g_1+s_{11}+s_{21}-\tau_{11}+\imath)\,
  \delta\bigl(\g_1+\g_2+s_{11}+s_{22}-\tau_{21}-\tau_{22}\bigr)
 \ee
 \be
  =\,-e^{\frac{\pi}{2}(\tau_{11}-\tau_{21}-\tau_{22})}\,
  \frac{\Gamma\bigl(\imath(\g_2-\g_1-s_{11})+\frac{1}{2}\bigr)}
  {\Gamma(-\imath s_{21}+\frac{1}{2})}
  \prod_{j=1}^2\Big(\imath(\tau_{11}-\tau_{2,j}-s_{21})\Big)\\
  \times\prod_{j=1}^2\Gamma\Big(\imath(\tau_{11}-\tau_{2,j}-s_{21})\Big)\,
  \Gamma\Big(\imath(\tau_{2,j}-\tau_{11})+\frac{1}{2}\Big)\,
  \Gamma\Big(\imath(\g_3-\tau_{2,j})+\frac{1}{2}\Big)\\
  \times\,\sqrt{2\pi}\,\delta(\g_1+s_{11}+s_{21}-\tau_{11}+\imath)\,
  \delta\bigl(\g_1+\g_2+s_{11}+s_{22}-\tau_{21}-\tau_{22}\bigr)\,,
 \ee
and
 \be
  \CE'_{12}\cdot C_R(\tau;s)\,
  =\,\Big(\imath s_{21}-\frac{1}{2}\Big)\,C_R(\tau\,;s_{11},s_{21}+\imath,s_{22})\\
  =\,\Big(\imath s_{21}-\frac{1}{2}\Big)\,
  e^{\frac{\pi}{2}(\tau_{11}-\tau_{21}-\tau_{22})}\,
  \frac{\Gamma\bigl(\imath(\g_2-\g_1-s_{11})+\frac{1}{2}\bigr)}
  {(-\imath s_{21}+\frac{1}{2})\,\Gamma(-\imath s_{21}+\frac{1}{2})}\\
  \times\prod_{j=1}^2\imath(\tau_{11}-\tau_{2,j}-s_{21})\,
  \Gamma\Big(\imath(\tau_{11}-\tau_{2,j}-s_{21})\Big)\,
  \Gamma\Big(\imath(\tau_{2,j}-\tau_{11})+\frac{1}{2}\Big)\\
  \times\,\sqrt{2\pi}\prod_{j=1}^2\Gamma\Big(\imath(\g_3-\tau_{2,j})+\frac{1}{2}\Big)\\
  \times\,\delta(\g_1+s_{11}+s_{21}-\tau_{11}+\imath)\,
  \delta\bigl(\g_1+\g_2+s_{11}+s_{22}-\tau_{21}-\tau_{22}\bigr)\,.
 \ee
 \be
  =\,-e^{\frac{\pi}{2}(\tau_{11}-\tau_{21}-\tau_{22})}\,
  \frac{\Gamma\bigl(\imath(\g_2-\g_1-s_{11})+\frac{1}{2}\bigr)}
  {\Gamma(-\imath s_{21}+\frac{1}{2})}
  \prod_{j=1}^2\Big(\imath(\tau_{11}-\tau_{2,j}-s_{21})\Big)\\
  \times\prod_{j=1}^2\Gamma\Big(\imath(\tau_{11}-\tau_{2,j}-s_{21})\Big)\,
  \Gamma\Big(\imath(\tau_{2,j}-\tau_{11})+\frac{1}{2}\Big)\\
  \times\,\sqrt{2\pi}\prod_{j=1}^2\Gamma\Big(\imath(\g_3-\tau_{2,j})+\frac{1}{2}\Big)\\
  \times\,\delta(\g_1+s_{11}+s_{21}-\tau_{11}+\imath)\,
  \delta\bigl(\g_1+\g_2+s_{11}+s_{22}-\tau_{21}-\tau_{22}\bigr)\,.
 \ee


The intertwining relation
$E_{21}\,C_R(\tau;s)=\hat{\CE}'_{21}\,C_R(\tau;s)$ may be verified
as follows. On the one hand, we have
 \be\label{E21C}
  E_{21}C_R(\tau;s)
  =-\imath\,C_R(\tau_{11}+\imath,\tau_{21},\tau_{22}\,;s)\,
  =\,\sqrt{2\pi}\,e^{\frac{\pi}{2}(\tau_{11}-\tau_{21}-\tau_{22})}\\
  \times\prod_{j=1}^2\frac{\imath(\tau_{2,j}-\tau_{11})+\frac{1}{2}}
  {\imath(\tau_{11}-\tau_{2,j}-s_{21})-1}
  \Gamma\Big(\imath(\tau_{11}-\tau_{2,j}-s_{21})\Big)\,
  \Gamma\Big(\imath(\tau_{2,j}-\tau_{11})+\frac{1}{2}\Big)\\
  \times\,
  \frac{\Gamma\bigl(\imath(\g_2-\g_1-s_{11})+\frac{1}{2}\bigr)}
  {\Gamma(-\imath s_{21}+\frac{1}{2})}\prod_{j=1}^2\Gamma\Big(\imath(\g_3-\tau_{2,j})+\frac{1}{2}\Big)\\
  \times\,\delta(\g_1+s_{11}+s_{21}-\tau_{11}-\imath)\,
  \delta\bigl(\g_1+\g_2+s_{11}+s_{22}-\tau_{21}-\tau_{22}\bigr)\,,
 \ee
and on the other hand we obtain
 \be
  \hat{\CE}'_{21}\cdot C_R(\tau;s)\,
  =\,\Big(\imath(\g_2-\g_1-s_{11})-\frac{1}{2}\Big)\,
  C_R(\tau;\,s_{11}-\imath,s_{21},s_{22}+\imath)\\
  +\,\Big(\imath(\g_2-\g_1-s_{11}+s_{22}-s_{21})-\frac{1}{2}\Big)\,
  C_R(\tau\,;s_{11},s_{21}-\imath,s_{22})
 \ee
 \be
  =\,\Big\{\Big(\imath(\g_2-\g_1-s_{11})-\frac{1}{2}\Big)\,
  \frac{e^{\frac{\pi}{2}(\tau_{11}-\tau_{21}-\tau_{22})}\,
  \Gamma\bigl(\imath(\g_2-\g_1-s_{11})+\frac{1}{2}\bigr)}
  {\bigl(\imath(\g_2-\g_1-s_{11})-\frac{1}{2}\bigr)\,\Gamma(-\imath s_{21}+\frac{1}{2})}\\
  \times\prod_{j=1}^2\Gamma\Big(\imath(\tau_{11}-\tau_{2,j}-s_{21})\Big)\,
  \Gamma\Big(\imath(\tau_{2,j}-\tau_{11})+\frac{1}{2}\Big)\\
  +\,\Big(\imath(\g_2-\g_1-s_{11}-s_{21}+s_{22})-\frac{1}{2}\Big)
  \frac{(-\imath s_{21}-\frac{1}{2})\,\Gamma\bigl(\imath(\g_2-\g_1-s_{11})+\frac{1}{2}\bigr)}
  {\Gamma(-\imath s_{21}+\frac{1}{2})}\\
  \times\,
  e^{\frac{\pi}{2}(\tau_{11}-\tau_{21}-\tau_{22})}
  \prod_{j=1}^2\frac{\Gamma\Big(\imath(\tau_{11}-\tau_{2,j}-s_{21})\Big)\,
  \Gamma\Big(\imath(\tau_{2,j}-\tau_{11})+\frac{1}{2}\Big)}
  {\imath(\tau_{11}-\tau_{2,j}-s_{21})-1}\Big\}\\
  \times\,\sqrt{2\pi}\prod_{j=1}^2\Gamma\Big(\imath(\g_3-\tau_{2,j})+\frac{1}{2}\Big)\\
  \times\,\delta(\g_1+s_{11}+s_{21}-\tau_{11}-\imath)\,
  \delta\bigl(\g_1+\g_2+s_{11}+s_{22}-\tau_{21}-\tau_{22}\bigr)
 \ee
 \be
  =\,\sqrt{2\pi}\,\Big\{-\,\Big(\imath(\g_2-\g_1-s_{11}-s_{21}+s_{22})-\frac{1}{2}\Big)
  \Big(\imath s_{21}+\frac{1}{2}\Big)\\
  +\prod\limits_{j=1}^2\bigl(\imath(\tau_{11}-\tau_{2,j}-s_{21})-1\bigr)\,
  \Big\}\,
  e^{\frac{\pi}{2}(\tau_{11}-\tau_{21}-\tau_{22})}\,
  \frac{\Gamma\bigl(\imath(\g_2-\g_1-s_{11})+\frac{1}{2}\bigr)}
  {\Gamma(-\imath s_{21}+\frac{1}{2})}\\
  \times\prod_{j=1}^2\frac{\Gamma\Big(\imath(\tau_{11}-\tau_{2,j}-s_{21})\Big)\,
  \Gamma\Big(\imath(\tau_{2,j}-\tau_{11})+\frac{1}{2}\Big)}
  {\imath(\tau_{11}-\tau_{2,j}-s_{21})-1}\,
  \Gamma\Big(\imath(\g_3-\tau_{2,j})+\frac{1}{2}\Big)\\
  \times\,\delta(\g_1+s_{11}+s_{21}-\tau_{11}-\imath)\,
  \delta\bigl(\g_1+\g_2+s_{11}+s_{22}-\tau_{21}-\tau_{22}\bigr)\,.
 \ee
The expression inside $\{\ldots\}$ takes the following form:
 \be\noindent
  -\,\Big[\imath(\g_2-\g_1-s_{11}+s_{22})-\Big(\imath s_{21}+\frac{1}{2}\Big)\Big]
  \Big(\imath s_{21}+\frac{1}{2}\Big)\\
  +\,\Big[\imath(\tau_{11}-\tau_{21})-\frac{1}{2}-\Big(\imath s_{21}+\frac{1}{2}\Big)\Big]
  \Big[\imath(\tau_{11}-\tau_{22})-\frac{1}{2}-\Big(\imath
  s_{21}+\frac{1}{2}\Big)\Big]\\
  =\,\prod\limits_{j=1}^2\Big(\bigl(\imath(\tau_{11}-\tau_{2,j}-\frac{1}{2}\Big)\,
  +\,2\Big(\imath s_{21}+\frac{1}{2}\Big)^2\\
  -\Big(\imath s_{21}+\frac{1}{2}\Big)
  \Big[2\imath(\tau_{11}+\imath)
  -\imath(\tau_{21}+\tau_{22})+1+\imath(\g_2-\g_1-s_{11}+s_{22})\Big]
 \ee
 \be\label{E21id}
  =\,\prod\limits_{j=1}^2\Big(\imath(\tau_{11}-\tau_{2,j}-\frac{1}{2}\Big)\,
  -\,\imath(\g_1+\g_2+s_{11}+s_{22}-\tau_{21}-\tau_{22})
  \Big(\imath s_{21}+\frac{1}{2}\Big)\,,
 \ee
which coincides with \eqref{E21C}.

The intertwining relation
$E_{32}\,C_R(\tau;s)=\hat{\CE}'_{32}\,C_R(\tau;s)$ holds due to
 \be
  E_{32}\cdot C_R(\tau;s)\,
  =\,\frac{\imath}{\imath(\tau_{21}-\tau_{22})}\Big\{
  \Big(\imath(\tau_{11}-\tau_{21})+\frac{1}{2}\Big)\,
  C_R(\tau_{11},\tau_{21}+\imath,\tau_{22};\,s)\\
  -\,\Big(\imath(\tau_{11}-\tau_{22})+\frac{1}{2}\Big)\,
  C_R(\tau_{11},\tau_{21},\tau_{22}+\imath;\,s)\Big\}\\
  =\,\frac{1}{\imath(\tau_{21}-\tau_{22})}\Big\{
  -\imath(\tau_{11}-s_{21}-\tau_{21})\,\Big(\imath(\g_3-\tau_{21})+\frac{1}{2}\Big)\\
  +\,\imath(\tau_{11}-s_{21}-\tau_{22})\,\Big(\imath(\g_3-\tau_{21})+\frac{1}{2}\Big)\Big\}
  e^{\frac{\pi}{2}(\tau_{11}-\tau_{21}-\tau_{22})}\\
  \times\,
  \frac{\Gamma\bigl(\imath(\g_2-\g_1-s_{11})+\frac{1}{2}\bigr)}
  {\Gamma(-\imath s_{21}+\frac{1}{2})}
  \prod_{j=1}^2\Gamma\Big(\imath(\tau_{11}-\tau_{2,j}-s_{21})\Big)\,
  \Gamma\Big(\imath(\tau_{2,j}-\tau_{11})+\frac{1}{2}\Big)\\
  \times\,\sqrt{2\pi}\prod_{j=1}^2\Gamma\Big(\imath(\g_3-\tau_{2,j})+\frac{1}{2}\Big)\\
  \times\,\delta(\g_1+s_{11}+s_{21}-\tau_{11})\,
  \delta\bigl(\g_1+\g_2+s_{11}+s_{22}-\tau_{21}-\tau_{22}-\imath\bigr)
 \ee
 \be\label{E32C}
  =\,\Big\{\imath(\tau_{11}-s_{21})
  +\imath\g_3+\frac{1}{2}
  -\imath(\tau_{21}+\tau_{22})\Big\}
  e^{\frac{\pi}{2}(\tau_{11}-\tau_{21}-\tau_{22})}\\
  \times\,
  \frac{\Gamma\bigl(\imath(\g_2-\g_1-s_{11})+\frac{1}{2}\bigr)}
  {\Gamma(-\imath s_{21}+\frac{1}{2})}
  \prod_{j=1}^2\Gamma\Big(\imath(\tau_{11}-\tau_{2,j}-s_{21})\Big)\,
  \Gamma\Big(\imath(\tau_{2,j}-\tau_{11})+\frac{1}{2}\Big)\\
  \times\,\sqrt{2\pi}\prod_{j=1}^2\Gamma\Big(\imath(\g_3-\tau_{2,j})+\frac{1}{2}\Big)\\
  \times\,\delta(\g_1+s_{11}+s_{21}-\tau_{11})\,
  \delta\bigl(\g_1+\g_2+s_{11}+s_{22}-\tau_{21}-\tau_{22}-\imath\bigr)\\
  =\,\Big(\imath(\g_3-\g_2-s_{22})-\frac{1}{2}\Big)\,C(\tau;\,s_{22}-\imath)\,
  =\,\hat{\CE}'_{32}\cdot C(\tau;s)\,.
 \ee


Finally we  verify whether the integral kernel \eqref{gl3CR}
intertwines the remaining generators $E_{23}$ and $\hat{\CE}_{23}$.
On the one hand, we have
 \be
  E_{23}\cdot C_R(\tau;\,s)\,
  =\,\frac{-\imath}{\imath(\tau_{21}-\tau_{22})}\Big\{
  \prod_{i=1}^3\Big(\imath(\tau_{21}-\g_i)+\frac{1}{2}\Big)\,
  C_R(\tau_{11},\tau_{21}-\imath,\tau_{22};\,s)\\
  -\,\prod_{i=1}^3\Big(\imath(\tau_{22}-\g_i)+\frac{1}{2}\Big)\,
  C_R(\tau_{11},\tau_{21},\tau_{22}-\imath;\,s)\Big\}
 \ee
 \be\label{GTE23C}
  =\,\frac{1}{\imath(\tau_{21}-\tau_{22})}\Big\{
  \frac{\imath(\tau_{21}-\tau_{11})+\frac{1}{2}}
  {\imath(\tau_{11}-\tau_{21}-s_{21})-1}
  \prod_{i=1}^2\Big(\imath(\tau_{21}-\g_i)+\frac{1}{2}\Big)\\
  -\,\frac{\imath(\tau_{22}-\tau_{11})+\frac{1}{2}}
  {\imath(\tau_{11}-\tau_{22}-s_{21})-1}
  \prod_{i=1}^2\Big(\imath(\tau_{22}-\g_i)+\frac{1}{2}\Big)\Big\}\,
  e^{\frac{\pi}{2}(\tau_{11}-\tau_{21}-\tau_{22})}\\
  \times\,
  \frac{\Gamma\bigl(\imath(\g_2-\g_1-s_{11})+\frac{1}{2}\bigr)}
  {\Gamma(-\imath s_{21}+\frac{1}{2})}
  \prod_{j=1}^2\Gamma\Big(\imath(\tau_{11}-\tau_{2,j}-s_{21})\Big)\,
  \Gamma\Big(\imath(\tau_{2,j}-\tau_{11})+\frac{1}{2}\Big)\\
  \times\,\sqrt{2\pi}\prod_{j=1}^2\Gamma\Big(\imath(\g_3-\tau_{2,j})+\frac{1}{2}\Big)\\
  \times\,\delta(\g_1+s_{11}+s_{21}-\tau_{11})\,
  \delta\bigl(\g_1+\g_2+s_{11}+s_{22}-\tau_{21}-\tau_{22}+\imath\bigr)\,.
 \ee
On the other hand,
 \be\label{GGE23C}
  \hat{\CE}'_{23}\cdot C_R(\tau;\,s)\,
  =\,\Big\{\Big(\imath s_{11}-\frac{1}{2}\Big)\,
  C_R(\tau;\,s_{11}+\imath,s_{21}-\imath,s_{22})\\
  +\Big(\imath(s_{11}-s_{21}+s_{22})-\frac{1}{2}\Big)\,
  C_R(\tau;\,s_{11},s_{21},s_{22}+\imath)\Big\}\\
  =\,\Big\{-\frac{\bigl(\imath s_{11}-\frac{1}{2}\bigr)\,
  \bigl(\imath s_{21}+\frac{1}{2}\bigr)\,\bigl(\imath(\g_2-\g_1-s_{11})+\frac{1}{2}\bigr)}
  {\bigl(\imath(\tau_{11}-\tau_{21}-s_{21})-1\bigr)\,
  \bigl(\imath(\tau_{11}-\tau_{22}-s_{21})-1\bigr)}\\
  +\,\imath(s_{11}-s_{21}+s_{22})-\frac{1}{2}\Big\}\,
  e^{\frac{\pi}{2}(\tau_{11}-\tau_{21}-\tau_{22})}\\
  \times\,
  \frac{\Gamma\bigl(\imath(\g_2-\g_1-s_{11})+\frac{1}{2}\bigr)}
  {\Gamma(-\imath s_{21}+\frac{1}{2})}
  \prod_{j=1}^2\Gamma\Big(\imath(\tau_{11}-\tau_{2,j}-s_{21})\Big)\,
  \Gamma\Big(\imath(\tau_{2,j}-\tau_{11})+\frac{1}{2}\Big)\\
  \times\,\sqrt{2\pi}\prod_{j=1}^2\Gamma\Big(\imath(\g_3-\tau_{2,j})+\frac{1}{2}\Big)\\
  \times\,\delta(\g_1+s_{11}+s_{21}-\tau_{11})\,
  \delta\bigl(\g_1+\g_2+s_{11}+s_{22}-\tau_{21}-\tau_{22}+\imath\bigr)\,.
 \ee
Let us  introduce the following variables
 \be\label{E23Cvr}
  s:=\imath s_{21}+\frac{1}{2}\,,\qquad t:=\imath\tau_{11}\,,\qquad
  a_j:=\imath\tau_{2,j}+\frac{1}{2}\,,\quad j=1,2\,.
 \ee
Then the equality of the expressions in \eqref{GTE23C} and
\eqref{GGE23C} follows via the identity
 \be\label{E23Cid}
  \frac{1}{a_1-a_2}\Big\{\frac{a_1-t}{t-s-a_1}\prod_{i=1}^2(a_1-\imath\g_i)\,
  -\,\frac{a_2-t}{t-s-a_2}\prod_{i=1}^2(a_2-\imath\g_i)\Big\}\\
  =\,a_1+a_2-\imath(\g_1+\g_2)-s\,+\,s\prod_{i=1}^2\frac{\imath\g_i-t+s}{t-s-a_i}\,.
 \ee

\subsection{Proof of Lemma \ref{BRLrels}}\label{BLintert}

Below we provide a proof of  the relations \eqref{gl3intert1}. At
first, we rewrite  the intertwining relations for the operator
$\CB^{\dag}_L$ defined in \eqref{gl3BLdag}, for any
$\phi\in\CV_{-\imath\g-\rho}$ in the following form:
 \be\label{gl3intertw2}
  \hat{\CE}_{ij}\cdot\!\int\limits_{\IR^3}\!\!d\tau\,B^{\dag}_L(s;\,\tau)\,\phi(\tau)\,
  =\!\int\limits_{\IR^3}\!\!d\tau\,B^{\dag}_L(s;\,\tau)\,
  (\KE_{ij}\cdot\phi)(\tau)\,,\quad1\leq i,j\leq3\,.
  \ee
Similarly to \eqref{gl3CR}, introduce the following integral kernel
 \be\label{gl3CL}
  C_L(s;\tau)\,
  =\,\mu_1(\tau)\times\,B^{\dag}_L(s;\,\tau_{11}+\imath\vk,\tau_{21},\tau_{22})\\
  =\,\frac{\sqrt{2\pi}\,e^{-\frac{\pi\tau_{11}}{2}}}
  {\Gamma\bigl(\imath\tau_{21}-\imath\tau_{22}\bigr)\,
  \Gamma\bigl(-\imath\tau_{21}+\imath\tau_{22}\bigr)}\,
  \frac{\Gamma\bigl(-\imath s_{21}+\frac{1}{2}\bigr)}
  {\Gamma\bigl(-\imath(\g_1-\g_2+s_{11})+\frac{1}{2}\bigr)}\\
  \times\prod_{j=1}^2\frac{\Gamma\bigl(-\imath(\tau_{11}-\tau_{2,j}-s_{21})\bigr)}
  {\Gamma\big(\imath(\tau_{2,j}-\tau_{11})+\frac{1}{2}\bigr)\,
  \Gamma\bigl(\imath(\g_3-\tau_{2,j})+\frac{1}{2}\bigr)}\\
  \times\,\delta(\g_1+s_{11}+s_{21}-\tau_{11})\,
  \delta\bigl(\g_1+\g_2+s_{11}+s_{22}-\tau_{21}-\tau_{22}\bigr)\,.
 \ee
The intertwining relations satisfied by $C_L$ to be checked are
given by the following Lemma.
\begin{lem}
The intertwining relations \eqref{gl3intertw2} are equivalent to the
following identity:
 \be\label{gl3intertwL}
  \hat{\CE}_{ij}\cdot C_L(s;\tau)\,
  =\,E'_{ij}\cdot C_L(s;\tau)\,,\qquad
  1\leq i,j\leq3\,,
 \ee
where the operators $E'_{ij}$ are given by
 \be\label{GTgl3right}
  E'_{11}=-\imath\tau_{11},\qquad
  E'_{22}=-\imath(\tau_{21}+\tau_{22})+\imath\tau_{11},\\
  E'_{33}=-\imath(\g_1+\g_2+\g_3)+\imath(\tau_{21}+\tau_{22}),\\
  E'_{12}=-\imath\Big(\imath(\tau_{11}-\tau_{21})-\frac{1}{2}\Big)
  \Big(\imath(\tau_{11}-\tau_{22})-\frac{1}{2}\Big)\,e^{\imath\pr_{\tau_{11}}},\\
  E'_{23}
  =-\imath\Big\{\frac{\prod\limits_{i=1}^3\bigl(\imath\tau_{21}-\imath\g_i-\frac{1}{2}\bigr)}
  {\imath\tau_{21}-\imath\tau_{22}-1}\,e^{\imath\pr_{\tau_{21}}}\,
  +\,\frac{\prod\limits_{j=1}^3\bigl(\imath\tau_{22}-\imath\g_i-\frac{1}{2}\bigr)}
  {\imath\tau_{22}-\imath\tau_{21}-1}\,e^{\imath\pr_{\tau_{22}}}\Big\}\,,\\
  E'_{21}=-\imath e^{-\imath\pr_{\tau_{11}}}\,,\\
  E'_{32}=\imath\Big\{\frac{\imath\tau_{11}-\imath\tau_{21}-\frac{1}{2}}
  {\imath\tau_{21}-\imath\tau_{22}+1}\,e^{-\imath\pr_{\tau_{21}}}\,
  +\,\frac{\imath\tau_{11}-\imath\tau_{22}-\frac{1}{2}}
  {\imath\tau_{22}-\imath\tau_{21}+1}\,e^{-\imath\pr_{\tau_{22}}}\Big\}\,.
 \ee
\end{lem}
\proof Similar to the proof of Lemma \eqref{GGgl3right}.  $\Box$

Similarly to \eqref{GGgl3Frep}, the operators \eqref{GTgl3right}
satisfy the ``opposite'' commutation relations:
  \be
   [E'_{ij}\,,E'_{k,l}]\,
   =\,-[E_{ij},\,E_{kl}]'\,,\qquad1\leq i,j,k,l\leq3\,.
   \ee
Verification of the  intertwining relations
\eqref{gl3intertwL} for the operator $\CC_L$ proceeds as follows.

For the Cartan generators $E_{ii}$ and $\CE_{ii}$, the
relations \eqref{gl3intertwL} hold due to the presence of the
delta-factor in $C_L(s;\,\tau)$,
$$
 \delta(\g_1+s_{11}+s_{21}-\tau_{11})\,
 \delta\bigl(\g_1+\g_2+s_{11}+s_{22}-\tau_{21}-\tau_{22}\bigr)\,.
$$
The relation $\hat{\CE}_{12}\cdot C_L(s;\tau)=E'_{12}\cdot
C_L(s;\,\tau)$ follows since
 \be
  \hat{\CE}_{12}\cdot C_L(s;\tau)\,
  =\,\Big(\imath
  s_{21}+\frac{1}{2}\Big)\,C_L(s_{11},s_{21}-\imath,s_{22};\,\tau)\\
  =\,-\frac{e^{-\frac{\pi\tau_{11}}{2}}}
  {\Gamma\bigl(\imath\tau_{21}-\imath\tau_{22}\bigr)\,
  \Gamma\bigl(-\imath\tau_{21}+\imath\tau_{22}\bigr)}\,
  \frac{\Gamma\bigl(-\imath s_{21}+\frac{1}{2}\bigr)}
  {\Gamma\bigl(-\imath(\g_1-\g_2+s_{11})+\frac{1}{2}\bigr)}\\
  \times\,\sqrt{2\pi}\prod_{j=1}^2
  \frac{-\imath(\tau_{11}-\tau_{2,j}-s_{21})\,
  \Gamma\bigl(-\imath(\tau_{11}-\tau_{2,j}-s_{21})\bigr)}
  {\Gamma\big(\imath(\tau_{2,j}-\tau_{11})+\frac{1}{2}\bigr)\,
  \Gamma\bigl(\imath(\g_3-\tau_{2,j})+\frac{1}{2}\bigr)}\\
  \times\,\delta(\g_1+s_{11}+s_{21}-\tau_{11}-\imath)\,
  \delta\bigl(\g_1+\g_2+s_{11}+s_{22}-\tau_{21}-\tau_{22}\bigr)
 \ee
equals to
 \be
  E'_{12}\cdot C_L(s;\,\tau)\,
  =\,-\imath\prod_{j=1}^2\Big(\imath(\tau_{11}-\tau_{2,j})-\frac{1}{2}\Big)\,
  C_L(s;\,\tau_{11}+\imath,\tau_{21},\tau_{22})\\
  =\,-\frac{e^{-\frac{\pi\tau_{11}}{2}}}
  {\Gamma\bigl(\imath\tau_{21}-\imath\tau_{22}\bigr)\,
  \Gamma\bigl(-\imath\tau_{21}+\imath\tau_{22}\bigr)}\,
  \frac{\Gamma\bigl(-\imath s_{21}+\frac{1}{2}\bigr)}
  {\Gamma\bigl(-\imath(\g_1-\g_2+s_{11})+\frac{1}{2}\bigr)}\\
  \times\,\sqrt{2\pi}\prod_{j=1}^2
  \frac{-\imath(\tau_{11}-\tau_{2,j}-s_{21})\,
  \Gamma\bigl(-\imath(\tau_{11}-\tau_{2,j}-s_{21})\bigr)}
  {\Gamma\big(\imath(\tau_{2,j}-\tau_{11})+\frac{1}{2}\bigr)\,
  \Gamma\bigl(\imath(\g_3-\tau_{2,j})+\frac{1}{2}\bigr)}\\
  \times\,\delta(\g_1+s_{11}+s_{21}-\tau_{11}-\imath)\,
  \delta\bigl(\g_1+\g_2+s_{11}+s_{22}-\tau_{21}-\tau_{22}\bigr)\,.
 \ee


For the relation $\hat{\CE}_{21}\cdot C_L(s;\tau)=E'_{21}\cdot
C_L(s;\,\tau)$, we have
 \be
  \hat{\CE}_{21}\cdot C_L(s;\tau)\,
  =\,\Big(\imath(\g_2-\g_1-s_{11})+\frac{1}{2}\Big)\,
  C(s_{11}+\imath,s_{21},s_{22}-\imath;\,\tau)\\
  +\,\Big(\imath(\g_2-\g_1-s_{11}-s_{21}+s_{22})+\frac{1}{2}\Big)\,
  C(s_{11},s_{21}+\imath,s_{22};\,\tau)\\
  =\,\sqrt{2\pi}\,\Big\{-\Big(\imath(\g_2-\g_1-s_{11}-s_{21}+s_{22})-\frac{1}{2}\Big)\Big(\imath s_{21}+\frac{1}{2}\Big)\\
  +\prod\limits_{j=1}^2\bigl(\imath(\tau_{11}-\tau_{2,j}-s_{21})-1\bigr)\Big\}
  \frac{e^{-\frac{\pi\tau_{11}}{2}}}
  {\Gamma\bigl(\imath\tau_{21}-\imath\tau_{22}\bigr)\,
  \Gamma\bigl(-\imath\tau_{21}+\imath\tau_{22}\bigr)}\\
  \times\frac{\Gamma\bigl(-\imath s_{21}+\frac{1}{2}\bigr)}
  {\Gamma\bigl(-\imath(\g_1-\g_2+s_{11})+\frac{1}{2}\bigr)}
  \prod_{j=1}^2
  \frac{\Gamma\bigl(-\imath(\tau_{11}-\tau_{2,j}-s_{21})\bigr)}
  {\Gamma\big(\imath(\tau_{2,j}-\tau_{11})+\frac{1}{2}\bigr)\,
  \Gamma\bigl(\imath(\g_3-\tau_{2,j})+\frac{1}{2}\bigr)}\\
  \times\,\delta(\g_1+s_{11}+s_{21}-\tau_{11}+\imath)\,
  \delta\bigl(\g_1+\g_2+s_{11}+s_{22}-\tau_{21}-\tau_{22}\bigr)
 \ee
and
 \be
  E'_{21}\cdot C_L(s;\,\tau)\,
  =\,-\imath\,C_L(s;\,\tau_{11}-\imath,\tau_{21},\tau_{22})\\
  =\,\frac{\sqrt{2\pi}\,e^{-\frac{\pi\tau_{11}}{2}}}
  {\Gamma\bigl(\imath\tau_{21}-\imath\tau_{22}\bigr)\,
  \Gamma\bigl(-\imath\tau_{21}+\imath\tau_{22}\bigr)}\,
  \frac{\Gamma\bigl(-\imath s_{21}+\frac{1}{2}\bigr)}
  {\Gamma\bigl(-\imath(\g_1-\g_2+s_{11})+\frac{1}{2}\bigr)}\\
  \times\prod_{j=1}^2
  \frac{\imath(\tau_{2,j}-\tau_{11})+\frac{1}{2}}
  {\imath(\tau_{11}-\tau_{2,j}-s_{21})-1}
  \frac{\Gamma\bigl(-\imath(\tau_{11}-\tau_{2,j}-s_{21})\bigr)}
  {\Gamma\big(\imath(\tau_{2,j}-\tau_{11})+\frac{1}{2}\bigr)\,
  \Gamma\bigl(\imath(\g_3-\tau_{2,j})+\frac{1}{2}\bigr)}\\
  \times\,\delta(\g_1+s_{11}+s_{21}-\tau_{11}+\imath)\,
  \delta\bigl(\g_1+\g_2+s_{11}+s_{22}-\tau_{21}-\tau_{22}\bigr)\,,
 \ee
so that the equality follows by \eqref{E21id}.


The relation $\hat{\CE}_{32}\cdot C_L(s;\tau)=E'_{32}\cdot
C_L(s;\,\tau)$ follows due to
 \be
  \hat{\CE}_{32}\cdot C_L(s;\tau)\,
  =\,\Big(\imath(\g_3-\g_2-s_{22})+\frac{1}{2}\Big)\,
  C_L(s_{11},s_{21},s_{22}+\imath;\,\tau)\\
  =\,\frac{e^{-\frac{\pi\tau_{11}}{2}}\,
  \bigl(\imath(\g_3-\g_2-s_{22})+\frac{1}{2}\bigr)}
  {\Gamma\bigl(\imath\tau_{21}-\imath\tau_{22}\bigr)\,
  \Gamma\bigl(-\imath\tau_{21}+\imath\tau_{22}\bigr)}\,
  \frac{\Gamma\bigl(-\imath s_{21}+\frac{1}{2}\bigr)}
  {\Gamma\bigl(-\imath(\g_1-\g_2+s_{11})+\frac{1}{2}\bigr)}\\
  \times\,\sqrt{2\pi}\prod_{j=1}^2
  \frac{\Gamma\bigl(-\imath(\tau_{11}-\tau_{2,j}-s_{21})\bigr)}
  {\Gamma\big(\imath(\tau_{2,j}-\tau_{11})+\frac{1}{2}\bigr)\,
  \Gamma\bigl(\imath(\g_3-\tau_{2,j})+\frac{1}{2}\bigr)}\\
  \times\,\delta(\g_1+s_{11}+s_{21}-\tau_{11})\,
  \delta\bigl(\g_1+\g_2+s_{11}+s_{22}-\tau_{21}-\tau_{22}+\imath\bigr)\,,
 \ee
so that
 \be
  E'_{32}\cdot C_L(s;\,\tau)\,
  =\,\imath\Big\{\frac{\imath(\tau_{11}-\tau_{21})-\frac{1}{2}}
  {\imath(\tau_{21}-\tau_{22})+1}\,
  C_L(s;\,\tau_{11},\tau_{21}-\imath,\tau_{22})\\
  -\,\frac{\imath(\tau_{11}-\tau_{22})-\frac{1}{2}}
  {\imath(\tau_{22}-\tau_{21})+1}\,
  C_L(s;\,\tau_{11},\tau_{21},\tau_{22}+\imath)\Big\}\\
  =\,\frac{1}{\imath(\tau_{21}-\tau_{22})}\Big\{
  -\imath(\tau_{11}-s_{21}-\tau_{21})\,\Big(\imath(\g_3-\tau_{21})+\frac{1}{2}\Big)\\
  +\,\imath(\tau_{11}-s_{21}-\tau_{22})\,\Big(\imath(\g_3-\tau_{21})+\frac{1}{2}\Big)\Big\}\\
  \times\,\frac{e^{-\frac{\pi\tau_{11}}{2}}}
  {\Gamma\bigl(\imath\tau_{21}-\imath\tau_{22}\bigr)\,
  \Gamma\bigl(-\imath\tau_{21}+\imath\tau_{22}\bigr)}\,
  \frac{\Gamma\bigl(-\imath s_{21}+\frac{1}{2}\bigr)}
  {\Gamma\bigl(-\imath(\g_1-\g_2+s_{11})+\frac{1}{2}\bigr)}\\
  \times\,\sqrt{2\pi}\prod_{j=1}^2
  \frac{\Gamma\bigl(-\imath(\tau_{11}-\tau_{2,j}-s_{21})\bigr)}
  {\Gamma\big(\imath(\tau_{2,j}-\tau_{11})+\frac{1}{2}\bigr)\,
  \Gamma\bigl(\imath(\g_3-\tau_{2,j})+\frac{1}{2}\bigr)}\\
  \times\,\delta(\g_1+s_{11}+s_{21}-\tau_{11})\,
  \delta\bigl(\g_1+\g_2+s_{11}+s_{22}-\tau_{21}-\tau_{22}+\imath\bigr)\\
  =\,\Big(\imath(\g_3-\g_2-s_{22})+\frac{1}{2}\Big)\,C_L(s_{11},s_{21},s_{22}+\imath;\,\tau)\,
  =\,\hat{\CE}_{32}\cdot C_L(s;\tau)\,,
 \ee
similarly to \eqref{E32C}.


As for  the relation $\hat{\CE}_{23}\cdot
C_L(s;\tau)=E'_{23}\cdot C_L(s;\,\tau)$, on the one hand, we have
 \be
  \hat{\CE}_{23}\cdot C_L(s;\tau)\,
  =\,\Big(\imath s_{11}+\frac{1}{2}\Big)\,
  C_L(s_{11}-\imath,s_{21}+\imath,s_{22};\,\tau)\\
  +\,\Big(\imath(s_{11}-s_{21}+s_{22})+\frac{1}{2}\Big)\,
  C_L(s_{11},s_{21},s_{22}-\imath;\,\tau)\\
  =\,\Big\{-\frac{\bigl(\imath s_{11}-\frac{1}{2}\bigr)\,
  \bigl(\imath s_{21}+\frac{1}{2}\bigr)\,\bigl(\imath(\g_2-\g_1-s_{11})+\frac{1}{2}\bigr)}
  {\bigl(\imath(\tau_{11}-\tau_{21}-s_{21})-1\bigr)\,
  \bigl(\imath(\tau_{11}-\tau_{22}-s_{21})-1\bigr)}\\
  +\,\imath(s_{11}-s_{21}+s_{22})-\frac{1}{2}\Big\}\,
  \frac{e^{-\frac{\pi\tau_{11}}{2}}}
  {\Gamma\bigl(\imath\tau_{21}-\imath\tau_{22}\bigr)\,
  \Gamma\bigl(-\imath\tau_{21}+\imath\tau_{22}\bigr)}\\
  \times\frac{\Gamma\bigl(-\imath s_{21}+\frac{1}{2}\bigr)}
  {\Gamma\bigl(-\imath(\g_1-\g_2+s_{11})+\frac{1}{2}\bigr)}\\
  \times\,\sqrt{2\pi}\prod_{j=1}^2\frac{\Gamma\bigl(-\imath(\tau_{11}-\tau_{2,j}-s_{21})\bigr)}
  {\Gamma\big(\imath(\tau_{2,j}-\tau_{11})+\frac{1}{2}\bigr)\,
  \Gamma\bigl(\imath(\g_3-\tau_{2,j})+\frac{1}{2}\bigr)}\\
  \times\,\delta(\g_1+s_{11}+s_{21}-\tau_{11})\,
  \delta\bigl(\g_1+\g_2+s_{11}+s_{22}-\tau_{21}-\tau_{22}-\imath\bigr)\,,
 \ee
and on the other hand,
 \be
  E'_{23}\cdot C_L(s;\,\tau)\,
  =\,-\imath\Big\{\frac{\prod\limits_{i=1}^3\bigl(\imath\tau_{21}-\imath\g_i-\frac{1}{2}\bigr)}
  {\imath\tau_{21}-\imath\tau_{22}-1}\,C_L(s;\,\tau_{11},\tau_{21}+\imath,\tau_{22})\\
  +\,\frac{\prod\limits_{j=1}^3\bigl(\imath\tau_{22}-\imath\g_i-\frac{1}{2}\bigr)}
  {\imath\tau_{22}-\imath\tau_{21}-1}\,C_L(s;\,\tau_{11},\tau_{21},\tau_{22}+\imath)\Big\}\\
  =\,\frac{1}{\imath(\tau_{21}-\tau_{22})}\Big\{
  \frac{\imath(\tau_{21}-\tau_{11})+\frac{1}{2}}
  {\imath(\tau_{11}-\tau_{21}-s_{21})-1}
  \prod_{i=1}^2\Big(\imath(\tau_{21}-\g_i)+\frac{1}{2}\Big)\\
  -\,\frac{\imath(\tau_{22}-\tau_{11})+\frac{1}{2}}
  {\imath(\tau_{11}-\tau_{22}-s_{21})-1}
  \prod_{i=1}^2\Big(\imath(\tau_{22}-\g_i)+\frac{1}{2}\Big)\Big\}\\
  \times\,\frac{e^{-\frac{\pi\tau_{11}}{2}}}
  {\Gamma\bigl(\imath\tau_{21}-\imath\tau_{22}\bigr)\,
  \Gamma\bigl(-\imath\tau_{21}+\imath\tau_{22}\bigr)}\,
  \frac{\Gamma\bigl(-\imath s_{21}+\frac{1}{2}\bigr)}
  {\Gamma\bigl(-\imath(\g_1-\g_2+s_{11})+\frac{1}{2}\bigr)}\\
  \times\,\sqrt{2\pi}\prod_{j=1}^2\frac{\Gamma\bigl(-\imath(\tau_{11}-\tau_{2,j}-s_{21})\bigr)}
  {\Gamma\big(\imath(\tau_{2,j}-\tau_{11})+\frac{1}{2}\bigr)\,
  \Gamma\bigl(\imath(\g_3-\tau_{2,j})+\frac{1}{2}\bigr)}\\
  \times\,\delta(\g_1+s_{11}+s_{21}-\tau_{11})\,
  \delta\bigl(\g_1+\g_2+s_{11}+s_{22}-\tau_{21}-\tau_{22}-\imath\bigr)\,.
 \ee
Then applying the substitution \eqref{E23Cvr} we deduce the required
equality via \eqref{E23Cid}.


\noindent {\small {\bf A.A.G.} {\sl Laboratory for Quantum Field
Theory
and Information},\\
\hphantom{xxxx} {\sl Institute for Information
Transmission Problems, RAS, 127994, Moscow, Russia};\\
\hphantom{xxxx} {\it E-mail address}: {\tt anton.a.gerasimov@gmail.com}}\\
\noindent{\small {\bf D.R.L.} {\sl Laboratory for Quantum Field
Theory
and Information},\\
\hphantom{xxxx}  {\sl Institute for Information
Transmission Problems, RAS, 127994, Moscow, Russia};\\
\hphantom{xxxx} {\sl Moscow Center for Continuous Mathematical
Education,\\
\hphantom{xxxx} 119002,  Bol. Vlasyevsky per. 11, Moscow, Russia};\\
\hphantom{xxxx} {\it E-mail address}: {\tt lebedev.dm@gmail.com}}\\
\noindent{\small {\bf S.V.O.} {\sl
 Beijing Institute of Mathematical Sciences and Applications\,,\\
\hphantom{xxxx} \textnumero544, Hefangkou Village, Huaibei Town, Huairou District, Beijing 101408, China};\\
\hphantom{xxxx} {\it E-mail address}: {\tt oblezin@gmail.com}}

\begin{thebibliography}{120}
\frenchspacing \smallbreak
\bibitem[BE]{BE} Bateman Manuscript Project, ed. A.~Erdelyi,
Higher Transcendental Functions, vol. 1, McGraw-Hill, 1953.

\bibitem[Ba]{Ba} E.W.~Barnes, {\it A new development of the theory of hypergeometric functions},
Proc. London Math. Soc. (2) 6 (1908) 141-177.

\bibitem[GG]{GG} I.M. Gelfand, M.I. Graev,
  {\it GG-functions and their relation with general hypergeometric functions},
  Rus. Math. Surveys 52:4 (1997) 3--48.

\bibitem[GS]{GS} I.M. Gelfand, G.E. Shilov,
  Generalized Functions, vol. 1: Properties and Operations,
  Academic Press, 2014.


\bibitem[GT]{GT} I.M. Gelfand, M.L. Tsetlin,
  {\it Finite-dimensional representations of the group of unimodular matrices},
  Dokl. Akad. Nauk SSSR 71 (1950) 825-828 (Russian). English
  transl. in I.M.Gelfand, Collected papers, vol. II, Springer 1988,
  653-656.

\bibitem[GKL]{GKL} A. Gerasimov, S. Kharchev, D. Lebedev,
{\it Representation theory and quantum inverse scattering method:
open Toda chain and hyperbolic Sutherland model}, Int. Math. Res.
Notices 17 (2004) 823-854; {\tt[arXiv:math.QA/0204206]}.

\bibitem[GKLO]{GKLO} A. Gerasimov, S.Kharchev, D. Lebedev, S. Oblezin,
{\it On a Gauss-Givental of quantum Toda chain wave function}, Int.
Math. Res. Notices (2006), Art. ID 96489;
{\tt[arXiv:math.RT/0505310]}.

\bibitem[GLO08]{GLO08} A.~Gerasimov, D.~Lebedev, S.~Oblezin, {\it Baxter operator
and Archimedean Hecke algebra}, Commun. Math. Phys. 284 (2008)
867--896; {\tt[arXiv:0706.3476]}.

\bibitem[GLO11]{GLO11} A. Gerasimov, D. Lebedev, S. Oblezin,
  {\it Parabolic Whittaker functions and topological field theories I},
Commun. Number Theory Phys. 5:1 (2011) 135-201; [{\tt
arXiv:1002.2622}].

\bibitem[GLO23]{GLO23} A. Gerasimov, D. Lebedev, S. Oblezin,
  {\it On matix element representation for the GKZ hypergeometric functions},
  Lett. Math. Phys. 113:43 (2023); [{\tt arXiv:2209.2886}].

\bibitem[Giv]{Giv} A.~Givental,
{\it Stationary Phase Integrals, Quantum Toda Lattices, Flag
Manifolds and the Mirror Conjecture}. Topics in Singularity Theory,
AMS Transl. Ser. 2, 180, AMS, Providence RI 1997, 103--115;
{\tt[arXiv:9612001]}.

\bibitem [Gu]{Gu} R. A. Gustafson, {\it Some $q$-Beta and Mellin-Barnes integrals on complact
Lie groups and Lie algebras}, Trans. AMS 341:1 (1994) 69--119.



\bibitem[JK]{JK} D.~Joe, B.~Kim,
{\it Equivariant mirrors and the Virasoro conjecture
for flag manifolds}, Int. Math. Res. Notices {\bf 2003} No. 15
(2003), 859--882; {\tt[arXiv: math.AG/0210377]}.

\bibitem[KL]{KL} S.~Kharchev, D.~Lebedev,
{\it Eigenfunctions of $GL(N,R)$ Toda chain:  The Mellin-Barnes
representation}, JETP Lett. 71, (2000), 235-238,
{\tt[arXiv:hep-th/0004065]}.


\bibitem[Ko]{Ko} B. Kostant, {\it On Whittaker vectors and representation theory},
Invent. Math. 48:2 (1978) 101--184.

\bibitem[Ra]{Ra} S.~Ramanujan, {\it A class of definite integrals},
Quart. J. Math. 48 (1920) 294-310.

\bibitem[STS]{STS} M. Semenov-Tian-Shansky,
{\it Quantization of open Toda lattices},
Chapter III in Encyclopedia of Mathematics 1987, 226-258.

\bibitem[Skl]{Skl} E.K.~Sklyanin,
{\it Separation of variables: new trends}, in Quantum field theory,
integrable models and beyond, Kyoto 1994, Progr. Theor. Phys. Suppl.
118 (1995) 35-60; {\tt[arXiv:solv-int/9504001]}.


\bibitem[SYZ]{SYZ} A.Strominger, S.-T.Yau, E.Zaslow, {\it Mirror symmetry is
T-duality}, Nucl. Phys. B479(1-2) (1996) 243-259;
{\tt[arXiv:hep-th/9606040]}.

\bibitem[V]{V} N.Ja.Vilenkin, Special functions and the theory of group representations, AMS 1968.

\end{thebibliography}
\end{document}